\documentclass[12pt,reqno,a4paper]{amsart}
\usepackage[english]{babel}

\usepackage{amsmath,amsthm,amsfonts,amssymb}
\usepackage{mathtools}
\usepackage{color}
\usepackage{lineno}
\usepackage{mathrsfs}
\usepackage[all]{xy}
\usepackage{tikz-cd}
\usepackage[normalem]{ulem} 
\usepackage{datetime}

\usepackage{hyperref}

\theoremstyle{definition}
\newtheorem{theorem}{Theorem}[section]
\newtheorem{definition}[theorem]{Definition}
\newtheorem{lemma}[theorem]{Lemma}
\newtheorem{proposition}[theorem]{Proposition}
\newtheorem{corollary}[theorem]{Corollary}
\newtheorem{remark}[theorem]{Remark}
\newtheorem{example}[theorem]{Example}
\newtheorem{construction}{Construction}[section]
\newcommand{\mc}{\mathcal}

\newcommand{\xra}{\xrightarrow}
\newcommand{\ra}{\rightarrow}
\newcommand{\rra}{\rightrightarrows}
\calclayout

\title{On two notions of a Gerbe over a stack}
\author{Saikat Chatterjee, Praphulla Koushik } 
\address{School of Mathematics, Indian Institute of Science Education and Research-Thiruvananthapuram \\
	Maruthamala PO, Vithura\\ Kerala-695551\\
	India}

\email{saikat.chat01@gmail.com\\
	koushik16@iisertvm.ac.in}

\begin{document}
	\begin{abstract}
		Let $\mc{G}$ be a Lie groupoid. The category $B\mc{G}$ of principal $\mc{G}$-bundles defines a differentiable stack. On the other hand, given a differentiable stack $\mc{D}$, there exists a Lie groupoid $\mc{H}$ such that $B\mc{H}$ is isomorphic to $\mc{D}$. Define a gerbe over a stack as a morphism of stacks $F\colon \mc{D}\ra \mc{C}$, such that $F$ and the diagonal map $\Delta_F\colon \mc{D}\ra \mc{D}\times_{\mc{C}}\mc{D}$ are epimorphisms. This paper explores the relationship between a gerbe defined above and a Morita equivalence class of a Lie groupoid extension. 
	\end{abstract} 
	
	\maketitle 
	
	\section{Introduction} 
	The paper was inspired by two different notions of \textit{a differentiable gerbe over a differentiable stack} that we have encountered. One notion is, as a morphism of stacks satisfying some ``additional" properties, defined by Behrend and Xu in \cite{MR2817778} and some others \cite{Ginot, MR2206877}. The second notion is, as a Morita equivalence class of a Lie groupoid extension, defined by Laurent-Gengoux, Stienon, and Xu in \cite{MR2493616}.
	Indeed, it is well known that given a Lie groupoid $\mc{G}$, the category of principal $\mc{G}$-bundles, denoted by $B\mc{G}$, is a differentiable stack \cite{MR2493616}. On the other hand, given a differentiable stack $\mc{D}$, there exists a Lie groupoid $\mc{H}$ such that $\mc{D}$ is isomorphic to $B\mc{H}$ \cite{MR2778793}. This suggests a possibility of correspondence between the notions of a gerbe over a stack mentioned above. However, to the best of our knowledge, nowhere, this possible correspondence between the two definitions has been explored. In this paper, we investigate this correspondence. To be precise, our observations are the following. 
	\begin{itemize}
		\item Let $F\colon \mc{D}\ra \mc{C}$ be a morphism of stacks, such that $F$ is an `epimorphism' and the diagonal morphism $\Delta_F\colon \mc{D}\ra \mc{D}\times_{\mc{C}}\mc{D}$ is a `representable surjective submersion'. Then there exists a Morita equivalence class of Lie groupoid extension (Theorem \ref{Theorem:GoSgivesLiegroupoidextension}).
		\item Let $f \colon \mc{G}\to \mc{H}$ be a Lie groupoid extension. Then there exists an epimorphism of stacks $F\colon B\mc{G}\to B\mc{H}$, such that the corresponding diagonal map $\Delta_F\colon B\mc G\ra B\mc G \times_{B\mc H} B\mc G$ is an epimorphism (Theorem \ref{Theorem:BGBHisaGoS}).
	\end{itemize}
	
	Since the introduction of gerbes by Giraud to study the nonabelian cohomology \cite{MR0344253}, the subject has grown rapidly in various directions and various forms. Nonabelian gerbes (or co-cycle gerbes) over a manifold and connection structures on it have appeared in several articles in the last few years, such as \cite{{MR3415507},{MR2764890},{MR2520993},{MR3084724},{MR3357822}}, just to mention a few. The cocycle gerbe and the connection structure on it is closely related to the so-called higher principal bundles and higher gauge theory, which is also quite well studied \cite{{MR3357822}, {MR3566125}, {MR3415507}, {MR3504595}, {MR2661492}, {MR2764890}, {MR2146290}, {MR2342821}, {MR3294946}} etc. A more algebraic geometry flavored treatment of the connection structure on the nonabelian gerbes can be found in the work of Breen and Messing \cite{MR2183393}. The same paper also describes the stack of gauge transformations of a differential gerbe. The abelian bundle gerbe was introduced by Murray \cite{MR2681698}, and subsequently, a nonabelian version of the same was proposed by Aschieri et al. \cite{ACJ2005}. The paper \cite{MR3089401} discusses an equivalence between various related notions of a nonabelian gerbe over a manifold. As mentioned earlier, in this paper, we will be concerned with the notion of a gerbe over a stack. The central idea for our constructions is the correspondence between Lie groupoids and differentiable stacks. We were aided by several articles available on the topic, such as \cite{{MR2817778}, {MR2493616}, {Moerdijk2}, {Matias}, {Ginot}, {Zhu}, {MR3480061}}. For a general exposition on stacks and differentiable stacks we refer to \cite{{Metzler}, {MR2206877}, {Carchedi}, {Noohi}}. For discussion on properties of Lie groupoids and Lie groupoid extensions, we mostly relied on \cite{{MR2157566}, {MR2012261}}. The papers \cite{{MR1950948}, {MR1876068}, {MR3150770}, {Brylinski}, {Borceux}, {MR1301844}} are some of the other articles and books which we have consulted, and found useful for this paper.

	\subsection*{Outline and organization of the paper.} Sections \ref{Section:GerbeAsLieGroupoidExtension} and \ref{Section:AGerbeAsaMorphismofStacks} are mainly devoted to review the existing definitions, collect the known results, and set up our notations and conventions. We would like to point out that in the sections mentioned above, occasionally we have included alternate proofs of the already known results, and given our interpretations for the already existing materials. For example, Lemma \ref{Lemma:MorphismfromMtoBG} has already been stated in \cite{MR2778793}. In this paper, we have given an alternate proof for the same. The primary objective of Section \ref{Section:GerbeAsLieGroupoidExtension} is to introduce the differentiable gerbe over a stack as defined in \cite{MR2493616}. In this section, on most occasions, we have followed the notations of \cite{MR2493616}. The section starts with the definition of a Lie groupoid. After recalling some basic properties of a Lie groupoid, we introduce Lie groupoid extensions, Morita equivalence of Lie groupoids, and Morita equivalent Lie groupoid extensions. We state the definition of ``a gerbe over a stack" as a Morita equivalence class of Lie groupoid extensions in Definition  \ref{Definition:GerbeoverstackVersion1}. We end Section \ref{Section:GerbeAsLieGroupoidExtension} after describing the fiber product of a pair of Lie groupoids. Section \ref{Section:AGerbeAsaMorphismofStacks} introduces the `other' definition of ``a gerbe over a stack", which will be compared with the first definition. We refer to the paper by Behrend and Xu \cite{MR2817778} for this definition. The first few subsections of Section \ref{Section:AGerbeAsaMorphismofStacks} are spent on introducing and discussing various definitions, such as category fibered in groupoids and morphisms between categories fibered in groupoids; in particular Definitions \ref{Definition:Stack} and \ref{Definition:MorphismOfStacks} state the definitions of a stack and a morphism of stacks respectively. To avoid any ambiguity of terminologies, note that we will be interested in a special type of stacks, called differentiable stacks (Definition \ref{Definition:differentiablestack}). The notion of an atlas for a differentiable stack is mentioned in Definition \ref{Definition:differentiablestack}. Also in this section we introduce two of the most important ideas for this paper, (1) a principal $\mc{G}$-bundle over a manifold $B$ for a given Lie groupoid $\mc{G}$ and maps between them, in Definitions \ref{Definition:principalLiegroupoidbundle}--\ref{Definition:morphismofprincipalLiegroupoidbundles} and (2) the $\mc{G}-\mc {H}$ bibundle for a pair of Lie groupoids $\mc {G}, \mc {H}$ in Definition \ref{Definition:GHbibundle}. The crucial observation here is that the category $B\mc G$ of principal $\mc G$ bundles is a differentiable stack and any morphism of stacks from $B\mc G$ to $B\mc H$ is classified by a $\mc G-\mc H$-bibundle (Lemma \ref{Lemma:mapBGtoBHdeterminesGHbibundle}). Finally, we introduce the notion of a gerbe over a stack (Definition \ref{Definition:GerbeoverstackasinMR2817778}), as defined in \cite{MR2817778}. That is a morphism of stacks, which is an epimorphism and whose corresponding diagonal map is also an epimorphism. Section \ref{Section:gerbegivingLiegroupoid} deals with the correspondence between two definitions, respectively discussed in Sections \ref{Section:GerbeAsLieGroupoidExtension} and \ref{Section:AGerbeAsaMorphismofStacks}, in one direction. The section starts with considering a pair of differentiable stacks $\pi_{\mc{D}}\colon\mc{D}\ra \text{Man}$ and $\pi_{\mc{C}}\colon\mc{C}\ra \text{Man}$, along with a morphism of stacks $F\colon\mc{D}\rightarrow \mc{C}$, such that $F$ is a gerbe over a stack as per the definition in \cite{MR2817778}; that is Definition \ref{Definition:GerbeoverstackasinMR2817778} in this paper. Note that this means both $F$ and the corresponding diagonal morphism $\Delta_F\colon \mc{D}\ra \mc{D}\times_{\mc{C}}\mc{D}$ are epimorphisms. This section aims to explore the possibility of finding a Morita equivalence class of a Lie groupoid extension (i.e., a gerbe over a stack as per the definition given in \cite{MR2493616} and Definition  \ref{Definition:GerbeoverstackVersion1} in this paper) from the above morphism of stacks $F\colon \mc{D}\rightarrow \mc{C}$. The notion of an atlas for a stack plays a pivotal role in our construction in this section. The first few subsections discuss the existence of `compatible' atlases for $\mc{D}$ and $\mc{C}$. However, we were required to assume that the diagonal morphism $\Delta_F\colon \mc{D}\ra\mc{D}\times_{\mc{C}}\mc{D}$ is a ``representable surjective submersion" (which is a slightly stronger condition than an epimorphism) to find a Lie groupoid extension. Theorem \ref{Theorem:GoSgivesLiegroupoidextension} is the main result of this section, which demonstrates the existence of the desired Morita equivalent Lie groupoid extension. In Section \ref{Section:GoSassociatedtoALiegrupoidExtension} we describe the construction in the other direction. We consider a Morita equivalence class of a Lie groupoid extension and recover a morphism of stacks satisfying the required properties of a gerbe over a stack defined in Definition \ref{Definition:GerbeoverstackasinMR2817778}. The notion of a principal Lie groupoid bundle and a bibundle (described in Section \ref{Section:AGerbeAsaMorphismofStacks}) play a significant role in this section. 
	We execute the construction in the following three steps. In Construction \ref{Construction:morphismofLiegroupoidstoMorphismofstacks} we find a morphism of stacks $B\mc G\to B\mc H$ from a given morphism of Lie groupoids $f \colon \mc G\to \mc H$. Next, we consider the special case where $f\colon \mc G\to \mc H$ is a Lie groupoid extension. We conclude in Theorem \ref{Theorem:BGBHisaGoS} that (1) the corresponding morphism of stacks is, in fact, a gerbe over a stack; that is, it satisfies the conditions of Definition \ref{Definition:GerbeoverstackasinMR2817778}, (2) a pair of Morita equivalent Lie groupoid extensions yields the same gerbe over a stack. Thus, a Morita equivalence class of a Lie groupoid extension produces a gerbe over a stack in the sense of Definition \ref{Definition:GerbeoverstackasinMR2817778}. 
	
	\section{A Gerbe over a stack as a Lie groupoid extension}
	\label{Section:GerbeAsLieGroupoidExtension}
	The purpose of this section is to recall the notion of the Morita equivalence class of a Lie groupoid extension. In this section, we mostly rely on the paper \cite{MR2493616} and introduce the notion of a gerbe over a stack as the Morita equivalence class of a Lie groupoid extension as defined there. 
	\begin{definition}[Lie groupoid]\label{Definition:Liegroupoid}
		A groupoid $\mc{G}=(\mc{G}_1\rightrightarrows \mc{G}_0)$ is said to be \textit{a Lie groupoid} if the source map $s\colon \mc{G}_1\rightarrow \mc{G}_0$, the target map $t\colon \mc{G}_1\rightarrow \mc{G}_0$ are submersions and the composition map $m\colon \mc{G}_1\times_{\mc{G}_0} \mc{G}_1\rightarrow \mc{G}_1$, the inverse map $i\colon \mc{G}_1\rightarrow \mc{G}_1$, the unit map $u\colon \mc{G}_0\rightarrow \mc{G}_1$ are smooth.
		  \end{definition}
	\begin{example}\label{Example:LiegroupoidassociatedtoManifold}Given a smooth manifold $M$, we associate a Lie groupoid whose object set is $M$, and the morphism set is $M$; all the other structure maps are identity maps. We denote this Lie groupoid by $(M\rightrightarrows M)$.
		  \end{example}
	\begin{example}\label{Example:ActionLiegroupoid}
		Let $G$ be a Lie group acting on a smooth manifold $X$. Consider the Lie groupoid $(G\times X\rightrightarrows X)$, whose source map is the projection map, target map is the action map, and other structure maps are defined similarly. This Lie groupoid will be called the \textit{action Lie groupoid} associated with the action of $G$ on $X$.
	\end{example}
	\begin{definition}[Transitive Lie groupoid]\label{Definition:TransitiveLiegroupoid}
		A Lie groupoid $\mc{G}=(\mc{G}_1\rightrightarrows \mc{G}_0)$ is said to be a \textit{transitive Lie groupoid}, if for each $a, b\in \mc{G}_0$, there exists a $g\in \mc{G}_1$ such that $s(g)=a$ and $t(g)=b$. 
		  \end{definition}
	\begin{definition}[Isotropy group]\label{Definition:isotropygroup}
		Let $\mc{G}=(\mc{G}_1\rightrightarrows \mc{G}_0)$ be a Lie groupoid. For $x\in \mc{G}_0$, the set \[\mc{G}_x=\{g\in \mc{G}_1| g\in s^{-1}(x)\cap t^{-1}(x)\},\] is called \textit{the isotropy group of $x$}.
		  \end{definition}
	\begin{proposition}\label{Proposition:IsotropyIsLiegroup}
		Given a Lie groupoid $\mc{G}=(\mc{G}_1\rightrightarrows \mc{G}_0)$, the isotropy group $\mc{G}_x$ is a Lie group for each $x\in \mc{G}_0$ (\cite[Corollary $1.4.11$]{MR2157566}).
	\end{proposition}
	\begin{definition}[Morphism of Lie groupoids]\label{Definition:morphismofliegroupoids}
		Let $\mc{G}$ and $\mc{H}$ be a pair of Lie groupoids. 
		A \textit{morphism of Lie groupoids} $\phi\colon \mc{G}\rightarrow \mc{H}$ is given by a pair of smooth maps $\phi_0\colon \mc{G}_0\rightarrow \mc{H}_0$
		and $\phi_1\colon \mc{G}_1\rightarrow \mc{H}_1$ which are compatible with structure maps of the Lie groupoids $\mc{G}$
		and $\mc{H}$.
		  \end{definition} 
	We will express a morphism of Lie groupoids by the following diagram, 
	\begin{equation}\label{Diagram:MorphismofLiegroupoids} \begin{tikzcd}
	\mc{G}_1 \arrow[dd,xshift=0.75ex,"t_{\mc{G}}"]\arrow[dd,xshift=-0.75ex,"s_{\mc{G}}"'] \arrow[rr, "\phi_1"] & & \mc{H}_1 \arrow[dd,xshift=0.75ex,"t_{\mc{H}}"]\arrow[dd,xshift=-0.75ex,"s_{\mc{H}}"'] \\
	& & \\
	\mc{G}_0 \arrow[rr, "\phi_0"] & & \mc{H}_0
	\end{tikzcd}.
	\end{equation}
	
	Alternately we denote a  morphism 
	of Lie groupoids $\phi\colon \mc{G}\rightarrow \mc{H}$ by $\phi\colon (\mc{G}_1\rightrightarrows \mc{G}_0)\rightarrow (\mc{H}_1\rightrightarrows \mc{H}_0)$ or by $(\phi_1,\phi_0)\colon (\mc{G}_1\rightrightarrows \mc{G}_0)\rightarrow (\mc{H}_1\rightrightarrows \mc{H}_0)$.
	
	\subsection{Pullback of a Lie groupoid along a surjective submersion}
	\label{Subsection:PullbackLiegroupoid}
	Let $\Gamma=(\Gamma_1\rightrightarrows\Gamma_0)$ be a Lie groupoid and
	$J\colon P_0\rightarrow \Gamma_0$ be a surjective submersion, where $P_0$ is a smooth manifold. Here we define the notion of pullback of the Lie groupoid $(\Gamma_1\rightrightarrows \Gamma_0)$ along the map $J\colon P_0\ra \Gamma_0$.
	
	First we pullback the source map $s\colon \Gamma_1\rightarrow \Gamma_0$ along $J\colon P_0\rightarrow \Gamma_0$ to obtain $P_0\times_{\Gamma_0}\Gamma_1$. We then pullback the map $J\colon P_0\rightarrow \Gamma_0$ along $t\circ pr_2\colon P_0\times_{\Gamma_0}\Gamma_1\rightarrow \Gamma_0$ to obtain $(P_0\times_{\Gamma_0}\Gamma_1)\times_{\Gamma_0}P_0$.
	
	Above pullbacks can be expressed by the following diagram,
	\begin{equation}\label{Diagram:PullbackLiegroupoid}
	\begin{tikzcd}
	( P_0\times_{\Gamma_0}\Gamma_1)\times_{\Gamma_0} P_0 \arrow[dd] \arrow[rrrr] & & & & P_0 \arrow[dd, "J"] \\
	& & & & \\
	P_0\times_{\Gamma_0}\Gamma_1 \arrow[dd,"pr_1"'] \arrow[rr,"pr_2"] & & \Gamma_1 \arrow[dd, "s"] \arrow[rr, "t"] & & \Gamma_0 \\
	& & & & \\
	P_0 \arrow[rr, "J"] & & \Gamma_0 & & 
	\end{tikzcd}.
	\end{equation}
	Denote the manifold $ (P_0\times_{\Gamma_0}\Gamma_1)\times_{\Gamma_0} P_0$ by $P_1$. The manifold $P_1$ along with $P_0$ gives a Lie groupoid $(P_1\rightrightarrows P_0)$, whose \begin{enumerate}
		\item source map $s\colon P_1\ra P_0$ is given by $(p,x,q)\mapsto p$,
		\item target map $t\colon P_1\ra P_0$ is given by $(p,x,q)\mapsto q$,
		\item composition map $m\colon P_1\times_{P_0}P_1\ra P_1$ is given by $\big((p,x,q),(q,y,r)\big)\mapsto (p,x\circ y,r)$,
		\item unit map $u\colon P_0\ra P_1$ is given by $a\mapsto (a,1_{J(a)},a)$ and
		\item inverse map $i\colon P_1\ra P_1$ is given by $(a,\gamma,b)\mapsto (b,\gamma^{-1},a)$.
	\end{enumerate}
	
	We call the Lie groupoid $(P_1\rightrightarrows P_0)$ to be the \textit{pullback groupoid} of the Lie groupoid $(\Gamma_1\rightrightarrows \Gamma_0)$ along the map $J\colon P_0\rightarrow \Gamma_0$. 
	\begin{definition}[Morita morphism of Lie groupoids]\label{Definition:MoritamorphismofLiegroupoids} Let $(Q_1\rightrightarrows Q_0)$ and $(\Gamma_1\rightrightarrows \Gamma_0)$ be a pair of Lie groupoids. 
		A morphism of Lie groupoids $(\phi_1,\phi_0)\colon (Q_1\rightrightarrows Q_0)\rightarrow (\Gamma_1\rightrightarrows \Gamma_0)$ expressed by the following diagram,
		\begin{equation}
		\begin{tikzcd}
		Q_1 \arrow[dd,xshift=0.75ex,"t"]\arrow[dd,xshift=-0.75ex,"s"'] \arrow[rr, "\phi_1"] & & \Gamma_1 \arrow[dd,xshift=0.75ex,"t"]\arrow[dd,xshift=-0.75ex,"s"'] \\
		& & \\
		Q_0 \arrow[rr, "\phi_0"] & & \Gamma_0
		\end{tikzcd}
		\end{equation} 
		is said to be a \textit{Morita morphism of Lie groupoids}, if 
		\begin{enumerate}
			\item the map $\phi_0\colon Q_0\rightarrow \Gamma_0$ is a surjective submersion and 
			\item the Lie groupoid $(Q_1\rightrightarrows Q_0)$ is isomorphic to the pullback groupoid of $(\Gamma_1\rightrightarrows \Gamma_0)$ along $\phi_0\colon Q_0\ra \Gamma_0$.     
		\end{enumerate}
	\end{definition} 
	A detailed discussion on the motivation behind this particular  definition can be found in \cite[Definition $3.5$, Remark $3.10$]{MR2778793}.
	\begin{remark}\label{Remark:MoritamorphismIsEquivalence}
		A Morita morphism of Lie groupoids $\phi\colon (Q_1\rightrightarrows Q_0)\rightarrow (\Gamma_1\rightrightarrows \Gamma_0)$ is actually an equivalence of categories. 
		  \end{remark}
	\begin{definition}[Morita equivalent Lie groupoids]\label{Definition:MoritaEquivalentLiegroupoids}
		Let $(\Gamma_1\rightrightarrows \Gamma_0)$ and $(\Delta_1\rightrightarrows \Delta_0)$ be a pair of Lie groupoids. We say that $(\Gamma_1\rightrightarrows \Gamma_0)$ and $(\Delta_1\rightrightarrows \Delta_0)$ are \textit{Morita equivalent Lie groupoids}, if there exists a third Lie groupoid $(Q_1\rightrightarrows Q_0)$ and a pair of Morita morphisms of Lie groupoids $\phi\colon (Q_1\rightrightarrows Q_0) \rightarrow (\Gamma_1\rightrightarrows \Gamma_0)$ and $\psi\colon (Q_1\rightrightarrows Q_0)
		\rightarrow (\Delta_1\rightrightarrows \Delta_0)$.
		
		We will express the Morita equivalent Lie groupoids in the above definition by the following diagram,
		\begin{equation}
		\begin{tikzcd}
		\Delta_1\arrow[dd,xshift=0.75ex,"t"]\arrow[dd,xshift=-0.75ex,"s"'] & & Q_1 \arrow[ll,"\psi_1"'] \arrow[rr,"\phi_1"]\arrow[dd,xshift=0.75ex,"t"]\arrow[dd,xshift=-0.75ex,"s"'] \arrow[dd,xshift=0.75ex,"t"]\arrow[dd,xshift=-0.75ex,"s"'] & & \Gamma_1 \arrow[dd,xshift=0.75ex,"t"]\arrow[dd,xshift=-0.75ex,"s"'] \\
		& & & & \\
		\Delta_0 & & Q_0 \arrow[rr,"\phi_0"] \arrow[ll,"\psi_0"'] & & \Gamma_0
		\end{tikzcd}.    
		\end{equation}   
	\end{definition}
	
	\begin{definition}[Lie groupoid extension]\label{Definition:Liegroupidextension}
		Let $(Y_1\rightrightarrows M)$ be a Lie groupoid. A \textit{Lie groupoid extension of $(Y_1\rightrightarrows M)$} is given by a Lie groupoid $(X_1\rightrightarrows M)$ and a morphism of Lie groupoids $(\phi,\text{Id})\colon (X_1\rightrightarrows M)\ra (Y_1\rightrightarrows M)$ such that, $\phi\colon X_1\ra Y_1$ is a surjective submersion.
		  \end{definition}
	We denote a Lie groupoid extension by $\phi\colon X_1\ra Y_1\rightrightarrows M$. In \cite{MR2493616} a weaker notion of a Lie groupoid extension has been used, where $\phi\colon X_1\rightarrow Y_1$ is a fibration. Here we will consider $\phi$ to be a surjective submersion. Note that every (smooth) fibration is a surjective submersion.
	
	The notions of a Morita morphism of Lie groupoids and Morita equivalent Lie groupoids extends respectively to the notions of a Morita morphism of Lie groupoid extensions and Morita equivalent Lie groupoid extensions, as we explain in the following subsection.
	\subsection{A Morita morphism of Lie groupoid extensions}
	\label{Subsection:MoritamorphismofGroupoidExtensions} Let
	$\phi'\colon X_1'\rightarrow Y_1'\rightrightarrows M'$ and 
	$\phi\colon X_1\rightarrow Y_1\rightrightarrows M$ be a pair of Lie groupoid extensions as in the following diagrams, \begin{equation}\begin{tikzcd}
	X_1' \arrow[dd,xshift=0.75ex,"t"]\arrow[dd,xshift=-0.75ex,"s"'] \arrow[rr, "\phi'"] & & Y_1' \arrow[dd,xshift=0.75ex,"t"]\arrow[dd,xshift=-0.75ex,"s"'] \\
	& & \\
	M' \arrow[rr, "\text{Id}"] & & M'
	\end{tikzcd} ~~~ \text{ and } ~~~ \begin{tikzcd}
	X_1 \arrow[dd,xshift=0.75ex,"t"]\arrow[dd,xshift=-0.75ex,"s"'] \arrow[rr, "\phi"] & & Y_1 \arrow[dd,xshift=0.75ex,"t"]\arrow[dd,xshift=-0.75ex,"s"']\\
	& & \\
	M \arrow[rr, "\text{Id}"] & & M
	\end{tikzcd}. \end{equation} 
	The most natural way of defining a morphism of Lie groupoid extensions \[\big(\phi'\colon X_1'\rightarrow Y_1'\rightrightarrows M'\big) \ra \big(\phi\colon X_1\rightarrow Y_1\rightrightarrows M\big)\] 
	would be by giving a pair of morphisms of Lie groupoids $(\psi_X,f)\colon (X_1'\rightrightarrows M')\rightarrow (X_1\rightrightarrows M)$ and 
	$(\psi_Y,g)\colon (Y_1'\rightrightarrows M')\rightarrow (Y_1\rightrightarrows M)$
	such that they are compatible with the morphisms of Lie groupoids $\phi'\colon X_1'\rightarrow Y_1'\rightrightarrows M'$ and $\phi\colon X_1\rightarrow Y_1\rightrightarrows M$ as in the following commutative diagram,
	\begin{equation}
	\begin{tikzcd}
	X_1 \arrow[dd,xshift=0.75ex,"t"]\arrow[dd,xshift=-0.75ex,"s"'] \arrow[rrrrrrr, "\phi", bend left] & & X_1' \arrow[ll, "\psi_X"'] \arrow[rrr, "\phi'"] \arrow[dd,xshift=0.75ex,"t"]\arrow[dd,xshift=-0.75ex,"s"'] & & & Y_1' \arrow[rr, "\psi_Y"] \arrow[dd,xshift=0.75ex,"t"]\arrow[dd,xshift=-0.75ex,"s"'] & & Y_1 \arrow[dd,xshift=0.75ex,"t"]\arrow[dd,xshift=-0.75ex,"s"'] \\
	& & & & & & & \\
	M \arrow[rrrrrrr, "\text{Id}", bend right] & & M' \arrow[rrr, "\text{Id}"] \arrow[ll, "f"'] & & & M' \arrow[rr, "g"] & & M
	\end{tikzcd}.
	\end{equation}
	By compatible, we mean $\psi_Y\circ \phi'=\phi\circ\psi_X$ and $f=g$. 
	\begin{definition}[Morita morphism of Lie groupoid extensions]\label{Definition:MoritamorphismofLiegroupoidExtension}
		Let $\phi'\colon X_1'\rightarrow Y_1'\rightrightarrows M'$ and 
		$\phi\colon X_1\rightarrow Y_1\rightrightarrows M$ be a pair of Lie groupoid extensions.
		A \textit{Morita morphism of Lie groupoid extensions} from $\phi'\colon X_1'\rightarrow Y_1'\rightrightarrows M'$ to 
		$\phi\colon X_1\rightarrow Y_1\rightrightarrows M$ is given by a pair of Morita morphisms of Lie groupoids, 
		\begin{equation}\begin{tikzcd}
		X_1' \arrow[dd,xshift=0.75ex,"t"]\arrow[dd,xshift=-0.75ex,"s"'] \arrow[rr, "\psi_X"] & & X_1 \arrow[dd,xshift=0.75ex,"t"]\arrow[dd,xshift=-0.75ex,"s"'] \\
		& & \\
		M' \arrow[rr, "f"] & & M
		\end{tikzcd}\text{ and }\begin{tikzcd}
		Y_1' \arrow[dd,xshift=0.75ex,"t"]\arrow[dd,xshift=-0.75ex,"s"'] \arrow[rr, "\psi_Y"] & & Y_1 \arrow[dd,xshift=0.75ex,"t"]\arrow[dd,xshift=-0.75ex,"s"'] \\
		& & \\
		M' \arrow[rr, "f"] & & M
		\end{tikzcd}\end{equation}
		such that the  diagram
		\begin{equation}\begin{tikzcd}
		X_1' \arrow[dd,"\psi_X"'] \arrow[rr, "\phi'"] & & Y_1' 
		\arrow[dd,"\psi_Y"] \\
		& & \\
		X_1 \arrow[rr, "\phi"] & & Y_1
		\end{tikzcd}\end{equation} is commutative.
		  \end{definition}
	Similarly, we define Morita equivalent Lie groupoid extensions as follows:
	\begin{definition}[Morita equivalent Lie groupoid extensions]\label{Definition:MoritaequivalentLiegroupoidextensions}
		Let $\phi'\colon X_1'\rightarrow Y_1'\rightrightarrows M'$ and 
		$\phi\colon X_1\rightarrow Y_1\rightrightarrows M$ be a pair of Lie groupoid extensions. We say that $\phi'\colon X_1'\rightarrow Y_1'\rightrightarrows M'$ and 
		$\phi\colon X_1\rightarrow Y_1\rightrightarrows M$ are \textit{Morita equivalent Lie groupoid extensions}, if there exists a third Lie groupoid extension $\phi''\colon X_1''\rightarrow Y''\rightrightarrows M''$ and a pair of Morita morphisms of Lie groupoid extensions
		\[(\phi''\colon X_1''\rightarrow Y''\rightrightarrows M'')\rightarrow (\phi\colon X_1\rightarrow Y_1\rightrightarrows M)\] and \[(\phi''\colon X_1''\rightarrow Y''\rightrightarrows M'')\rightarrow (\phi'\colon X_1'\rightarrow Y_1'\rightrightarrows M').\]
		  \end{definition}
	The Remark $2.6$ in \cite{MR2493616} states that ``There is a $1-1$ correspondence between Morita equivalence classes of Lie groupoid extensions and (equivalence classes of) differentiable gerbes over stacks", without explaining the correspondence or probing it further. 
	In a personal communication, one of the authors has clarified that their idea was to ``\textit{redefine}" the notion of a gerbe over a stack in terms of Lie groupoid extensions. The objective of this paper is to compare two `` different" definitions of a differentiable gerbe over a stack. One of the definitions is given in \cite{MR2493616}, and we state the definition below.
	\begin{definition}[Gerbe over a stack as a Lie groupoid extension \cite{MR2493616}] \label{Definition:GerbeoverstackVersion1} A \textit{gerbe over a stack} is the Morita equivalence class of a Lie groupoid extension.
		  \end{definition}
	In this paper, we will work with two notions of a $2$-fiber product. 
	The first notion is that of ``the $2$-fiber product of Lie groupoids" (Definition \ref{Definition:2-fiberproductinLiegroupoids}), and the second notion is that of ``the $2$-fiber product of categories fibered in groupoids" (Definition \ref{Definition:2-fiberproductinCFGs}). Here we will not discuss the ordinary fiber product in a category. However, before introducing $2$-fiber products, 
	it is necessary to recall the following fact regarding fiber product in the category of manifolds, which will be used frequently (for example, in Definition  \ref{Definition:2-fiberproductinLiegroupoids}). Note that, the issue of fiber product in category $\text{Man}$ will be revisited in Remark \ref{Remark:GrothendieckTopologyOnMan}.
	\begin{remark}\label{Remark:pullbackofmanifolds}
		Let $\text{Man}$ be  the category of smooth manifolds. Let $M,P,N$ be smooth manifolds with smooth maps $f\colon M\ra P, g\colon N\ra P$. Then the set theoretic pullback \[ M\times_PN=\{(m,n)\in M\times N \mid f(m)=g(n)\}\] may not have a nice smooth structure. In particular, $M\times_PN$ is not always an embedded submanifold of $M\times N$. However if $f$ and $g$ intersect transversally, in particular when one of $f$ or $g$ is a  submersion, then $M\times_PN$ is an embedded submanifold of $M\times N$.  
	\end{remark}
	
	Let $\mc{G},\mc{H},\mc{K}$ be Lie groupoids. Let $\phi\colon \mc{G}\ra \mc{K}, \psi\colon \mc{H}\ra \mc{K}$ be morphisms of Lie groupoids. We define the notion of fiber product of $\mc{G}$ and $\mc{H}$ with respect to morphisms of Lie groupoids $\phi\colon \mc{G}\ra \mc{K}, \psi\colon \mc{H}\ra \mc{K}$. As $\mc{G},\mc{H},\mc{K}$ are categories and $\phi\colon \mc{G}\ra \mc{K}, \psi\colon \mc{H}\ra \mc{K}$ are functors, we call the fiber product in this case to be the $2$-fiber product (\cite[I.2.2 Weak 2-pullbacks]{Carchedi}). Consider the groupoid $\mc{G}\times_{\mc{K}}\mc{H}$, whose object set is given by 
	\[(\mc{G}\times_{\mc{K}}\mc{H})_0=\big\{(a,\alpha,b)\mid a\in \mc{G}_0, b\in \mc{H}_0,\alpha\colon \phi(a)\ra \psi(b) \in \mc{K}_1\big\}.\]
	Given $(a,\alpha,b ),(a',\alpha',b')\in (\mc{G}\times_{\mc{K}}\mc{H})_0$,
	an arrow $(a,\alpha,b)\ra (a',\alpha',b')$ in $\mc{G}\times_{\mc{K}}\mc{H}$ is given by an arrow $u\colon a\ra a'$ in $\mc{G}$, an arrow $v\colon b\ra b'$ in $\mc{H}$ such that, $\alpha'\circ \phi(u)=\psi(v)\circ \alpha$; that is, 
	\[\text{Mor}_{\mc{G}\times_{\mc{K}}\mc{H}}\big((a,\alpha,b),(a',\alpha,b')\big)=\{u\colon a\ra a'\in \mc{G}_1, v\colon b\ra b'\in \mc{H}_1| \alpha'\circ \phi(u)=\psi(v)\circ \alpha\}.\]
	
	Observe that, the object set $(\mc{G}\times_{\mc{K}}\mc{H})_0$
	can be identified as follows: \begin{equation}\label{Equation:GKH0}(\mc{G}\times_{\mc{K}}\mc{H})_0= \mc{G}_0\times_{\phi,\mc{K}_0,s}\mc{K}_1\times_{t\circ pr_2,\mc{K}_0,\psi}\mc{H}_0.\end{equation}
	
	Under the above identification (Equation \ref{Equation:GKH0}), it would be convenient to view the object set of $\mc{G}\times_{\mc{K}}\mc{H}$ as the following pullback diagram,
	\begin{equation}\label{Diagram:ObjectSetofGKH}
	\begin{tikzcd}
	\mc{G}_0\times_{\phi,\mc{K}_0,s}\mc{K}_1\times_{t\circ pr_2,\mc{K}_0,\psi}\mc{H}_0 \arrow[dd,"pr_1"'] \arrow[rrrr,"pr_2"] & & & & \mc{H}_0 \arrow[dd, "\psi"] \\
	& & & & \\
	\mc{G}_0\times_{\phi,\mc{K}_0,s}\mc{K}_1 \arrow[dd,"pr_1"'] \arrow[rr, "pr_2"] & & \mc{K}_1 \arrow[dd, "s"] \arrow[rr, "t"] & & \mc{K}_0 \\
	& & & & \\
	\mc{G}_0 \arrow[rr, "\phi"] & & \mc{K}_0 & & 
	\end{tikzcd}.\end{equation}
	
	Likewise, the morphism set $(\mc{G}\times_{\mc{K}}\mc{H})_1$
	can be identified with:
	\begin{equation}\label{Equation:GKH1}
	(\mc{G}\times_{\mc{K}}\mc{H})_1=\mc{G}_1\times_{t\circ \phi,\mc{K}_0,s}\mc{K}_1\times_{t,\mc{K}_0,s\circ \psi}\mc{H}_1.\end{equation}
	
	Under the above identification (Equation \ref{Equation:GKH1}), we view the morphism set of $\mc{G}\times_{\mc{K}}\mc{H}$ as the following pullback diagram,
	\begin{equation}\label{Diagram:MorphismSetofGKH}
	\begin{tikzcd}
	{\mc{G}_1\times_{t\circ \phi,\mc{K}_0,s}\mc{K}_1\times_{t,\mc{K}_0,s\circ \psi}\mc{H}_1} \arrow[dd,"pr_1"'] \arrow[rrrr,"pr_2"] & & & & \mc{H}_1 \arrow[dd, "s\circ \psi"] \\
	& & & & \\
	\mc{G}_1\times_{\mc{K}_0}\mc{K}_1 \arrow[dd,"pr_1"'] \arrow[rr,"pr_2"] & & \mc{K}_1 \arrow[dd, "s"] \arrow[rr, "t"] & & \mc{K}_0 \\
	& & & & \\
	\mc{G}_1 \arrow[rr, "t\circ \phi"] & & \mc{K}_0 & & 
	\end{tikzcd}.\end{equation}
	
	It should be noted here that in general, $(\mc{G}\times_{\mc{K}}\mc{H})_0$ or $(\mc{G}\times_{\mc{K}}\mc{H})_1$, are not smooth manifolds. So, $\mc{G}\times_{\mc{K}}\mc{H}$ is not a Lie groupoid even when $\mc{G}$ and $\mc{H}$ are Lie groupoids. Here we state a sufficient condition for $\mc{G}\times_{\mc{K}}\mc{H}$ to be  a Lie groupoid 
	\cite[p.~5]{MR1950948}.
	
	Assume that the composition $t\circ pr_2\colon \mc{G}_0\times_{\mc{K}_0}\mc{K}_1\ra \mc{K}_0$ in Diagram \ref{Diagram:ObjectSetofGKH} is a submersion. Observe that 
	$\mc{G}_0\times_{\mc{K}_0}\mc{K}_1\times_{\mc{K}_0}\mc{H}_0$, is the pullback of $\psi\colon \mc{H}_0\ra \mc{K}_0$ along the submersion $t\circ pr_2\colon \mc{G}_0\times_{\mc{K}_0}\mc{K}_1\ra \mc{K}_0$. So, $(\mc{G}\times_{\mc{H}}\mc{K})_0=\mc{G}_0\times_{\mc{K}_0}\mc{K}_1\times_{\mc{K}_0}\mc{H}_0$ is a manifold. 
	
	Let $\Phi\colon \mc{G}_1\times_{\mc{K}_0}\mc{K}_1\ra \mc{G}_0\times_{\mc{K}_0}\mc{K}_1$ be the map given by $(g,k)\mapsto (t(g),k)$. This map is a submersion. So, the composition 
	\[ \mc{G}_1\times_{\mc{K}_0}\mc{K}_1\xra{\Phi} \mc{G}_0\times_{\mc{K}_0}\mc{K}_1\xra{t\circ pr_2} \mc{K}_0\] is also a submersion. 
	Observe that the composition $(t\circ pr_2)\circ \Phi\colon \mc{G}_1\times_{\mc{K}_0}\mc{K}_1\ra \mc{K}_0$ is equal to the map $t\circ pr_2\colon \mc{G}_1\times_{\mc{K}_0}\mc{K}_1\ra \mc{K}_0$ in Diagram \ref{Diagram:MorphismSetofGKH}. Thus, $t\circ pr_2\colon \mc{G}_1\times_{\mc{K}_0}\mc{K}_1\ra \mc{K}_0$ is a submersion. 
	Thenm
	$\mc{G}_1\times_{\mc{K}_0}\mc{K}_1\times_{\mc{K}_0}\mc{H}_1$  is the pullback of $s\circ\psi\colon \mc{H}_1\ra \mc{K}_0$ along the submersion $t\circ pr_2\colon \mc{G}_1\times_{\mc{K}_0}\mc{K}_1\ra \mc{K}_0$. So, $(\mc{G}\times_{\mc{H}}\mc{K})_1=\mc{G}_1\times_{\mc{K}_0}\mc{K}_1\times_{\mc{K}_0}\mc{H}_1$ is a manifold. Thus, assuming $t\circ pr_2 \colon \mc{G}_0\times_{\mc{K}_0} \mc{K}_1\ra \mc{K}_0$ is a submersion, we see that both the object set and the morphism set 
	\begin{align*}
	(\mc{G}\times_{\mc{K}}\mc{H})_0&= \mc{G}_0\times_{\phi,\mc{K}_0,s}\mc{K}_1\times_{t\circ pr_2,\mc{K}_0,\psi}\mc{H}_0,\\
	(\mc{G}\times_{\mc{K}}\mc{H})_1&= \mc{G}_1\times_{t\circ \phi,\mc{K}_0,s}\mc{K}_1\times_{t,\mc{K}_0,s\circ \psi}\mc{H}_1
	\end{align*}
	of $(\mc{G}\times_{\mc{K}}\mc{H})$ are manifolds. It is easy to see that this gives a Lie groupoid structure on $\mc{G}\times_{\mc{K}}\mc{H}$. 
	
	\begin{definition}[$2$-fiber product of Lie groupoids \cite{Carchedi}]
		\label{Definition:2-fiberproductinLiegroupoids} 
		Let $\mc{G},\mc{H},\mc{K}$ be Lie groupoids.
		Let $\phi\colon \mc{G}\ra \mc{K}$ and $\psi\colon \mc{H}\ra \mc{K}$ be a pair of morphisms of Lie groupoids. Assume further that, the composition $t\circ pr_2\colon \mc{G}_0\times_{\mc{K}_0}\mc{K}_1\ra \mc{K}_0$ in Diagram \ref{Diagram:ObjectSetofGKH} is a submersion. The Lie groupoid $\mc{G}\times_{\mc{K}}\mc{H}$ described above, is called the \textit{the $2$-fiber product of Lie groupoids} corresponding to morphisms of Lie groupoids $\phi\colon \mc{G}\ra \mc{K}$ and $\psi\colon \mc{H}\ra \mc{K}$. 
		  \end{definition}
	For future reference, we explicitly write down the source and target maps of this Lie groupoid $(\mc{G}\times_{\mc{H}}\mc{K})$,   
	\begin{align}
	&\begin{aligned}\label{Equation:DefinitionOfsourcemapforGHK}
	s\colon \mc{G}_1\times_{t\circ \phi,\mc{K}_0,s}\mc{K}_1\times_{t,\mc{K}_0,s\circ \psi}\mc{H}_1 &\ra \mc{G}_0\times_{\phi,\mc{K}_0,s}\mc{K}_1\times_{t\circ pr_2,\mc{K}_0,\psi}\mc{H}_0\\
	(u,\gamma,v)&\mapsto (s(u),\gamma\circ \phi(u),s(v)),\\
	\end{aligned}\\
	&\begin{aligned}\label{Equation:DefinitionOftargetmapforGHK}
	t\colon \mc{G}_1\times_{t\circ \phi,\mc{K}_0,s}\mc{K}_1\times_{t,\mc{K}_0,s\circ \psi}\mc{H}_1 &\ra \mc{G}_0\times_{\phi,\mc{K}_0,s}\mc{K}_1\times_{t\circ pr_2,\mc{K}_0,\psi}\mc{H}_0\\
	(u,\gamma,v)&\mapsto (t(u),\psi(v)\circ \gamma, t(v)).
	\end{aligned}
	\end{align}

	\section{A Gerbe over A stack as a morphism of stacks}\label{Section:AGerbeAsaMorphismofStacks}
	We have stated the definition of a gerbe over a stack as the Morita equivalence class of a Lie groupoid extension in Definition  \ref{Definition:GerbeoverstackVersion1}. Another definition of a gerbe over a stack, which is commonly used in literature (for example, in \cite{MR2817778}), is in terms of a morphism of stacks. The purpose of this section is to introduce the second definition of a gerbe over a stack given in \cite{MR2817778}. Eventually, we will compare the two definitions. Before stating the definition given in \cite{MR2817778}, we introduce a few more definitions. Firstly we recall the notion of a category fibered in groupoids. More details about categories fibered in groupoids can be found in \cite{MR2223406}.
	\begin{definition}[Category fibered in groupoids]\label{Definition:CategoryfiberedinGroupoids} Let $\mc{S}$ be a category. A \textit{category fibered in groupoids (CFG) over $\mc{S}$} is a category $\mc{D}$ with a functor $\pi\colon \mc{D}\rightarrow \mc{S}$ such that the following conditions hold:
		\begin{enumerate}
			\item Given an arrow $f\colon S'\rightarrow S$ in $\mc{S}$ and an object $\xi$ in $\mc{D}$ with $\pi(\xi)=S$, there exists an arrow $\tilde{f}\colon \xi'\rightarrow \xi$ in $ \mc{D}$ with $\pi(\xi')=S'$ and $\pi(\tilde{f})=f$. We call $\xi'$ to be a pullback of $\xi$ along $f$.
			\item Given a diagram 
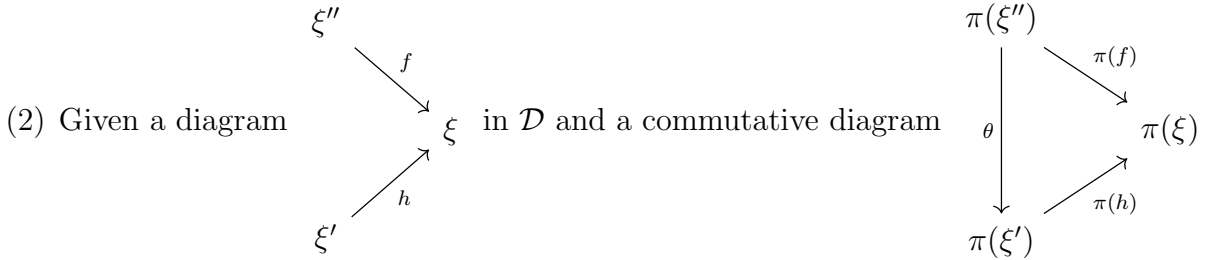
			\begin{tikzcd}
				\xi'' \arrow[rd,"f"] & \\
				& \xi \\
				\xi' \arrow[ru,"h"'] & 
			\end{tikzcd} in $\mc{D}$ and a commutative diagram \begin{tikzcd}
				\pi(\xi'') \arrow[rd,"\pi(f)"]\arrow[dd,"\theta"'] & \\
				& \pi(\xi) \\
				\pi(\xi') \arrow[ru,"\pi(h)"']& 
			\end{tikzcd}
			in $\mc{S}$ there exists a unique arrow $g\colon \xi''\rightarrow \xi'$ in $\mc{D}$ such that $h\circ g=f$ and $\pi(g)=\theta$.
		\end{enumerate}
		  \end{definition}
	We denote a category fibered in groupoids by the triple $(\mc{D},\pi,\mc{S})$.
	We write $\xi\mapsto S$ to mean $\pi(\xi)=S$. 
	We will express the second condition of the definition by the following diagram,
	\begin{equation}\label{Diagram:CategoryfiberedinGroupoids}
	\begin{tikzcd}
	\xi'' \arrow[rd,"f"] \arrow[dd, dotted,"g"] \arrow[rrr,maps to] & & & \pi(\xi'') \arrow[rd,"\pi(f)"] \arrow[dd,"\theta"'] & \\
	& \xi \arrow[rrr, maps to] & & & \pi(\xi) \\
	\xi' \arrow[ru,"h"'] \arrow[rrr, maps to] & & & \pi(\xi') \arrow[ru,"\pi(h)"'] & 
	\end{tikzcd}.
	\end{equation}
	
	\begin{definition}[Morphism of categories fibered in groupoids]\label{Definition:morphismofCFGs}
		Let $\mc{S}$ be a category. Let $\pi_{\mc{D}}\colon \mc{D}\rightarrow \mc{S}$ and $\pi_{\mc{C}}\colon \mc{C}\rightarrow \mc{S}$ be categories fibered in groupoids over $\mc{S}$. \textit{A morphism of categories fibered in groupoids from $(\mc{D},\pi_{\mc{D}},\mc{S})$ to $(\mc{C},\pi_{\mc{C}},\mc{S})$} is given by a functor $F\colon \mc{D}\rightarrow \mc{C}$ such that $\pi_{\mc{C}}\circ F=\pi_{\mc{D}}$.
		We will express the above morphism of categories fibered in groupoids by the following diagram,
		\begin{equation}\label{Diagram:morphismofCFGs}
		\begin{tikzcd}
		\mc{D} \arrow[rd, "\pi_{\mc{D}}"'] \arrow[rr, "F"] & & \mc{C} \arrow[ld, "\pi_{\mc{C}}"] \\
		& \mc{S} & 
		\end{tikzcd}.
		\end{equation}
		  \end{definition}
	\begin{definition}[Fiber of an object]\label{Definition:FiberOveranObject} 
		Let $(\mc{F},\pi,\mc{S})$ be a category fibered in groupoids. Given an object $U$ of $\mc{S}$, \textit{the fiber over $U$ in $\mc{F}$} is the subcategory $\mc{F}(U)$ of $\mc{F}$, whose 
		\begin{align*}
		\text{Ob}(\mc{F}(U))&=\{\xi\in \mc{F}_0\mid \pi(\xi)=U\}.\\
		\text{Mor}_{\mc{F}(U)}(\xi,\xi')&=\{\phi\colon \xi\rightarrow \xi' \text{ in } \mc{F}_1\text{ such that } \pi(\phi)=\text{Id}_U\}.
		\end{align*}  \end{definition}
	\begin{remark}\label{Remark:fiberpreserving}
		Let $(\mc{D},\pi_{\mc{D}},\mc{S})$ and $(\mc{C},\pi_{\mc{C}},\mc{S})$ be categories fibered in groupoids over the category $\mc{S}$. Let $F\colon \mc{D}\ra \mc{C}$ be a morphism of categories fibered in groupoids over the category $\mc{S}$. 
		For each object $U$ of $\mc{S}$, the morphism $F\colon \mc{D}\ra \mc{C}$ induces a functor $F(U)\colon \mc{D}(U)\ra \mc{C}(U)$. 
		  \end{remark}
	Let $\pi_{\mc{F}}\colon \mc{F}\rightarrow \text{Man}$ be a category fibered in groupoids over the category of manifolds. Let $M$ be a manifold. Given an open cover $\{U_i\rightarrow M\}$ of $M$, there is  a notion of the descent category associated to 
	$\{U_i\rightarrow M\}$, denoted by $\mc{F}(\{U_i\rightarrow M\})$, and the notion of pullback functor $\mc{F}(M)\rightarrow \mc{F}(\{U_i\rightarrow M\})$ (\cite[Definition $4.9$, Remark $4.10$]{MR2778793}). More details about this can be found in \cite{MR2778793} and \cite{MR2223406}.
	\begin{definition}[Stack]\label{Definition:Stack}
		Let $\pi_{\mc{D}}\colon \mc{D}\rightarrow \text{Man}$ be a category fibered in groupoids over the category of manifolds. We call $\pi_{\mc{D}}\colon \mc{D}\rightarrow \text{Man}$  \textit{a stack over the category of manifolds} if for any manifold $M$ and any open cover $\{U_i\rightarrow M\}$ of $M$, the pullback functor 
		\[\mc{D}(M)\rightarrow \mc{D}(\{U_i\rightarrow M\})\] is an equivalence of categories. 
		  \end{definition}
	
	A more general notion of a stack over a site $\mc{C}$ (a site is a category with a specified Grothendieck topology) can be found in \cite{MR2223406}. In this paper, we restrict our attention to stacks over the category of manifolds as defined in Definition \ref{Definition:Stack}. However, these two notions are related, as we explain in the following remark. 
	\begin{remark}\label{Remark:GrothendieckTopologyOnMan}
		In order to define a Grothendieck topology on the category $\text{Man}$, we need to associate, for each manifold $U$, a collection $\mc{O}_U$ of covering families which behave well under the pullback operation. Recall, we have observed in Remark \ref{Remark:pullbackofmanifolds} that, though the category $\text{Man}$ does not admit arbitrary fiber product when the collection of arrows are submersions to a given manifold, they are well behaved under the pullback operations. 
		
		Now, given an object $M$ of the category $\text{Man}$, we declare an open cover $\{U_\alpha\}$ of $M$ to be a covering family of $M$. It is straight forward to see that this gives a Grothendieck topology on the category $\text{Man}$. We call this the \textit{open cover site} on Man. Then, a stack over Man defined in Definition \ref{Definition:Stack} is the same as the stack over Man (with the open cover site) defined in \cite{MR2223406}.
		  \end{remark}
	
	\begin{definition}[Morphism of stacks]\label{Definition:MorphismOfStacks}
		A \textit{morphism of stacks} from a stack $(\mc{D},\pi_{\mc{D}}, \text{Man})$ to another stack $(\mc{C},\pi_{\mc{C}}, \text{Man})$ is a functor $F\colon \mc{D}\rightarrow \mc{C}$ such that $\pi_{\mc{C}}\circ F=\pi_{\mc{D}}$. We call a morphism of stacks $F\colon \mc{D}\rightarrow \mc{C}$ \textit{an isomorphism of stacks}, if $F\colon \mc{D}\rightarrow \mc{C}$ is an equivalence of categories. 
		  \end{definition}
	The definition of a gerbe over a stack requires the notion of the diagonal morphism associated to a morphism of categories fibered in groupoids.
	\begin{definition}[$2$-fiber product of categories fibered in groupoids]\label{Definition:2-fiberproductinCFGs}Let $\mc{S}$ be a category.
		Let $\pi_{\mc{X}}\colon \mc{X}\rightarrow \mc{S}$, $\pi_{\mc{Y}}\colon{\mc{Y}}\rightarrow \mc{S}$ and $\pi_{\mc{Z}}\colon  {\mc{Z}}\rightarrow \mc{S}$ be categories fibered in groupoids. 
		Let    $f\colon \mc{Y}\rightarrow \mc{X}, g\colon \mc{Z}\rightarrow \mc{X}$ be a pair of morphisms of categories fibered in groupoids. We define the \textit{2-fiber product of $ Y $ and $ Z $ with respect to morphisms $f,g$} to be the groupoid $\mc{Y}\times_{\mc{X}}\mc{Z}$, whose object set is given by
		\[(\mc{Y}\times_{\mc{X}}\mc{Z})_0=\big\{(y,z,\alpha)\in \mc{Y}_0\times \mc{Z}_0 \times \mc{X}_1 \mid~ \pi_Y(y)=\pi_Z(z), \alpha\colon f(y)\rightarrow g(z)\big\}.\]
		
		Given $(y,z,\alpha), (y',z',\alpha')\in (\mc{Z}\times_{\mc{X}}\mc{Y})_0$,
		an arrow $(y,z,\alpha)\ra (y',z',\alpha')$ in $\mc{Z}\times_{\mc{X}}\mc{Y}$ is given by an arrow $u\colon y\ra y'$ in $\mc{Y}$ and an arrow $v\colon z\ra z'$ in $\mc{Z}$ such that, $\alpha'\circ f(u)=g(v)\circ \alpha$ ; that is, 
		\[\text{Mor}_{\mc{Y}\times_{\mc{X}}\mc{Z}}
		\big((y,z,\alpha), (y',z',\alpha')\big)=\{u\colon y\ra y'\in \mc{Y}_1, v\colon z\ra z'\in \mc{Z}_1| \alpha'\circ f(u)=g(v)\circ \alpha\}.\]
		The groupoid $\mc{Y}\times_{\mc{X}}\mc{Z}$ comes with the following $2$-commutative diagram,
		\begin{equation}
		\begin{tikzcd}
		\mc{Y}\times_{\mc{X}}\mc{Z} \arrow[dd, "pr_2"'] \arrow[rr, "pr_1"] & & \mc{Y} \arrow[dd, "f"] \arrow[rddd, "\pi_{\mc{Y}}"] & \\
		& & & \\
		\mc{Z} \arrow[rr, "g"] \arrow[Rightarrow, shorten >=20pt, shorten <=20pt, uurr] \arrow[rrrd, "\pi_{\mc{Z}}"'] & & \mc{X} \arrow[rd, "\pi_{\mc{X}}"] & \\
		& & & \mc{S}
		\end{tikzcd}\end{equation}
		The functor $\pi_{f,g}\colon \mc{Y}\times_{\mc{X}}\mc{Z}\ra \mc{S}$ given by composition $\mc{Y}\times_{\mc{X}}\mc{Z}\xra{pr_1}\mc{Y}
		\xra{\pi_{\mc{Y}}}\mc{S}$ or 
		$\mc{Y}\times_{\mc{X}}\mc{Z}\xra{pr_2}\mc{Z}
		\xra{\pi_{\mc{Z}}}\mc{S}$ turns $\mc{Y}\times_{\mc{X}}\mc{Z}$ into a category fibered in groupoids over $\mc{S}$. We call this category $\mc{Y}\times_{\mc{X}}\mc{Z}$ \textit{the $2$-fiber product of $Y$ and $Z$ with respect to the morphisms
			$f\colon \mc{Y}\rightarrow \mc{X}$ and
			$g\colon \mc{Z}\rightarrow \mc{X}$}.   \end{definition}
	See \cite{MR2778793} and \cite{MR2223406} for a more extensive discussion on $2$-fiber product. 
	
	The definition of $2$-fibered product for categories fibered in groupoids naturally extend for stacks. Let $\pi_{\mc{D}}\colon \mc{D}\ra \text{Man}, \pi_{\mc{D}'}\colon \mc{D}'\ra \text{Man}$ and $\pi_{\mc{C}}\colon \mc{C}\ra \text{Man}$ be stacks.
	Let $F\colon \mc{D}\ra \mc{C}$ and $G\colon \mc{D}'\ra \mc{C}$ be a pair of morphisms of stacks. Then we have the category fibered in groupoids $\mc{D}\times_{\mc{C}}\mc{D}'$ over $\text{Man}$, as described in Definition \ref{Definition:2-fiberproductinCFGs} and the subsequent passage. In fact, $\mc{D}\times_{\mc{C}} \mc{D}'\ra \text{Man}$ is a stack.
	
	\begin{definition}[$2$-fiber product stack]\label{Definition:2-fiberproductstack}
		The stack $\mc{D}\times_{\mc{C}} \mc{D}'\ra \text{Man}$ mentioned above is called \textit{the $2$-fiber product stack} of $\mc{D}$ and $\mc{D'}$ with respect to the morphisms $F\colon \mc{D}\ra \mc{C}$ and $F'\colon \mc{D}'\ra \mc{C}$.
		  \end{definition}
	
	To define the notion of a gerbe over a stack, we need the notion of the diagonal morphism associated to a morphism of categories fibered in groupoids (stacks).
	\subsection{Diagonal morphism associated to a morphism of CFGs}\label{Subsection:DefinitionOfDiagonalmap} Let $\mc{S}$ be a category. Let $\pi_{\mc{D}}\colon \mc{D}\ra \mc{S}$ and $\pi_{\mc{C}}\colon \mc{C}\ra \mc{S}$ be categories fibered in groupoids over $\mc{S}$. Let $F\colon \mc{D}\rightarrow \mc{C}$ be a morphism of categories fibered in groupoids (Definition \ref{Definition:morphismofCFGs}). Consider the $2$-fiber product $\mc{D}\times_{\mc{C}}\mc{D}$ of $\mc{D}$ with itself with respect to the morphism $F:\mc{D}\ra \mc{C}$. For this $F\colon \mc{D}\ra \mc{C}$, we associate a morphism of categories fibered in groupoids $\Delta_F\colon \mc{D}\rightarrow \mc{D}\times_{\mc{C}}\mc{D}$ as follows:
	
	Given an object $a$ of $\mc{D}$ we associate the object $\big(a,a,\text{Id}\colon F(a)\rightarrow F(a)\big)$ 
	in $\mc{D}\times_{\mc{C}}\mc{D}$. 
	Given an arrow $\theta\colon a\rightarrow b$ of $\mc{D}$ we associate the arrow $(\theta,\theta)\colon \big(a,a,F(a)\rightarrow F(a)\big)\rightarrow \big(b,b,F(a)\rightarrow F(b)\big)$ in $\mc{D}\times_{\mc{C}}\mc{D}$.
	This gives a morphism of categories fibered in groupoids $\Delta_F\colon \mc{D}\rightarrow \mc{D}\times_{\mc{C}}\mc{D}$. We call this morphism $\Delta_F\colon \mc{D}\rightarrow \mc{D}\times_{\mc{C}}\mc{D}$ \textit{the diagonal morphism associated to the morphism $F\colon \mc{D}\rightarrow \mc{C}$}. 
	\begin{lemma}\label{Lemma:Diagonalmorphismofstacks}
		If $F\colon \mc{D}\ra \mc{C}$ is a morphism of stacks, then the diagonal morphism $\Delta_F\colon \mc{D}\rightarrow \mc{D}\times_{\mc{C}}\mc{D}$ is a morphism of stacks.
		  \end{lemma}
	
	We give a couple of examples that will be frequently recalled in this paper.
	\begin{example}\label{Example:StackAssociatedtoaManifold}
		Let $M$ be a smooth manifold, that is, an object in $\text{Man}$. Let $\underline{M}$ be the category whose objects are smooth maps of the form $f\colon N\rightarrow M$ for some manifold $N$. For convenience, we denote the object $f\colon N\ra M$ of $\underline{M}$ by the triple $(N,f,M)$. An arrow from an object $(N,f,M)$ of $\underline{M}$ to another object $(N',f',M)$ of $\underline{M}$ is given by a smooth map $g\colon N\rightarrow N'$ such that $f'\circ g=f$. We denote the arrow by $g\colon (N,f,M)\ra (N',f',M)$. Consider the functor $\pi_M\colon \underline{M}\rightarrow \text{Man}$ defined by 
		$(N,f,M)\mapsto N$ (at the level of objects) and $\big(g\colon (N,f,M)\ra (N',f',M)\big)\mapsto (g\colon N\ra N')$ (at the level of arrows). Then, $(\underline{M},\pi_M,\text{Man})$ is a stack over the category of manifolds. We call this $(\underline{M},\pi_M,\text{Man})$ to be the stack associated to the manifold $M$.
		  \end{example}
	\begin{example}\label{Example:CFGassociatedtoanObject}
		Let $\mc{C}$ be a category. Given an object $C$ of $\mc{C}$, we define a category $\underline{C}$ and a functor $\pi_C\colon \underline{C}\ra \mc{C}$ in the same way as we did in Example \ref{Example:StackAssociatedtoaManifold}. Then $(\underline{C},\pi_C,\mc{C})$ is a category fibered in groupoids.
		  \end{example}
	\begin{remark}\label{Remark:identificationoffiberproductstack}
		Let $X,Y$ be manifolds. Let $\pi_{\mc{D}}\colon \mc{D}\ra \text{Man}$ be a stack. Consider a pair of morphisms of stacks $p\colon \underline{X}\ra \mc{D}$ and $q\colon \underline{Y}\ra \mc{D}$. The $2$-fiber products $\underline{X}\times_{\mc{D}}\underline{Y}$ and $\underline{Y}\times_{\mc{D}}\underline{X}$ are identified as follows: 
		
		Let $\big(x,y,\alpha\colon p(x)\ra q(y)\big)$ be an object of $\underline{X}\times_{\mc{D}}\underline{Y}$.
		Let $C=\pi_X(x)=\pi_Y(y)$. As the fiber $\mc{D}(C)$ is a groupoid, every arrow in $\mc{D}(C)$ is invertible. In particular,
		$\alpha\colon p(x)\ra q(y)$ in $\mc{D}(C)$ gives an arrow $\alpha^{-1}\colon q(y)\ra p(x)$ in $\mc{D}(C)$. Thus, we have an object $(y,x,\alpha^{-1}\colon q(y)\ra p(x))$ in $\underline{Y}\times_{\mc{D}}\underline{X}$. We define a functor $\underline{X}\times_{\mc{D}}\underline{Y}\ra \underline{Y}\times_{\mc{D}}\underline{X}$ at the level of objects by 
		\[ \big(x,y,\alpha\colon p(x)\ra q(y)\big)\mapsto \big(y,x,\alpha^{-1}\colon q(y)\ra p(x)\big).\]
		Now, consider an arrow $(u,v)\colon \big(x,y,\alpha\colon p(x)\ra q(y)\big)\ra \big(x',y',\beta\colon p(x')\ra q(y')\big)$ in $\underline{X}\times_{\mc{D}}\underline{Y}$.
		For this $(u,v)$, we have the following commutative diagram,
		\[\begin{tikzcd}
		p(x) \arrow[rr, "\alpha"] \arrow[dd, "p(u)"'] & & q(y) \arrow[dd, "q(v)"] \\
		& & \\
		p(x') \arrow[rr, "\beta"] & & q(y') 
		\end{tikzcd}.\]
		As $q(v)\circ \alpha=\beta\circ p(u)$, we have $p(u)\circ \alpha^{-1}=\beta^{-1}\circ q(v)$ giving the following commutative diagram,
		\[\begin{tikzcd}
		p(y) \arrow[rr, "\alpha^{-1}"] \arrow[dd, "q(v)"'] & & q(y) \arrow[dd, "p(u)"] \\
		& & \\
		q(y') \arrow[rr, "\beta^{-1}"] & & p(x') 
		\end{tikzcd}.\]
		This diagram gives an arrow $(v,u)\colon \big(y,x,\alpha^{-1}\colon q(y)\ra p(x)\big)\ra \big(y',x',\beta^{-1}\colon q(y')\ra p(x')\big)$ in $\underline{Y}\times_{\mc{D}}\underline{X}$. 
		For the arrow $(u,v)\colon \big(x,y,\alpha\colon p(x)\ra q(y)\big)\ra \big(x',y',\beta\colon p(x')\ra q(y')\big)$ in $\underline{X}\times_{\mc{D}}\underline{Y}$ we associate the arrow $(v,u)\colon \big(y,x,\alpha^{-1}\colon q(y)\ra p(x)\big)\ra \big(y',x',\beta^{-1}\colon q(y')\ra p(x')\big)$ in $\underline{Y}\times_{\mc{D}}\underline{X}$.
		This gives a morphism of stacks $\Phi\colon \underline{X}\times_{\mc{D}}\underline{Y}\ra \underline{Y}\times_{\mc{D}}\underline{X}$. It turns out that this morphism of stacks is an isomorphism of stacks. 
		In this way we identify the $2$-fiber products $\underline{X}\times_{\mc{D}}\underline{Y}$ and 
		$\underline{Y}\times_{\mc{D}}\underline{X}$. Note that this identification has nothing to do with stacks $\underline{X}$ and $\underline{Y}$ being representable by manifolds. In fact, for any arbitrary morphisms of stacks $\mc{E}\ra \mc{D}$ and $\mc{E}'\ra \mc{D}$, we have an isomorphism of stacks $\mc{E}\times_{\mc{D}}\mc{E}'\cong \mc{E}'\times_{\mc{D}}\mc{E}$.
		  \end{remark}
	The other class of examples, of stacks that we consider next, will be arising from a Lie groupoid. For that, first we need to introduce the notion of a principal $\mc{G}$-bundle over a smooth manifold.
	\begin{definition}[Left action of a Lie groupoid on a manifold]\label{Definition:leftactionofaLiegroupoid}
		Let $\mc{G}$ be a Lie groupoid and $P$ be a smooth manifold. A \textit{left action of $\mc{G}$ on $P$} 
		consists of, 
		\begin{enumerate}
			\item a smooth map $a\colon P\rightarrow \mc{G}_0$ (called the anchor map) and
			\item a smooth map $\mu\colon \mc{G}_1\times_{s,\mc{G}_0,a}P\rightarrow P$ 
			with $(g,p)\mapsto g.p$ (called the action map)
		\end{enumerate}
		such that 
		\begin{enumerate}
			\item $a(g.p)=t(g)$ for $p\in P$ and $g\in \mc{G}_1$ with $s(g)=a(p)$,
			\item $g.(g'.p)=(g\circ g').p$ for $p\in P$ and $g,g'\in \mc{G}_1$ with $t(g')=s(g)=a(p)$ and 
			\item $1_{a(p)}.p=p$ for all $p\in P$.
		\end{enumerate}
		We will express a left action of $\mc{G}$ on $P$ by the following diagram,
		\[\begin{tikzcd}
		\mc{G}_1\arrow[dd,xshift=0.75ex,"t"]\arrow[dd,xshift=-0.75ex,"s"']& \\
		& P \arrow[ld, "a"'] \\
		\mc{G}_0 & 
		\end{tikzcd}.\]
		  \end{definition}
	A right action of a Lie groupoid on a manifold is defined likewise. But, it should be noted that,
	for a right action of $\mc{G}$ on $P$, the action map is given by \begin{equation}\label{rightact}\mu\colon P\times_{a,\mc{G}_0,t}\mc{G}_1\rightarrow P\end{equation} using the target map of $\mc{G}$ whereas, for a left action of $\mc{G}$ on $P$, the action map is given by \begin{equation}\label{Equation:leftactionofLiegroupoid}
	\mu\colon \mc{G}_1\times_{s,\mc{G}_0,a}P\rightarrow P
	\end{equation} using the source map of $\mc{G}$.
	
	We will express a right action of $\mc{G}$ on $P$ by the following diagram,
	\[\begin{tikzcd}
	& \mc{G}_1 \arrow[dd,xshift=0.75ex,"t"]\arrow[dd,xshift=-0.75ex,"s"'] \\
	P \arrow[rd, "a"] & \\
	& \mc{G}_0 
	\end{tikzcd}.\]
	\begin{definition}[principal $\mc{H}$-bundle]\label{Definition:principalLiegroupoidbundle}
		Let $\mc{H}$ be a Lie groupoid and $B$ be a smooth manifold. A \textit{principal right $\mc{H}$-bundle over $B$} consists of,
		\begin{enumerate}
			\item a smooth manifold $P$ with a right action of $\mc{H}$ on $P$ 
			and
			\item a surjective submersion $\pi\colon P\rightarrow B$,
		\end{enumerate}
		such that, 
		\begin{enumerate}
			\item the map $\pi\colon P\rightarrow B$ is $\mc{H}$-invariant; that is, 
			$\pi(p.h)=\pi(p)$ for all $p\in P$, $h\in \mc{H}_1$ with $a(p)=t(h)$ and 
			\item the map $P\times_{a,\mc{H}_0,t}\mc{H}_1\rightarrow P\times_{\pi,B,\pi} P$
			given by $(p,h)\mapsto (p,p.h)$ is a diffeomorphism.
		\end{enumerate}
		  \end{definition}
	We denote above principal $\mc{H}$-bundle by the triple $(P,\pi,B)$. We will express a principal $\mc{H}$-bundle by the following diagram, 
	\begin{equation}\label{Diagram:principalLiegroupoidbundle}
	\begin{tikzcd}
	& & \mc{H}_1 \arrow[dd,xshift=0.75ex,"t"]\arrow[dd,xshift=-0.75ex,"s"'] \\
	& P \arrow[ld, "\pi"'] \arrow[rd, "a"] & \\
	B & & \mc{H}_0 
	\end{tikzcd}.\end{equation}
	\begin{example}\label{Example:tG1G0isprincipalbundle}
		Let $\mc{G}=(\mc{G}_1\rra \mc{G}_0)$ be a Lie groupoid. Then,  the target map $t:\mc{G}_1\ra \mc{G}_0$ is a principal $\mc{G}$-bundle over the manifold $\mc{G}_0$; as in the following diagram 
		\begin{equation}\label{Diagram:tG1G0isprincipalbundle}
		\begin{tikzcd}
		& & \mc{H}_1 \arrow[dd,xshift=0.75ex,"t"]\arrow[dd,xshift=-0.75ex,"s"'] \\
		& \mc{G}_1 \arrow[ld, "t"'] \arrow[rd, "s"] & \\
		\mc{G}_0 & & \mc{G}_0 
		\end{tikzcd}.\end{equation}
	\end{example}
	\begin{remark}\label{Remark:DefinitionOfLeftprincipalHbundle}
		A principal left $\mc{H}$-bundle is defined similarly replacing a right action in Definition \ref{Definition:principalLiegroupoidbundle} by a left action.
		  \end{remark}
	\begin{remark}
		In this paper, we mostly work with the right principal bundles. So, unless otherwise stated, all principal bundles would be right principal bundles.
		  \end{remark}
	\begin{definition}[morphism of principal $\mc{H}$-bundles]\label{Definition:morphismofprincipalLiegroupoidbundles}
		Let $\mc{H}$ be a Lie groupoid. Let $(Q,\pi,M)$ and $(Q',\pi',M')$ be principal $\mc{H}$-bundles. A \textit{morphism of principal $\mc{H}$-bundles} from 
		$(Q,\pi,M)$ to $(Q',\pi',M')$ consists of a pair of smooth maps $f\colon Q\rightarrow Q'$ and $\alpha\colon M\rightarrow M'$ such that $\pi'\circ f=\alpha\circ \pi$ and $f(q.h)=f(q).h$, for all $q\in Q$ and $h\in \mc{H}_1$ satisfying $a(q)=t(h)$.
		  \end{definition}
	We denote this morphism of principal $\mc{H}$-bundles by $(f,\alpha)\colon (Q,\pi,M)\ra (Q',\pi',M')$.
	We will express the a morphism of principal $\mc{H}$-bundles by the following diagram,
	\begin{equation}\label{Diagram:morphismofprincipalLiegroupoidbundles} \begin{tikzcd}
	Q \arrow[dd, "\pi"'] \arrow[rr, "f"] & & Q' \arrow[dd, "\pi'"] \\
	& & \\
	M \arrow[rr, "\alpha"] & & M' 
	\end{tikzcd}.
	\end{equation}
	Next, we give an example of a stack associated to a Lie groupoid $\mc{G}$.
	\begin{example}\label{Example:StackAssociatedtoaLiegroupoid}
		Let $\mc{G}$ be a Lie groupoid. Let $B\mc{G}$ denote the category whose objects are principal $\mc{G}$-bundles and arrows are morphisms of principal $\mc{G}$-bundles.
		Consider the functor $\pi_{\mc{G}}\colon B\mc{G}\rightarrow \text{Man}$ given by $(Q,\pi,M)\mapsto M$ (at the level of objects) and $\big((f,\alpha)\colon (Q,\pi,M)\ra (Q',\pi',M')\big)\mapsto (\alpha \colon M\ra M')$ (at the level of arrows). Then $(B\mc{G},\pi_{\mc{G}},\text{Man})$ is a stack over the category of manifolds.
		  \end{example}
	\begin{definition} 
		The stack $(B\mc{G},\pi_{\mc{G}},\text{Man})$ (in the Example \ref{Example:StackAssociatedtoaLiegroupoid}) is called \textit{the classifying stack associated to the Lie groupoid $\mc{G}$}.
	\end{definition}
	
	Let $\mc{G}$ be a Lie groupoid. Consider the weak presheaf on the category of manifolds $\text{Man}$, defined as
	$X\mapsto \text{Hom} ((X\rightrightarrows X),\mc{G})$. Let $\text{LieGpd}$ and $\text{Gpd}$ respectively be $2$-categories of Lie groupoids and groupoids. This defines an extended Yoneda $2$-functor $\tilde{y}\colon \text{LieGpd}\rightarrow \text{Gpd}^{\text{Man}^{\text{op}}}$. This $2$-functor $\tilde{y}$ preserves all weak-limits. Given a Lie groupoid $\mc{G}$, the stack $B\mc{G}$ is isomorphic to the stackification of $\tilde{y}(\mc{G})$ \cite[p.~27]{Carchedi}.
	
	We will mainly be interested in stacks of the form $(B\mc{G},\pi_{\mc{G}},\text{Man})$ for some Lie groupoid $\mc{G}$. We call such a stack a \textit{differentiable stack}. The precise definition of a differentiable stack will be given in Definition \ref{Definition:differentiablestack}. 
	
	The morphism of stacks, we will be most interested in, will be either of the form $\underline{M}\ra B\mc{G}$ or $B\mc{G}\ra B\mc{H}$, where $M$ is a manifold, and $\mc{G},\mc{H}$ are Lie groupoids. The necessary mathematical framework will be developed in the following sections.
	
	\subsection{Pullback of a principal $\mc{G}$-bundle}
	\label{Subsection:PullbackofPrincipalLiegroupoidbundle}
	Let $\mc{G}$ be a Lie groupoid and $\pi\colon P\rightarrow B$ be a principal $\mc{G}$-bundle (Definition \ref{Definition:principalLiegroupoidbundle}). Let $f\colon N\rightarrow B$ be a smooth map. As $\pi\colon P\rightarrow B$ is a submersion, the 
	pullback $f^*P=N\times_BP=\{(n,p)\colon f(n)=\pi(p)\}$ is a manifold (an embedded submanifold of $N\times P$).
	We will express the pullback by the following diagram,
	\begin{equation}\label{Diagram:PullbackofPrincipalLiegroupoidbundle}\begin{tikzcd}
	N\times_BP \arrow[dd,"pr_1"'] \arrow[rr,"pr_2"] & & P \arrow[dd, "\pi"] \\
	& & \\
	N \arrow[rr, "f"] & & B
	\end{tikzcd}.\end{equation} 
	For our convenience, we combine the Diagrams \ref{Diagram:principalLiegroupoidbundle} and \ref{Diagram:PullbackofPrincipalLiegroupoidbundle} to draw the following diagram,
	\begin{equation}\begin{tikzcd}
	& N\times_BP \arrow[ld,"pr_1"'] \arrow[rd, "pr_2"] & & \mc{G}_1 \arrow[dd,xshift=0.75ex,"t"]\arrow[dd,xshift=-0.75ex,"s"'] \\
	N \arrow[rd, "f"] & & P \arrow[ld, "\pi"'] \arrow[rd, "a"] & \\
	& B & & \mc{G}_0
	\end{tikzcd}.\end{equation} 
	For a manifold $P$ with anchor map $a\colon P\ra \mc{G}_0$, a right action of $\mc{G}$ on $P$ 
	is given by  a map $\mu\colon P\times_{a,\mc{G}_0,t} \mc{G}_1\ra P$ (Equation \ref{rightact}). Thus, for $N\times_B P$ with anchor map $a\circ pr_2\colon N\times_BP\ra \mc{G}_0$, a right action of $\mc{G}$ on $N\times_B P$ should be given by a map 
	\[\mu\colon (N\times_B P)\times_{a\circ pr_2,\mc{G}_0,t}\mc{G}_1 \rightarrow N\times_BP.\]
	For notational simplicity, we write $(N\times_B P)\times_{\mc{G}_0}\mc{G}_1$ for $(N\times_B P)\times_{a\circ pr_2,\mc{G}_0,t}\mc{G}_1$.
	We define the map $\mu\colon (N\times_B P)\times \mc{G}_1\ra N\times_BP$ as $\big((n,p),g\big)\mapsto (n,p.g)$. 
	This map gives a right action of $\mc{G}$ on $N\times_BP$. 
	In turn we get a principal $\mc{G}$-bundle $(N\times_BP, pr_1,N)$. We call $(N\times_BP, pr_1,N)$ \textit{the pullback of the principal $\mc{G}$-bundle $\pi\colon P\ra B$ along $f\colon N\ra B$}.
	
	Let $G$ be a Lie group and $\pi\colon P\ra B$ be a principal $G$-bundle. Suppose $\pi\colon P\ra B$ has a global section, then we know that $\pi\colon P\ra B$ is a trivial $G$-bundle; that is, there exists an isomorphism of principal $G$-bundles $(P,\pi,B)\ra (B\times G,pr_1,B)$. That is to say that $(P,\pi,B)$ is isomorphic to the pullback of the (trivial) principal $G$-bundle $G\ra *$ along the map $B\ra *$. In case of Lie groupoids, the principal $\mc{G}$-bundle $t\colon \mc{G}_1\ra \mc{G}_0$ plays the role of the principal $G$-bundle $G\ra *$. Thus, we have the following result \cite[Lemma $3.19$]{MR2778793}. 
	\begin{lemma}\label{Lemma:GlobalSectionForaPrincipalLiegroupoidbundle}
		Let $\mc{G}$ be a Lie groupoid. A principal $\mc{G}$-bundle $\pi\colon P\ra B$ has a global section if and only if $(P,\pi,B)$ is isomorphic to the pullback of the (trivial) principal $\mc{G}$-bundle $(\mc{G}_1,t,\mc{G}_0)$ along a smooth map $B\ra \mc{G}_0$.
		  \end{lemma}
	Let $\pi\colon P\ra B$ be a principal $\mc{G}$-bundle. As $\pi$ is a surjective submersion, there exists an open cover $\{U_\alpha\}$ of $B$ and sections $\sigma_\alpha\colon U_\alpha\ra P$ of $\pi\colon P\ra B$. Restricting the principal $\mc{G}$-bundle $\pi\colon P\ra B$ to $U_\alpha$ gives a principal $\mc{G}$-bundle $\pi|_{\pi^{-1}(U_\alpha)}\colon \pi^{-1}(U_\alpha)\ra U_\alpha$ admitting a global section, for each $\alpha$. Thus, by Lemma \ref{Lemma:GlobalSectionForaPrincipalLiegroupoidbundle} we see that the principal $\mc{G}$-bundle $\pi|_{\pi^{-1}(U_\alpha)}\colon \pi^{-1}(U_\alpha)\ra U_\alpha$ is isomorphic to the pullback of the principal $\mc{G}$-bundle $t\colon \mc{G}_1\ra \mc{G}_0$ along a smooth map $B\ra \mc{G}_0$ for each $\alpha$. We have the following result.
	\begin{corollary}
		\label{Corollary:LocalPropertyofprincipalLiegroupoidbundle}
		Given a principal $\mc{G}$-bundle $\pi\colon P\ra B$, there exists an open cover $\{U_\alpha\}$ of $B$ such that the principal $\mc{G}$-bundle $\pi|_{\pi^{-1}(U_\alpha)}\colon \pi^{-1}(U_\alpha)\ra U_\alpha$ is the pullback of the principal $\mc{G}$-bundle $t\colon \mc{G}_1\ra \mc{G}_0$ along a smooth map $U_\alpha\ra \mc{G}_0$ for each $\alpha$.
		  \end{corollary}
	Let $(P,\pi, M)$ be a principal $\mc{G}$-bundle. Now, for 
	this principal $\mc{G}$-bundle, we associate a morphism of stacks $\underline{M}\rightarrow B\mc{G}$.
	\begin{construction}\label{Construction:BPforaprncipalLiegroupoidbundle}
		Let $M$ be a manifold and $\pi_M\colon \underline{M}\rightarrow \text{Man}$ be the stack
		associated to $M$ as in Example \ref{Example:StackAssociatedtoaManifold}. Let $\mc{G}$ be a Lie groupoid and $\pi_{\mc{G}}\colon B\mc{G}\rightarrow \text{Man}$ be the stack associated to $\mc{G}$ as in Example \ref{Example:StackAssociatedtoaLiegroupoid}. Given a principal $\mc{G}$-bundle
		$\theta\colon P\rightarrow M$, we associate a morphism of stacks $BP\colon \underline{M}\rightarrow B\mc{G}$. 
		By Remark \ref{Remark:fiberpreserving}, this morphism $BP\colon \underline{M}\ra B\mc{G}$ has to be fiber preserving, in particular, it should induces a functor $\underline{M}(N)\ra B\mc{G}(N)$ for every manifold $N$; that is, for each object $(N,g,M)$ of $\underline{M}(N)$, we should associate a principal $\mc{G}$-bundle over $N$.
		
		Let $g\colon N\rightarrow M$ be an object in $\underline{M}(N)$. We pullback the principal $\mc{G}$-bundle $\theta\colon P\rightarrow M$ along $g$ to obtain a principal $\mc{G}$-bundle $g^*P\rightarrow N$ over $N$.
		Let $BP(g)\colon g^*P\rightarrow N$ denote the projection to first coordinate and $g^*\colon g^*P\ra P$ denote the projection to second coordinate. We have the following pullback diagram,
		\begin{equation}\label{Diagram:BP(g)}\begin{tikzcd}
		g^*P \arrow[dd, "BP(g)"'] \arrow[rr,"g^*"] & & P \arrow[dd, "\theta"] \\
		& & \\
		N \arrow[rr, "g"] & & M
		\end{tikzcd}.\end{equation} 
		This gives a map $BP\colon \text{Obj}(\underline{M})\rightarrow \text{Obj}(B\mc{G})$ defined by $(N,g,M)\mapsto (g^*P,BP(g),N)$.
		Let $\Psi\colon (N,g,M)\ra (N',g',M)$ be an arrow in $\underline{M}$; that is, a smooth map $\Psi\colon N\ra N'$ such that $g=g'\circ \Psi$. We associate an arrow $BP(\Psi)\colon (g^*P,BP(g),N)\ra (g'^*P,BP(g'),N')$ in $B\mc{G}$; that is, a morphism of principal $\mc{G}$-bundles.
		
		We pullback the principal $\mc{G}$-bundle $\theta\colon P\rightarrow M$ along $g'$ to obtain the principal $\mc{G}$-bundle $BP(g')\colon g'^*P\rightarrow N'$. These principal $\mc{G}$-bundles can be expressed by the following diagram,
		\begin{equation}
		\begin{tikzcd}
		& & g'^*P \arrow[rr, "g'^*"] \arrow[dd, "BP(g')"] & & P \arrow[dd, "\theta"] \\
		& g^*P \arrow[dd, "BP(g)"] \arrow[rrru, "g^*"] & & & \\
		& & N' \arrow[rr, "g'"] & & M \\
		& N \arrow[rrru, "g"] & & & 
		\end{tikzcd}.\end{equation}
		Similarly by pulling back the principal $\mc{G}$-bundle $BP(g')\colon g'^*P\rightarrow N'$ along $\Psi\colon N\rightarrow N'$ we get the principal $\mc{G}$-bundle $\Psi^*(BP(g'))\colon \Psi^*(g'^*P)\rightarrow N$. We express the successive pullbacks by the following diagram, 
		\begin{equation}
		\begin{tikzcd}
		\Psi^*(g'^*P) \arrow[rr, "\Psi^*"] \arrow[rddd, "\Psi^*(BP(g'))"'] & & g'^*P \arrow[rr, "g'^*"] \arrow[dd, "BP(g')"] & & P \arrow[dd, "\theta"] \\
		& g^*P \arrow[dd, "BP(g)"] \arrow[rrru, "g^*"] & & & \\
		& & N' \arrow[rr, "g'"] & & M \\
		& N \arrow[ru, "\Psi"] \arrow[rrru, "g"] & & & 
		\end{tikzcd}.\end{equation}
		As $g=g'\circ \Psi$, 
		we have $BP(g'\circ\Psi)=\Psi^*(BP(g'))$. 
		Let $\Phi\colon \Psi^*(BP(g'))\rightarrow g^*P$ be the isomorphism of principal $\mc{G}$-bundles (pullback bundles are unique up to unique isomorphism). 
		The various principal $\mc{G}$-bundles and morphisms mentioned above can be expressed by a composite diagram as follows,
		\begin{equation}
		\begin{tikzcd}
		\Psi^*(g'^*P) \arrow[rr, "\Psi^*"] \arrow[rd, "\Phi"] \arrow[rddd, "\Psi^*(BP(g'))"'] & & g'^*P \arrow[rr, "g'^*"] \arrow[dd, "BP(g')"] & & P \arrow[dd, "\theta"] \\
		& g^*P \arrow[dd, "BP(g)"] \arrow[rrru, "g^*"] & & & \\
		& & N' \arrow[rr, "g'"] & & M \\
		& N \arrow[ru, "\Psi"] \arrow[rrru, "g"] & & & 
		\end{tikzcd}.\end{equation}
		Consider the composition $\Psi^*\circ \Phi^{-1}\colon g^*P\rightarrow \Psi^*(g'^*P)\ra g'^*P$. We have 
		\begin{align*}
		BP(g')\circ (\Psi^*\circ\Phi^{-1})
		&=(BP(g')\circ \Psi^*)\circ \Phi^{-1}\\
		& =(\Psi\circ \Psi^*(BP(g')))\circ\Phi^{-1}\\
		&=\Psi\circ (BP(g)\circ \Phi)\circ \Phi^{-1}\\
		&=\Psi\circ BP(g).
		\end{align*}
		Thus, we have a morphism of principal $\mc{G}$-bundles, given by the maps $\Psi^*\circ \Phi^{-1}\colon g^*P\rightarrow g'^*P$  and $\Psi\colon N\rightarrow N'$, as in the following diagram, 
		\begin{equation}\begin{tikzcd}
		g^*P \arrow[dd, "BP(g)"'] \arrow[rr, "\Psi^*\circ ~\Phi^{-1}"] & & g'^*P \arrow[dd, "BP(g')"] \\
		& & \\
		N \arrow[rr, "\Psi"] & & N'
		\end{tikzcd}.\end{equation}
		The assignments $g\mapsto BP(g)$ at the level of objects and 
		\[\big(\Psi\colon (N,g,M)\ra (N',g',M)\big)\mapsto \big((\Psi^*\circ \Phi^{-1},\Psi)\colon (g^*P,BP(g),N)\ra (g'^*P,BP(g'),N')\big)\] at the level of morphisms define a functor $BP\colon \underline{M}\rightarrow B\mc{G}$. 
		This is a morphism of stacks.  
	\end{construction}
	In conclusion, given a principal $\mc{G}$-bundle $\theta\colon P\ra M$ we have associated a morphism of stacks $BP\colon \underline{M}\rightarrow B\mc{G}$. In fact, any morphism of stacks $\underline{M}\rightarrow B\mc{G}$ is of the form $BP$ for some principal $\mc{G}$-bundle $\theta\colon P\ra M$. This result is Lemma $4.15$ in \cite{MR2778793}. Here we give an alternate proof.
	\begin{lemma}\label{Lemma:MorphismfromMtoBG}
		Let $M$ be a manifold and $\mc{G}$ be a Lie groupoid. Then, any morphism of stacks $F\colon \underline{M}\rightarrow B\mc{G}$ is of the form $BP\colon \underline{M}\rightarrow B\mc{G}$ for some principal $\mc{G}$-bundle $\theta\colon P\rightarrow M$; that is, there exists a natural isomorphism $F\Rightarrow BP\colon \underline{M}\ra B\mc{G}$.
		\begin{proof}
			Let $F\colon \underline{M}\rightarrow B\mc{G}$ be a morphism of stacks. Since $F$ is fiber preserving, as mentioned in Remark \ref{Remark:fiberpreserving}, $F\colon \underline{M}\rightarrow B\mc{G}$ induces a functor $F(N)\colon \underline{M}(N)\ra B\mc{G}(N)$ for each manifold $N$. Thus, for an object $f$ of $\underline{M}(N)$; that is, a smooth map $f\colon N\ra M$, $F(f)$ is an object of $B\mc{G}(N)$; that is, a principal $\mc{G}$-bundle over $N$. For the identity map $\text{Id}_M\colon M\rightarrow M$, $F(\text{Id}_M)$ is a principal $\mc{G}$-bundle over $M$ of the form $F(\text{Id}_M)\colon P\rightarrow M$ for some manifold $P$. For notational convenience, we denote $F(\text{Id}_M)\colon P\rightarrow M$ by $\theta\colon P\rightarrow M$. 
			
			As discussed in Construction \ref{Construction:BPforaprncipalLiegroupoidbundle}, the principal $\mc{G}$-bundle $\theta\colon P\rightarrow M$ defines the morphism of stacks $BP\colon \underline{M}\rightarrow B\mc{G}$, given by the  pullback
			of $\theta\colon P\rightarrow M$ along 
			an object/morphism of $\underline{M}$.
			We show that there is a natural isomorphism $F\Rightarrow BP\colon \underline{M}\rightarrow B\mc{G}$. 
			
			Given an object $(N,f,M)$ of $\underline{M}$, we will assign a morphism $F(f)\rightarrow BP(f)$ of principal $\mc{G}$-bundles. Given $f\colon N\ra M$ of $\underline{M}$, the Diagram \ref{Diagram:BP(g)}
			gives a morphism of principal $\mc{G}$-bundles $BP(f)\rightarrow F(\text{Id}_M)$; that is, an arrow in $B\mc{G}$. Observe that the arrow $BP(f)\ra F(\text{Id}_M)$ in $B\mc{G}$ projects to the arrow $f\colon N\ra M$ in $\text{Man}$ under the functor $\pi_{\mc{G}}\colon B\mc{G}\ra \text{Man}$.
			On the other hand, the map $f\colon N\rightarrow M$ trivially gives an arrow 
			$f\rightarrow \text{Id}_M$ in $\underline{M}$,  
			which in turn gives an arrow $F(f)\rightarrow F(\text{Id})$ in $B\mc{G}$. Observe that the arrow $F(f)\rightarrow F(\text{Id}_M)$ in $B\mc{G}$ projects to the arrow $f\colon N\rightarrow M$ in $\text{Man}$ under the functor $\pi_{\mc{G}}\colon B\mc{G}\ra \text{Man}$.
			
			The arrows $BP(f)\rightarrow F(\text{Id}_M)$ and $F(f)\rightarrow F(\text{Id}_M)$ in $B\mc{G}$ (which projects to $f\colon N\rightarrow M$ in $\text{Man}$) along with $\text{Id}_N\colon N\rightarrow N$ gives following diagram, 
			\begin{equation}\label{Diagram:F(f)toBP(f)}
			\begin{tikzcd}
			F(f) \arrow[rrrr, maps to] \arrow[rd] & & & & N \arrow[rd,"f"] \arrow[dd,"\text{Id}_N"'] & \\
			& F(\text{Id}_M) \arrow[rrrr, maps to] & & & & M \\
			BP(f) \arrow[rrrr, maps to] \arrow[ru] & & & & N \arrow[ru,"f"'] & 
			\end{tikzcd}.\end{equation}
			From Definition \ref{Definition:CategoryfiberedinGroupoids} and the Diagram \ref{Diagram:CategoryfiberedinGroupoids}, we see that the Diagram \ref{Diagram:F(f)toBP(f)} gives a unique arrow $F(f)\rightarrow BP(f)$ in $B\mc{G}$ (which projects to $\text{Id}_N\colon 
			N\rightarrow N$ in $\text{Man}$). This produces the following diagram,
			\begin{equation} \begin{tikzcd}\label{Diagram:F(f)toBP(f)part2}
			F(f) \arrow[rrrr, maps to]\arrow[dd,dotted] \arrow[rd] & & & & N \arrow[rd,"f"] \arrow[dd,"\text{Id}_N"'] & \\
			& F(\text{Id}_M) \arrow[rrrr, maps to] & & & & M \\
			BP(f) \arrow[rrrr, maps to] \arrow[ru] & & & & N \arrow[ru,"f"'] & 
			\end{tikzcd}.\end{equation} It is straightforward to see that this association of the arrow $F(f)\rightarrow BP(f)$ in $B\mc{G}$ for each object $f$ of $\underline{M}$ gives a natural transformation of functors $F\Rightarrow BP\colon \underline{M}\ra B\mc{G}$.
			Interchanging $F(f)$ and $BP(f)$ in Diagram \ref{Diagram:F(f)toBP(f)} gives an arrow $BP(f)\rightarrow F(f)$ in $B\mc{G}$ for each object $f$ of $\underline{M}$. It is easy to see that the arrows $F(f)\rightarrow BP(f)$ and $BP(f)\rightarrow F(f)$ are inverses to each other for each object $f$ of $\underline{M}$. Thus, the natural transformation $F\Rightarrow BP\colon \underline{M}\ra B\mc{G}$ is a natural isomorphism; that is, $F$ and $BP$ are naturally isomorphic functors.
		\end{proof}
	\end{lemma}
	Let $M,M'$ be manifolds and $F\colon \underline{M}\ra \underline{M'}$ be a morphism of stacks. Let $\mc{G}$ be the Lie groupoid associated to $M'$; that is, $\mc{G}=(M'\rightrightarrows M')$. Then $\underline{M'}=B\mc{G}=B(M'\rightrightarrows M')$. By Lemma \ref{Lemma:MorphismfromMtoBG}, the morphism of stacks $F\colon \underline{M}\rightarrow \underline{M'}=B\mc{G}$ is determined by a unique principal $\mc{G}$-bundle over $M$; that is, a map $f\colon M\rightarrow M'$. Explicitly, $f=F(\text{Id}_M\colon M\rightarrow M)\colon M\rightarrow M'$ determines the morphism of stacks $F\colon \underline{M}\ra \underline{M'}$. We have the following result.
	\begin{lemma}\label{Lemma:Uniquemapofmanifoldsassociatedtomapofstacks}
		Let $M,M'$ be smooth manifolds. Given a morphism of stacks $F\colon \underline{M}\rightarrow \underline{M'}$, there exists a unique map of manifolds $f\colon M\rightarrow M'$ determining $F$. 
		  \end{lemma}
	
	Suppose that $F\colon \underline{M}\rightarrow \underline{M'}$ is an isomorphism of stacks. Let $G\colon \underline{M'}\rightarrow \underline{M}$ be the inverse of $F\colon \underline{M}\rightarrow \underline{M'}$. Let $f\colon M\rightarrow M'$ be the map of manifolds associated to the morphism of stacks $F\colon \underline{M}\rightarrow \underline{M'}$ and $g\colon M'\rightarrow M$ be the map of manifolds associated to the morphism of stacks $G\colon \underline{M}\rightarrow \underline{M'}$. These maps $f\colon M\rightarrow M'$ and $g\colon M'\rightarrow M$ are such that $f\circ g=1_{M'}$ and 
	$g\circ f=1_M$; that is $M$ and $M'$ are diffeomorphic. We have the following result.
	\begin{lemma}\label{Lemma:EmbeddinOfManInCat}
		The functor $\text{Man}\rightarrow CFG$ which sends $M$ to $\underline{M}$ is an embedding of categories. 
		  \end{lemma}

	\begin{remark}\label{Remark:Uniquemanifoldrepresentingastack}
		Let $\pi_{\mc{D}}\colon \mc{D}\ra \text{Man}$ be a stack. We say that the stack $\mc{D}$ is \textit{representable by a manifold} $M$ if there exists an isomorphism of stacks $\mc{D}\cong \underline{M}$. By Lemma \ref{Lemma:EmbeddinOfManInCat}, this $M$ is unique up to diffeomorphism. We say that the stack $\mc{D}$ is \textit{representable by a Lie groupid} $\mc{G}$ if there exists an isomorphism of stacks $\mc{D}\cong B\mc{G}$. Unlike the case of manifolds, a Lie groupoid representing a stack is not uniquely determined up to an isomorphism/equivalence of categories (which was diffeomorphism in the category of manifolds). However, it is unique up to a Morita equivalence (Definition \ref{Definition:MoritaEquivalentLiegroupoids}). 
		  \end{remark} 
	The following theorem is a part of Theorem $2.26$ in \cite{MR2817778}. The converse of the Theorem below also holds. We prove it in the Section \ref{Section:GoSassociatedtoALiegrupoidExtension} (Proposition \ref{Proposition:Moritaequivalentimpliesisomorphicstacks}).
	\begin{theorem}\label{Theorem:BGBHisomorphicimpliesGHareME}
		Let $\mc{G}$ and $\mc{H}$ be Lie groupoids. If the stacks $B\mc{G}$ and $B\mc{H}$ are isomorphic, then the Lie groupoids $\mc{G}$ and $\mc{H}$ are Morita equivalent.
	\end{theorem}
	Let $\mc{G},\mc{H}$ be Lie groupoids and $B\mc{G},B\mc{H}$ be the stacks associated to $\mc{G},\mc{H}$ respectively. Our next goal is to construct an example of a morphism of stacks of the form $B\mc{G}\rightarrow B\mc{H}$. Before that we need the notion of a $\mc{G}-\mc{H}$ bibundle. For that we will mainly follow the definitions given in \cite{MR2778793}.
	\begin{definition}[$\mc{G}-\mc{H}$ bibundle]\label{Definition:GHbibundle}
		Let $\mc{G},\mc{H}$ be Lie groupoids. A \textit{$\mc{G}-\mc{H}$ bibundle} consists of,
		\begin{enumerate}
			\item a smooth manifold $P$,
			\item a left action of $\mc{G}$ on $P$ (Definition \ref{Definition:leftactionofaLiegroupoid}), with anchor map 
			$a_{\mc{G}}\colon P\rightarrow \mc{G}_0$,
			\item a right action of $\mc{H}$ on $P$ (Equation \ref{rightact}), with anchor map $a_{\mc{H}}\colon P\rightarrow \mc{H}_0$,
		\end{enumerate}
		such that,
		\begin{enumerate}
			\item the anchor map $a_{\mc{G}}\colon P\rightarrow \mc{G}_0$ is a principal $\mc{H}$-bundle,
			\item the anchor map $a_{\mc{H}}\colon P\rightarrow \mc{H}_0$ is a $\mc{G}$-invariant map; that is, $a_{\mc{H}}(g.p)=a_{\mc{H}}(p)$ for $p\in P$ and $g\in \mc{G}_1$ with $s(g)=a_{\mc{G}}(p)$,
			\item the action of $\mc{G}$ on $P$ is compatible with the action of $\mc{H}$ on $P$; that is, $(g.p).h=g.(p.h)$ for $g\in \mc{H}, p\in P$ and $h\in \mc{H}_1$ with $s(g)=a_{\mc{G}}(p)$ and $t(h)=a_{\mc{H}}(p)$.
		\end{enumerate}
		We will express a $\mc{G}-\mc{H}$ bibundle by the following diagram,
		\begin{equation} \begin{tikzcd}
		\mc{G}_1 
		\arrow[dd,xshift=0.75ex,"t"]
		\arrow[dd,xshift=-0.75ex,"s"'] & & \mc{H}_1 
		\arrow[dd,xshift=0.75ex,"t"]
		\arrow[dd,xshift=-0.75ex,"s"'] \\
		& P \arrow[rd, "a_{\mc{H}}"] \arrow[ld, "a_{\mc{G}}"'] & \\
		\mc{G}_0 & & \mc{H}_0
		\end{tikzcd}.\end{equation}
		  \end{definition}
	We denote a $\mc{G}-\mc{H}$ bibundle by $P\colon \mc{G}\ra \mc{H}$. A $\mc{G}-\mc{H}$ bibundle is called \textit{a generalized morphism of Lie groupoids}, because, given a morphism of Lie groupoids $\mc{G}\ra \mc{H}$, one can associate a $\mc{G}-\mc{H}$ bibundle (Section \ref{SubSubsection:bibundleAssociatedtoMorphismofLiegroupoids}). This justifies the notation $P\colon \mc{G}\ra \mc{H}$ for a $\mc{G}-\mc{H}$ bibundle.
	\begin{remark}\label{Remark:GPrincipalbibunde}
		If the anchor map $a_{\mc{H}}\colon P\ra \mc{H}_0$, in a $\mc{G}-\mc{H}$ bibundle $P\colon \mc{G}\ra \mc{H}$, is a principal $\mc{G}$-bundle, then, we call $P\colon \mc{G}\ra \mc{H}$ \textit{a $\mc{G}$-principal bibundle}.
		  \end{remark}
	Heuristically, a $\mc{G}-\mc{H}$ bibundle is a right principal $\mc{H}$-bundle along with a compatible action of $\mc{G}$ from the left side. 
	Given a $\mc{G}-\mc{H}$ bibundle $P\colon \mc{G}\ra \mc{H}$, one can associate a morphism of stacks $BP\colon B\mc{G}\rightarrow B\mc{H}$. We will discuss this in the  Section \ref{Section:GoSassociatedtoALiegrupoidExtension}.

For virtually the same reasons as in Lemma
\ref{Lemma:MorphismfromMtoBG}, any morphism of stacks $B\mc{G}\rightarrow B\mc{H}$ is determined by a $\mc{G}-\mc{H}$ bibundle.
\begin{lemma}\label{Lemma:mapBGtoBHdeterminesGHbibundle}
	Let $\mc{G}$ and $\mc{H}$ be a pair of Lie groupoids.
	Then, any morphism of stacks $F\colon B\mc{G}\rightarrow B\mc{H}$ is of the form $BP\colon B\mc{G}\rightarrow B\mc{H}$ for some $\mc{G}-\mc{H}$ bibundle $P\colon \mc{G}\rightarrow \mc{H}$; that is, there exists a natural isomorphism $F\Rightarrow BP\colon B\mc{G}\ra B\mc{H}$. 
	  \end{lemma}
\begin{definition}[Differentiable stack]\label{Definition:differentiablestack}
	A stack $\pi_{\mc{D}}\colon \mc{D}\rightarrow \text{Man}$ is called \textit{a differentiable stack} if there exists a smooth manifold $X$ and 
	a morphism of stacks $p\colon \underline{X}\rightarrow \mc{D}$
	satisfying the following condition:
	
	Given a smooth manifold $M$ and a morphism of stacks $f\colon \underline{M}\rightarrow \mc{D}$,
	the $2$-fibered product $\underline{M}\times_{\mc{D}}\underline{X}$ is representable by a manifold $M\times_{\mc{D}}X$ (Remark \ref{Remark:Uniquemanifoldrepresentingastack})
	and the map of manifolds $M\times_{\mc{D}}X\rightarrow M$ associated to the morphism of stacks $pr_1\colon \underline{M}\times_{\mc{D}}\underline{X}\rightarrow \underline{M}$ (Lemma \ref{Lemma:Uniquemapofmanifoldsassociatedtomapofstacks}) is a surjective submersion. We call this morphism of stacks $p\colon \underline{X}\rightarrow \mc{D}$ \textit{an atlas for the stack} $\pi_{\mc{D}}\colon \mc{D}\ra\text{Man}$.
	  \end{definition}
\begin{remark}\label{Remark:compofSurSubwithAtlasgivesAtlas}
	An atlas for a stack $\pi_{\mc{D}}\colon \mc{D}\ra\text{Man}$
	is not uniquely defined. It is easy to see that, given an atlas $p\colon \underline{X}\rightarrow \mc{D}$ for $\mc{D}$ and a surjective submersion $g\colon Y\rightarrow X$, the composition $p\circ G\colon \underline{Y}\ra \underline{X}\ra \mc{D}$ is \textit{an atlas for $\mc{D}$}.
	  \end{remark}
\begin{remark}\label{Remark:allstacksaredifferentiable}
	In this paper, unless otherwise mentioned, all stacks are differentiable stacks.
	  \end{remark}
\begin{example}
	Given a manifold $M$, the stack $(\underline{M},\pi_M,\text{Man})$ is a differentiable stack. The morphism of stacks $Id:\underline{M}\ra \underline{M}$ induced by the identity map $Id:M\ra M$ can be taken as an atlas for the stack $\underline{M}$.
	  \end{example}
\begin{example}\label{Example:AtlasforBG}
	Given a Lie groupoid $\mc{G}$, the classifying stack $(B\mc{G},\pi_{\mc{G}},\text{Man})$ is a differentiable stack. The morphism of stacks $\underline{\mc{G}_0}\ra B\mc{G}$, associated to the principal $\mc{G}$-bundle $t:\mc{G}_1\ra \mc{G}_0$ (Construction \ref{Construction:BPforaprncipalLiegroupoidbundle}), can be viewed as an atlas for the stack $B\mc{G}$. In particular, the $2$-fiber product stack $\underline{\mc{G}_0}\times_{B\mc{G}}\underline{\mc{G}_0}$ is representable by the manifold $\mc{G}_1$. We refer to Example $4.24$ in \cite{MR2778793} for further details.
	  \end{example}
\begin{example}[Quotient stack]\label{Example:Quotientstack}
	Let $G$ be a Lie group acting on a smooth manifold $X$. Let $\mc{G}=[G\times X\rra X]$ be the corresponding action Lie groupoid (Example \ref{Example:ActionLiegroupoid}). The classifying stack $B\mc{G}$ of this Lie groupoid, denoted  $[X/G]$, is called \textit{the quotient stack}.
\end{example}
\begin{remark}[\cite{Ginot}]\label{Remark:2-fiberproductisnotalwaysgeometric}
	Given a pair of morphism of differentiable stacks $B\mc{G}\ra B\mc{H}$ and $B\mc{K}\ra B\mc{H}$, the $2$-fiber product $B\mc{G}\times_{B\mc{H}} B\mc{K}$ is not in general a differentiable stack. 
	  \end{remark}
To define the notion of a gerbe over a stack as a morphism of stacks, we need the notion of an epimorphism of stacks. 
\begin{definition}[Epimorphism of stacks \cite{MR2817778}]\label{Definition:epimorphismofStacks} A morphism of stacks $F\colon \mc{D}\rightarrow \mc{C}$ is said to be \textit{an epimorphism of stacks} if given a manifold $N$ and a morphism of stacks $q\colon \underline{N}\rightarrow \mc{C}$, there exists a surjective submersion $g\colon M\rightarrow N$ and a morphism of stacks $L\colon \underline{M}\rightarrow \mc{D}$ with the following $2$-commutative diagram,
	\begin{equation} \begin{tikzcd}
	\underline{M} \arrow[dd,"L"'] \arrow[rr, "G"] & & \underline{N} \arrow[dd, "q"] \\
	& & \\
	\mc{D} \arrow[Rightarrow, shorten >=20pt, shorten <=20pt, uurr] \arrow[rr, "F"] & & \mc{C}
	\end{tikzcd}. \end{equation}
	Equivalently, given a manifold $N$ and a morphism of stacks $q\colon \underline{N}\rightarrow \mc{C}$, there exists an open cover $\{U_\alpha\rightarrow N\}$ of $N$ and a morphism of stacks $L_\alpha\colon \underline{U}_\alpha\ra \mc{D}$ for each $\alpha$ with the following $2$-commutative diagram,
	\begin{equation} \begin{tikzcd}
	\underline{U_\alpha} \arrow[dd,"L_\alpha"'] \arrow[rr] & & \underline{N} \arrow[dd, "q"] \\
	& & \\
	\mc{D} \arrow[Rightarrow, shorten >=20pt, shorten <=20pt, uurr] \arrow[rr, "F"] & & \mc{C}
	\end{tikzcd}.\end{equation}
	  \end{definition}
\begin{definition}[Representable surjective submersion \cite{MR2206877}]\label{Definition:RepSurSub}
	A morphism of stacks $F\colon \mc{D}\rightarrow \mc{C}$ is said to be \textit{representable} if for any manifold $N$ and a morphism of stacks $q\colon \underline{N}\rightarrow \mc{C}$, the $2$- fibered product (Definition \ref{Definition:2-fiberproductinCFGs}) $\mc{D}\times_{\mc{C}}\underline{N}$ is representable by a manifold $\mc{D}\times_{\mc{C}}N$. Further, if the morphism $\mc{D}\times_{\mc{C}}\underline{N}\rightarrow \underline{N}$ induces surjective submersion at the level of manifolds, then we call the morphism $F\colon \mc{D}\rightarrow \mc{C}$ \textit{a representable surjective submersion}.
	  \end{definition}
\begin{remark}
	Let $F:\mc{D}\ra\mc{C}$  be a morphism of stacks.
	It is easy to see that if $F$ is a representable surjective submersion, then $F$ is an epimorphism. 
\end{remark}
\begin{example}
	Let $M, N$ be smooth manifolds and $f:M\ra N$ be  a surjective submersion. Then, the associated morphism of stacks $F:\underline{M}\ra \underline{N}$ is a representable surjective submersion.
\end{example}
\begin{definition}[A Gerbe over a stack \cite{MR2817778}]
	\label{Definition:GerbeoverstackasinMR2817778} Let $\mc{C}$ be a differentiable stack. A morphism of stacks $F\colon \mc{D}\rightarrow \mc{C}$ is said to be \textit{a gerbe over the stack $\mc{C}$}, if the morphism $F\colon \mc{D}\rightarrow \mc{C}$ and the diagonal morphism $\Delta_F\colon \mc{D}\rightarrow \mc{D}\times_{\mc{C}}\mc{D}$ associated to $F\colon \mc{D}\ra \mc{C}$ (Section \ref{Subsection:DefinitionOfDiagonalmap}) are epimorphisms of stacks.
	  \end{definition} 

We will give an equivalent description of a gerbe over a stack in Lemma \ref{Lemma:equivalentnotionfordifferentiablestack} and illustrate the definition with several standard examples. For that purpose, first, we recall (without proof) the $2$-Yoneda Lemma (\cite[Lemma $4.19$]{MR2778793}). Note that in Lemma \ref{Lemma:MorphismfromMtoBG}, we have already observed a special case of the $2$-Yoneda Lemma.

Let $\mc{S}$ be a category and $X$ be an object of $\mc{S}$. Let $\pi_{\mc{D}}\colon \mc{D}\rightarrow \mc{S}$ be a category fibered in groupoids. Consider the functor
\begin{align*}
\Phi \colon \text{Hom}_{CFG}(\underline{X} , \mc{D}) & \rightarrow \mc{D}(X)\\
(F\colon \underline{X}\rightarrow \mc{D}) & \mapsto F(X\xrightarrow{\text{Id}} X)\\
(\alpha\colon F\Rightarrow G) & \mapsto (\alpha(X\xrightarrow{\text{Id}}X)\colon F(X\xrightarrow{\text{Id}}X)\rightarrow G(X\xrightarrow{\text{Id}}X))
\end{align*} 
where $\mc{D}(X)$ denote the fiber of $X$ in $\mc{D}$ (Definition \ref{Definition:FiberOveranObject}) and $\text{Hom}_{CFG}(\underline{X},\mc{D})$ denotes the category whose objects are morphisms of categories fibered in groupoids from $\underline{X}$ to $\mc{D}$ and whose morphisms are natural transformations. Here $\underline{X}$ is the category fibered in groupoids over $\mc{S}$, as in Example \ref{Example:CFGassociatedtoanObject} .

\begin{lemma}[$2$-Yoneda]\label{Lemma:2-yoneda}
	The functor $\Phi \colon \text{Hom}_{CFG}(\underline{X} , \mc{D}) \rightarrow \mc{D}(X)$ mentioned above is an equivalence of categories.
	  \end{lemma}

In the following, we give an application of $2$-Yoneda Lemma.

Suppose that $F\colon \mc{D}\rightarrow \mc{C}$ is an epimorphism of stacks. Given a manifold $U$ and a morphism of stacks $q\colon \underline{U}\rightarrow \mc{C}$, there exists a cover $\{U_\alpha\rightarrow U\}$ of $U$ and a morphism of stacks $L_\alpha\colon \underline{U_\alpha}\ra \mc{D}$ with the following $2$-commutative diagram, 
\begin{equation}\label{Diagram:LalphaDtoC}\begin{tikzcd}
\underline{U_\alpha} \arrow[dd,"L_\alpha"'] \arrow[rr] & & \underline{U} \arrow[dd, "q"] \\
& & \\
\mc{D} \arrow[Rightarrow, shorten >=20pt, shorten <=20pt, uurr] \arrow[rr, "F"] & & \mc{C}
\end{tikzcd}.\end{equation}
As $\pi_{\mc{C}}\colon \mc{C}\ra \text{Man}$ and $\pi_{\mc{D}}\colon \mc{D}\ra \text{Man}$ are categories fibered in groupoids, we can use the $2$-Yoneda lemma. As $U$ is an object of $\text{Man}$, the morphism $q\colon \underline{U}\rightarrow \mc{C}$ corresponds 
to an object $a$ of $\mc{C}(U)$.
As $U_\alpha$ is an object of $\text{Man}$, the morphism $L_\alpha\colon \underline{U}_\alpha\rightarrow \mc{D}$ corresponds 
to an object $x_\alpha$ of $\mc{D}(U_\alpha)$. The $2$-commutative diagram \ref{Diagram:LalphaDtoC} corresponds to an isomorphism $F(x_\alpha)\rightarrow a|_{U_\alpha}$ in $\mc{C}(U_\alpha)$ for each $\alpha$.

Thus, if a morphism of stacks $F\colon \mc{D}\rightarrow \mc{C}$ is an epimorphism, then given a manifold $U$ and an object $a$ of $\mc{C}(U)$, there exists an open cover $\{U_\alpha\rightarrow U\}$ of $U$ and objects $x_\alpha$ of $\mc{D}(U_\alpha)$ with an isomorphism $F(x_\alpha)\rightarrow a|_{U_\alpha}$ in $\mc{C}(U_\alpha)$ for each $\alpha$. It turns out that the converse is true as well.
That means the following.
%
Suppose that $F:\mc{D}\ra\mc{C}$ is a morphism of stacks with the following property: given a manifold $U$ and an object $a$ of $\mc{C}(U)$, there exists an open cover $\{U_\alpha\ra U\}$ of $U$ and objects $x_\alpha\in \mc{D}(U_\alpha)$ such that $F(x_\alpha)$ is isomorphic to $a|_{U_\alpha}$ for each $\alpha$. Then $F\colon \mc{D}\ra \mc{C}$ is an epimorphism of stacks.

Suppose that $\Delta_F\colon \mc{D}\rightarrow \mc{D}\times_{\mc{C}}\mc{D}$ is an epimorphism of stacks. Given a manifold $U$ and a morphism of stacks $q\colon \underline{U}\rightarrow \mc{D}\times_{\mc{C}}\mc{D}$, there exists a cover $\{U_\alpha\rightarrow U\}$ and a morphism of stacks $L_\alpha\colon \underline{U_\alpha}\rightarrow \mc{D}$ such that we have following $2$-commutative diagram,
\begin{equation} \begin{tikzcd}
\underline{U_\alpha} \arrow[dd,"L_\alpha"'] \arrow[rr] & & \underline{U} \arrow[dd, "q"] \\
& & \\
\mc{D} \arrow[rr, "\Delta"] \arrow[Rightarrow, shorten >=30pt, shorten <=30pt, uurr] & & \mc{D}\times_{\mc{C}}\mc{D}
\end{tikzcd}.\end{equation}
By $2$-Yoneda lemma,
the map $q\colon \underline{U}\rightarrow \mc{D}\times_{\mc{C}}\mc{D}$ corresponds to an object $\big(a,b,p\colon F(a)\rightarrow F(b)\big)$ in $(\mc{D}\times_{\mc{C}}\mc{D})(U)$; that is, $a\in \mc{D}(U)_0,b\in \mc{D}(U)_0 \text { and } p\in \mc{C}(U)_1$. 
The morphism of stacks $L_\alpha\colon \underline{U_\alpha}\rightarrow \mc{D}$ corresponds to an object $c$ of $\mc{D}(U_\alpha)$. We have $\Delta_F(c)=\big(c,c,\text{Id}\colon F(c)\rightarrow
F(c)\big)$. The $2$-commutative diagram yields an isomorphism \[\big(c,c,\text{Id}\colon F(c)\rightarrow F(c)\big)\rightarrow \big(a|_{U_\alpha},b|_{U_\alpha},p|_{U_\alpha}\colon F(a|_{U_\alpha})\rightarrow F(b|_{U_\alpha})\big).\] 
That is, there exists isomorphisms $a_{\alpha}\colon c\rightarrow a|_{U_\alpha},b_{\alpha}\colon c\rightarrow b|_{U_\alpha}$ in $\mc{D}(U_\alpha)$ satisfying the following commutative diagram,
\[\begin{tikzcd}
F(c) \arrow[dd, "F(a_\alpha)"'] \arrow[rr, "F(\text{Id})"] & & F(c) \arrow[dd, "F(b_\alpha)"] \\
& & \\
F(a|_{U_\alpha}) \arrow[rr,"p|_{U_\alpha}"] & & F(b|_{U_\alpha})
\end{tikzcd}.\]
In other words, we have $p|_{U_\alpha}\circ F(a_\alpha)=F(b_\alpha)\circ F(\text{Id})=F(b_\alpha)$. As $a_\alpha$ is an isomorphism, we have 
$p|_{U_\alpha}=F(b_\alpha)\circ F(a_{\alpha}^{-1})=F(b_\alpha\circ a_{\alpha}^{-1})$; that is, $p|_{U_\alpha}\colon F(a|_{U_\alpha})\rightarrow F(b|_{U_\alpha})$ is equal to $F(\tau_\alpha)$ for some isomorphism $\tau_\alpha\colon a|_{U_\alpha}\rightarrow b|_{U_\alpha}$.

Thus, if the diagonal morphism $\Delta_F\colon \mc{D}\rightarrow\mc{D}\times_{\mc{C}}\mc{D}$ is an epimorphism, then given a manifold $U$ and an arrow $p\colon F(a)\rightarrow F(b)$ in $\mc{C}(U)$, there exists an open cover $\{U_\alpha\rightarrow U\}$ of $U$ and a family of isomorphisms $\{\tau_\alpha\colon a|_{U_\alpha}\rightarrow b|_{U_\alpha}\}$ such that $F(\tau_\alpha)=p|_{U_\alpha}$. 
Again the converse holds.
Thus, we have the following result.
\begin{lemma}\label{Lemma:equivalentnotionfordifferentiablestack}
	A morphism of stacks $F\colon \mc{D}\rightarrow \mc{C}$ is a gerbe over a stack if and only if the following two conditions holds:
	\begin{enumerate}
		\item Given a manifold $U$ and an object $a$ of $ \mc{C}(U)$, there exists an open cover $\{U_\alpha\rightarrow U\}$ of $U$ and objects $x_\alpha$ of $\mc{D}(U_\alpha)$ with an isomorphism $F(x_\alpha)\rightarrow a|_{U_\alpha}$ in $\mc{C}(U_\alpha)$ for each $\alpha$. 
		\item Given a manifold $U$ and an arrow $p\colon F(a)\rightarrow F(b)$ in $\mc{C}(U)$, there exists an open cover $\{U_\alpha\rightarrow U\}$ of $U$ and isomorphisms $\tau_\alpha\colon a|_{U_\alpha}\rightarrow b|_{U_\alpha}$ in $\mc{D}(U_\alpha)$ such that $F(\tau_\alpha)=p|_{U_\alpha}$ in $\mc{C}(U_\alpha)$ for each $\alpha$.
	\end{enumerate}
	  \end{lemma}
\begin{example}
	Let $X$ be a smooth manifold. Let $G$ be a Lie group acting on the smooth manifold $X$. Consider a central extension of Lie groups $1\ra S^1\ra \hat{G}\xra{\pi} G\ra 1$. Let $[X/G]$
	and $[X/\hat{G}]$ respectively be the quotient stacks (Example \ref{Example:Quotientstack}) associated to the actions of $\hat{G}$ and $G$ on $X$. The morphism of Lie groups $\pi:\hat{G}\ra G$ defines a morphism of Lie groupoids $(X\times \hat{G}\rra X)\ra (X\times G\rra X)$, given by $x\mapsto x$ and $(x,\hat{g})\mapsto (x,\pi(\hat{g}))$. Then, as we will see in Section \ref{Section:GoSassociatedtoALiegrupoidExtension}, this morphism of Lie groupoids $(X\times \hat{G}\rra X)\ra (X\times G\rra X)$ associates a morphism of stacks $[X/\hat{G}]\xra{\pi} [X/G]$. Infact, this morphism of stacks $[X/\hat{G}]\xra{\pi} [X/G]$ is a gerbe over the quotient stack $[X/G]$. 
\end{example}
\begin{example}
	Let $M,N$ be manifolds. Let $f:M\ra N$ be a diffeomorphism. Then, the associated morphism of stacks $F:\underline{M}\ra \underline{N}$ is a gerbe over the stack $\underline{N}$. More over, a morphism of stacks $G:\underline{M}\ra \underline{N}$ is a gerbe over the stack $\underline{N}$ implies that the associated map of manifolds $g:M\ra N$ is a diffeomorphism.
\end{example}
\begin{remark} 
	Let $\mc{D}\ra\mc{C}$ be a gerbe over the stack $\mc{C}$. When the stack $\mc{C}$ is representable by a manifold; that is $\mc{C}\cong \underline{M}$ for a manifold $M$, we recover the notion of \textit{a gerbe over a manifold $M$}. It is immediate from Lemma \ref{Lemma:equivalentnotionfordifferentiablestack} that, a gerbe over a manifold $M$, associates a groupoid $\mc{G}(U)$ with each open set $U\subseteq M$, such that the following conditions are satisfied:
	\begin{itemize}
		\item  given $x\in M$ there is an open subset $U\subseteq M$ containing $x$ such that $\mc{G}(U)$ is non empty,
		\item  given $a,b\in\mc{G}(U)$ and $x\in U\subseteq M$, there exists an open subset $V$ of $U$ containing $x$ such that $a|_V$ is isomorphic to $b|_V$.
	\end{itemize}
	These two properties respectively, are called ``locally non empty" and ``locally connected". 
	For further details on this topic, we
	refer to the Section $3$ of \cite{Moerdijk2}.
\end{remark}
\begin{example}
	Let $M$ be a manifold and $\mc{O}(M)$ be the category of open sets of the manifold $M$. Let $G$ be a Lie group. For an open set $U\subseteq M$, let $Tor(G)|_U$ denote the groupoid of principal $G$ bundles over the manifold $U$. Then, the assignment $U\mapsto Tor(G)|_U$ for the open set $U\subseteq M$ gives a gerbe over the manifold $M$.
	  \end{example}
\section{A Lie groupoid extension associated to a Gerbe over a stack}\label{Section:gerbegivingLiegroupoid}
Let $\pi_{\mc{D}}\colon \mc{D}\ra \text{Man}$ and $\pi_{\mc{C}}\colon \mc{C}\ra \text{Man}$ be a pair of differentiable stacks. Let $F\colon \mc{D}\rightarrow \mc{C}$ be a gerbe over a stack; that is, the morphism $F\colon \mc{D}\ra \mc{C}$ and the diagonal morphism $\Delta_F\colon \mc{D}\ra \mc{D}\times_{\mc{C}}\mc{D}$ are epimorphisms of stacks. We further assume that the diagonal morphism $\Delta_F\colon \mc{D}\ra \mc{D}\times_{\mc{C}}\mc{D}$ is a representable surjective submersion. With this gerbe  $F\colon \mc{D}\ra \mc{C}$ we associate a (Morita equivalence class of) Lie groupoid extension.

The outline of this section is as follows: 
\begin{enumerate}
	\item Given an atlas $r\colon \underline{X}\rightarrow \mc{D}$ for the stack $\mc{D}$, we associate a Lie groupoid $(X\times_\mc{D}X\rightrightarrows X)$. We denote this Lie groupoid $(X\times_{\mc{D}}X\rightrightarrows X)$ by $\mc{G}_r$. We further prove that, if $r\colon \underline{X}\ra \mc{D}$ and $l\colon \underline{Y}\ra \mc{D}$ are atlases for the stack $\mc{D}$, then the corresponding Lie groupoids $\mc{G}_r=(X\times_{\mc{D}}X\rightrightarrows X)$ and $\mc{G}_l=(Y\times_{\mc{D}}Y\rightrightarrows Y)$ are Morita equivalent (Lemma \ref{Lemma:Liegroupoidrepresentingstack} and Lemma \ref{Lemma:atlasindependent}).
	\item We use the fact that $F\colon \mc{D}\ra \mc{C}$ is an epimorphism of stacks to prove that there exists an atlas $q\colon \underline{X}\ra \mc{C}$ for $\mc{C}$ and a morphism of stacks $p\colon \underline{X}\ra \mc{D}$ satisfying the following $2$-commutative diagram (Lemma \ref{Lemma:existenceofatlasforC}),
	\[\begin{tikzcd}
	\underline{X}\arrow[dd, "p"']\arrow[rrdd, "q"{name=M}] & & \\
	& & \\
	\mc{D}\arrow[rr, "F", swap]
	\arrow[Rightarrow, shorten >=10pt, shorten <=10pt, to=M] & & \mc{C}
	\end{tikzcd}.\]
	By $2$ commutative diagram above, $F\circ p$ and $q$ can be identified upto a $2$-isomorphism.
	\item The fact that the diagonal morphism $\Delta_F\colon \mc{D}\ra \mc{D}\times_{\mc{C}}\mc{D}$ is an epimorphism of stacks implies that the morphism of stacks $p\colon \underline{X}\ra \mc{D}$ obtained in step $2$ is an epimorphism of stacks (Lemma \ref{Lemma:pisanepimorphism}).
	\item Under the ``assumption" that the diagonal morphism $\Delta_F\colon \mc{D}\ra \mc{D}\times_{\mc{C}}\mc{D}$ is a representable surjective submersion, we prove that, the morphism of stacks $p\colon \underline{X}\ra \mc{D}$ mentioned in step $3$ is an atlas for the stack $\pi_{\mc{D}}\colon \mc{D}\ra \text{Man}$ (Lemma \ref{Lemma:extraassumptiononDiagonalmorphism}).
	\item For the choices made in step $(2)$ and step $(4)$ for atlases $p\colon \underline{X}\ra \mc{D}$ and $q\colon \underline{X}\ra\mc{C}$ we respectively obtain Lie groupoids $\mc{G}_p=(X\times_{\mc{D}}X\rightrightarrows X)$
	and $\mc{G}_q=(X\times_{\mc{C}}X\rightrightarrows X)$. In Lemma \ref{Lemma:MoSgivingLiegroupoidExtension} we prove that, the morphism of stacks $F\colon \mc{D}\ra \mc{C}$ gives a Lie groupoid extension $\mc{G}_p\ra \mc{G}_q$.
	\item Finally we prove that the above construction does not depend on the choice of $q\colon \underline{X}\rightarrow \mc{C}$ (Lemma \ref{Lemma:LiegroupoidExtnIsIndependentofq}).
\end{enumerate}
\subsection{A stack with an atlas giving a Lie groupoid}
\label{Subsection:Liegroupoidassociatedtoatlas}
Let $\pi_{\mc{D}}\colon \mc{D}\rightarrow \text{Man}$ be a differentiable stack. Let $r\colon \underline{X}\rightarrow \mc{D}$ be an atlas 
for the stack $\mc{D}$. 
Let $\underline{X}\times_{\mc{D}}\underline{X}$ be the $2$-fiber product expressed in the following diagram,
\[\begin{tikzcd}
\underline{X}\times_{\mc{D}}\underline{X} \arrow[dd, "pr_2"'] \arrow[rr, "pr_1"] & & \underline{X} \arrow[dd, "r"] \\
& & \\
\underline{X} 
\arrow[Rightarrow, shorten >=30pt, shorten <=30pt, uurr]
\arrow[rr, "r"] & & \mc{D} 
\end{tikzcd}.\]
As $r\colon \underline{X}\rightarrow \mc{D}$ is an atlas for the stack $\mc{D}$, the $2$-fiber product $\underline{X}\times_{\mc{D}}\underline{X} $ is representable by a manifold, which we denote by $X\times_{\mc{D}}X$; that is, there exists an isomorphism of stacks $\underline{X\times_{\mc{D}}X}\cong\underline{X}\times_{\mc{D}}\underline{X}$. Let $s\colon X\times_{\mc{D}}X\ra X$ and $t\colon X\times_{\mc{D}}X\ra X$ respectively be the morphisms of manifolds associated to the morphism of stacks $pr_1\colon \underline{X}\times_{\mc{D}}\underline{X}\rightarrow \underline{X}$ and $ pr_2\colon \underline{X}\times_{\mc{D}}\underline{X}\rightarrow \underline{X}$. 
These maps $s,t\colon X\times_{\mc{D}}X\ra X$ along with the following structure maps
gives a Lie groupoid $\mc{G}_r=(X\times_{\mc{D}}X\rightrightarrows X)$
\begin{enumerate}
	\item The composition is given by the morphism of stacks \[m\colon (\underline{X}\times_{\mc{D}}\underline{X})\times_{\underline{X}}(\underline{X}\times_{\mc{D}}\underline{X})\ra \underline{X}\times_{\mc{D}}\underline{X}\] defined (at the level of objects) as \[m\big((a,b,\alpha\colon r(a)\ra r(b)),(b,c,\beta\colon r(b)\ra r(c))\big)=\big(a,c,\beta\circ \alpha\colon r(a)\ra r(c)\big).\]
	\item The unit map is given by the morphism of stacks
	\[u\colon \underline{X}\ra \underline{X}\times_{\mc{D}}\underline{X}\] defined (at the level of objects) as \[i(a)= \big(a,a, \text{Id}\colon r(a)\ra r(a)\big).\]
	\item The inverse map is given by the morphism of stacks \[i\colon \underline{X}\times_{\mc{D}}\underline{X}\ra \underline{X}\times_{\mc{D}}\underline{X}\] defined (at the level of objects) as \[i\big((a,b,\alpha\colon r(a)\ra r(b))\big)=\big(b,a,\alpha^{-1}\colon r(b)\ra r(a)\big).\]
\end{enumerate}
To be precise, the morphism of stacks $m\colon (\underline{X}\times_{\mc{D}}\underline{X})\times_{\underline{X}}(\underline{X}\times_{\mc{D}}\underline{X})\ra \underline{X}\times_{\mc{D}}\underline{X}$ induces the map of manifolds $m\colon (X\times_{\mc{D}}X)\times_X(X\times_{\mc{D}}X)\ra X\times_{\mc{D}}X$. View this $m$ as the composition map
$(\mc{G}_r)_1\times_{(\mc{G}_r)_0}(\mc{G}_r)_1\ra (\mc{G}_r)_1$ for $\mc{G}_r$.
Similarly, the morphism of stacks $i\colon \underline{X}\times_{\mc{D}}\underline{X}\ra \underline{X}\times_{\mc{D}}\underline{X}$ and $u\colon \underline{X}\ra \underline{X}\times_{\mc{D}}\underline{X}$ induce the map of manifolds $i\colon X\times_{\mc{D}}X\ra X\times_{\mc{D}}X$ and $u\colon X\ra X\times_{\mc{D}}X$ respectively. View this $i$ as the inverse map $i\colon (\mc{G}_r)_1\ra (\mc{G}_r)_1$ and $u$ as the unit map $u\colon (\mc{G}_r)_0\ra (\mc{G}_r)_1$ for $\mc{G}_r$.

It turns out that there is an isomorphism of stacks $\mc{D}\cong B\mc{G}_r$. More details about this isomorphism $\mc{D}\cong B\mc{G}_r$ can be found in \cite{MR2778793}. 
For our purpose, we note the following result \cite[Proposition $4.31$]{MR2778793}.

\begin{lemma}\label{Lemma:Liegroupoidrepresentingstack} Let $\pi_{\mc{D}}\colon \mc{D}\rightarrow \text{Man}$ be a differentiable stack and $r\colon \underline{X}\rightarrow \mc{D}$ be an atlas for $\mc{D}$. 
	Then,
	there exists a Lie groupoid $\mc{G}$ with an isomorphism of stacks $\mc{D}\cong B\mc{G}$. Moreover, we may take $\mc{G}_0=X$ and $\mc{G}_1=X\times_{\mc{D}}X$, where $X\times_{\mc{D}}X$ is the manifold representing the $2$-fiber product $\underline{X}\times_{\mc{D}}\underline{X}$. 
	  \end{lemma}
It should be noted here that the Lie groupoid associated to a differentiable stack is independent of the choice of an atlas, up to a Morita equivalence. This can be argued as follows:

Let $p\colon \underline{X}\ra \mc{D}$ be an atlas for the stack $\mc{D}$ and $B\mc{G}_p\cong \mc{D}$ be the isomorphism of stacks.
Let $q\colon \underline{Y}\ra \mc{D}$ be another atlas for $\mc{D}$ and $B\mc{G}_q\cong \mc{D}$ be the isomorphism of stacks.
Thus, we have an isomorphism of stacks $B\mc{G}_p\ra B\mc{G}_q$.
Recall that an isomorphism of stacks $B\mc{G}\ra B\mc{H}$ gives a Morita equivalence of Lie groupoids $\mc{G}\ra \mc{H}$ (Theorem \ref{Theorem:BGBHisomorphicimpliesGHareME}, \cite[Theorem $2.26$]{MR2817778}). Using this, we can conclude that the isomorphism of stacks $B\mc{G}_p\ra B\mc{G}_q$ gives a Morita equivalence of Lie groupoids $\mc{G}_p\ra \mc{G}_q$. Thus, $\mc{G}_p$ and $\mc{G}_q$ are Morita equivalent Lie groupoids. So, we have shown the following:
\begin{lemma}\label{Lemma:atlasindependent}
	Let $\pi_{\mc{D}}\colon \mc{D}\ra \text{Man}$ be a differentiable stack. Let $p\colon \underline{X}\ra \mc{D}$ and $q\colon \underline{Y}\ra \mc{D}$ be atlases for the stack $\mc{D}$. Then, the Lie groupoids $\mc{G}_p=(X\times_{\mc{D}}X\rightrightarrows X)$ and 
	$\mc{G}_q=(Y\times_{\mc{D}}Y\rightrightarrows Y)$ are Morita equivalent.
	  \end{lemma}
\subsection{Existence of an atlas for $\mc{C}$}\label{Subsection:existenceofatlas} As $F:\mc{D}\ra \mc{C}$ is an epimorphism, given a manifold $M$ and a morphism of stacks $\tilde{q}\colon \underline{M}\ra \mc{C}$, there exists a surjective submersion $G\colon \underline{X}\ra \underline{M}$
and a morphism of stacks $p\colon \underline{X}\ra \mc{D}$ with the following $2$-commutative diagram,
\[\begin{tikzcd}
\underline{X} \arrow[dd, "p"'] 
\arrow[rr, "G"] & & \underline{M} \arrow[dd, "\tilde{q}"] \\
& & \\
\mc{D} \arrow[rr, "F"] \arrow[Rightarrow, shorten >=20pt, shorten <=20pt, uurr] & & \mc{C}
\end{tikzcd}.\]
Now, choose $\tilde{q}\colon \underline{M}\ra \mc{C}$ to be an atlas for $\mc{C}$. Since $g\colon X\ra M$ is a surjective submersion, as per the Remark \ref{Remark:compofSurSubwithAtlasgivesAtlas}, the composition $q:=(\tilde{q}\circ G)\colon \underline{X}\ra \mc{C}$ is an atlas for $\mc{C}$. Thus, we have obtained an atlas $q\colon \underline{X}\ra \mc{C}$ for $\mc{C}$ and a morphism of stacks $p\colon \underline{X}\ra \mc{D}$ with the following $2$-commutative diagram,
\begin{equation}\label{Diagram:FpToq}
\begin{tikzcd}
\underline{X}\arrow[dd, "p"']\arrow[rrdd, "q"{name=M}] & & \\
& & \\
\mc{D}\arrow[rr, "F", swap]
\arrow[Rightarrow, shorten >=10pt, shorten <=10pt, to=M] & & \mc{C}
\end{tikzcd}.
\end{equation}
\begin{lemma}\label{Lemma:existenceofatlasforC} 
	Let $F\colon \mc{D}\ra \mc{C}$ be an epimorphism of stacks. Then, there exists an atlas $q\colon \underline{X}\ra \mc{C}$ for $\mc{C}$ and a morphism of stacks $p\colon \underline{X}\ra \mc{D}$ satisfying the $2$-commutative diagram \ref{Diagram:FpToq}.
	  \end{lemma}
\subsection{Proof that $p\colon \underline{X}\ra \mc{D}$ is an epimorphism of stacks}\label{Subsection:pIsAnEpimorphism}
Let $F\colon \mc{D}\ra \mc{C}$ be a gerbe over a stack. Let $p\colon \underline{X}\ra \mc{D}$ and $q\colon \underline{X}\ra \mc{C}$ be as in Lemma \ref{Lemma:existenceofatlasforC}. Let $M$ be a manifold and $r\colon \underline{M}\ra \mc{D}$ be a morphism of stacks. Consider the following set up of morphism of stacks,
\[\begin{tikzcd}
& & \underline{M} \arrow[dd, "r"] & & \\
& & & & \\
\underline{X} \arrow[rr, "p"] & & \mc{D} \arrow[rr, "F"] & & \mc{C}
\end{tikzcd}.\]
This gives a morphism of stacks $F\circ r\colon \underline{M}\ra \mc{C}$. As $F\circ p=q\colon \underline{X}\ra \mc{C}$ is an atlas for $\mc{C}$, the $2$-fiber product $\underline{X}\times_{\mc{C}}\underline{M}$ in the following diagram, \begin{equation}\label{Diagram:pr2XMtoM}
\begin{tikzcd}
\underline{X}\times_{\mc{C}}\underline{M} \arrow[dd, "pr_1"'] \arrow[rr, "pr_2"] & & \underline{M} \arrow[dd, "F\circ r"] \\
& & \\
\underline{X} \arrow[Rightarrow, shorten >=30pt, shorten <=30pt, uurr] \arrow[rr, "F\circ p=q"] & & \mc{C} 
\end{tikzcd}.\end{equation} is representable by a manifold and the projection map $pr_2\colon \underline{X}\times_{\mc{C}}\underline{M}\ra \underline{M}$ is a surjective submersion at the level of manifolds.     

Consider the morphism of stacks $(p,r)\colon \underline{X}\times_{\mc{C}}\underline{M}\ra \mc{D}\times_{\mc{C}}\mc{D}$ given by (at the level of objects)
\begin{equation}\label{Equation:definitionof(p,r)}
\big(x,m,\alpha\colon (F\circ p)(x)\ra (F\circ r)(m)\big)\mapsto\big (p(x),r(m),\alpha\colon F(p(x)\ra F(r(m))\big).\end{equation}
As the diagonal morphism $\Delta_F\colon \mc{D}\ra \mc{D}\times_{\mc{C}}\mc{D}$ is an epimorphism of stacks, for the morphism of stacks $(p,r)\colon \underline{X}\times_{\mc{C}}\underline{M}\ra \mc{D}\times_{\mc{C}}\mc{D}$, there exists a surjective submersion $\Phi\colon \underline{W}\ra \underline{X}\times_{\mc{C}}\underline{M}$ and a morphism of stacks $\gamma\colon \underline{W}\ra \mc{D}$ producing the following $2$-commutative diagram,
\begin{equation}\label{Diagram:PhiWtoXCM}
\begin{tikzcd}
\underline{W} \arrow[dd, "\gamma"'] \arrow[rr, "\Phi"] & & \underline{X}\times_{\mc{C}}\underline{M} \arrow[dd, "{(p,r)}"] \\
& & \\
\mc{D} \arrow[Rightarrow, shorten >=30pt, shorten <=30pt, uurr] \arrow[rr, "\Delta_F"] & & \mc{D}\times_{\mc{C}}\mc{D}
\end{tikzcd}.\end{equation}
Extending the Diagram \ref{Diagram:PhiWtoXCM} along the first projections $pr_1\colon \underline{X}\times_{\mc{C}}\underline{M}\ra \underline{M}$ and 
$pr_1\colon \mc{D}\times_{\mc{C}}\mc{D}\ra \mc{D}$
we obtain the following diagram,
\begin{equation}\label{Diagram:extendalongpr1}
\begin{tikzcd}
\underline{W} 
\arrow[dd, "\gamma"'] \arrow[rr, "\Phi"] & & \underline{X}\times_{\mc{C}}\underline{M}
\arrow[rr, "pr_1"] \arrow[dd, "{(p,r)}"] & & \underline{X} \arrow[dd, "p"] \\
& & & & \\
\mc{D} \arrow[Rightarrow, shorten >=30pt, shorten <=30pt, uurr] \arrow[rr, "\Delta_F"] & & \mc{D}\times_{\mc{C}}\mc{D} \arrow[Rightarrow, shorten >=30pt, shorten <=30pt, uurr]
\arrow[rr, "pr_1"] & & \mc{D} 
\end{tikzcd}.
\end{equation}
Similarly, extending the Diagram \ref{Diagram:PhiWtoXCM} along the second projections $pr_2\colon \underline{X}\times_{\mc{C}}\underline{M}\ra \underline{M}$ and 
$pr_2\colon \mc{D}\times_{\mc{C}}\mc{D}\ra \mc{D}$
we obtain the following diagram,
\begin{equation}\label{Diagram:extendalongpr2}
\begin{tikzcd}
\underline{W} \arrow[dd, "\gamma"'] \arrow[rr, "\Phi"] & & \underline{X}\times_{\mc{C}}\underline{M} \arrow[rr, "pr_2"] \arrow[dd, "{(p,r)}"] & & \underline{M} \arrow[dd, "r"] \\
& & & & \\
\mc{D} \arrow[Rightarrow, shorten >=30pt, shorten <=30pt, uurr] \arrow[rr, "\Delta_F"] & & \mc{D}\times_{\mc{C}}\mc{D} \arrow[Rightarrow, shorten >=30pt, shorten <=30pt, uurr] \arrow[rr, "pr_2"] & & \mc{D} 
\end{tikzcd}.
\end{equation}
For the sake of convenience we combine the Diagrams \ref{Diagram:extendalongpr1} and \ref{Diagram:extendalongpr2} to draw the following diagram,
\begin{equation}\label{Diagram:WtoXMtoMandX}\begin{tikzcd}
\underline{W} \arrow[dd, "\gamma"'] \arrow[rr, "\Phi"] & & \underline{X}\times_{\mc{C}}\underline{M} \arrow[dd, "{(p,r)}"] \arrow[rr, "pr_2"] \arrow[rrrd,"pr_1"'] & & \underline{M} \arrow[dd, "r"] & \\
& & & & & \underline{X} \arrow[dd, "p"] \\
\mc{D} \arrow[rr, "\Delta_F"] & & \mc{D}\times_{\mc{C}}\mc{D} \arrow[rr, "pr_2"] \arrow[rrrd, "pr_1"'] & & \mc{D} & \\
&& & & & \mc{D} 
\end{tikzcd}.\end{equation}
Observe that the maps $pr_2\colon \underline{X}\times_{\mc{C}}\underline{M}\ra \underline{M}$ (from Diagram \ref{Diagram:pr2XMtoM}) and $\Phi\colon \underline{W}\ra \underline{X}\times_{\mc{C}}\underline{M}$ (from Diagram \ref{Diagram:PhiWtoXCM}) are surjective submersions. Thus, the composition 
$pr_2\circ \Phi\colon \underline{W}\ra \underline{M}$ is a surjective submersion. Consider the composition 
$pr_1\circ \Phi\colon \underline{W}\ra \underline{X}$. This gives the following diagram of morphism of stacks,
\begin{equation}\label{Diagram:CommutativeWMXD}\begin{tikzcd}
\underline{W} \arrow[dd, "pr_1\circ \Phi"'] \arrow[rr, "pr_2\circ \Phi"] & & \underline{M} \arrow[dd, "r"] \\
& & \\
\underline{X} \arrow[rr, "p"] & & \mc{D} 
\end{tikzcd}.\end{equation}
We further note from Diagram \ref{Diagram:WtoXMtoMandX} that, $r\circ pr_2\circ \Phi=pr_2\circ \Delta_F\circ \gamma$ and  
$p\circ pr_1\circ \Phi=pr_1\circ \Delta_F\circ \gamma$. 
As $pr_1\circ \Delta=pr_2\circ \Delta$, we see that 
$pr_1\circ \Delta\circ \gamma=pr_2\circ \Delta\circ \gamma$.
Thus, $r\circ pr_2\circ \Phi=p\circ pr_1\circ \Phi$. So, the Diagram \ref{Diagram:CommutativeWMXD} is a $2$-commutative diagram. Thus, given a morphism of stacks $r\colon \underline{M}\ra \mc{D}$, there exists a surjective submersion $\Gamma=pr_2\circ \Phi\colon \underline{W}\ra \underline{M}$ and a morphism of stacks $\Psi=pr_1\circ \Phi\colon \underline{W}\ra \underline{X}$ with following $2$-commutative diagram,
\begin{equation}\begin{tikzcd}
\underline{W} \arrow[dd, "\Psi"'] \arrow[rr, "\Gamma"] & & \underline{M} \arrow[dd, "r"] \\
& & \\
\underline{X} \arrow[Rightarrow, shorten >=20pt, shorten <=20pt, uurr] \arrow[rr, "p"] & & \mc{D} 
\end{tikzcd}.\end{equation}
Thus, we conclude that the morphism $p\colon \underline{X}\ra\mc{D}$ is an epimorphism of stacks.
\begin{lemma}\label{Lemma:pisanepimorphism} 
	The morphism of stacks $p\colon \underline{X}\ra \mc{D}$ mentioned in Diagram \ref{Diagram:FpToq} is an epimorphism of stacks.
	  \end{lemma}
\subsection{$p\colon \underline{X}\ra \mc{D}$ is an atlas for $\mc{D}$}
\label{Subsection:pIsAnAtlas} Let $F\colon \mc{D}\ra \mc{C}$ be a gerbe over a stack. We further assume that the diagonal morphism $\Delta_F\colon \mc{D}\ra \mc{D}\times_{\mc{C}}\mc{D}$ is a representable surjective submersion. Let $p\colon \underline{X}\ra \mc{D}$ be as in Lemma \ref{Lemma:pisanepimorphism}. In general $p\colon \underline{X}\ra \mc{D}$ is not an atlas for $\mc{D}$. We have proved that $p\colon \underline{X}\ra \mc{D}$ is an epimorphism of stacks. We use the following Proposition to conclude that $p\colon \underline{X}\ra \mc{D}$ is in fact an atlas for $\mc{D}$.
The following proposition is a variant of Proposition $2.16$ in \cite{MR2817778}.
\begin{proposition}\label{Proposition:epimorphismbeinganatlas}
	A morphism of stacks $r\colon \underline{X}\rightarrow \mc{D}$ is an atlas for $\mc{D}$ if
	\begin{enumerate}
		\item the morphism $r\colon \underline{X}\rightarrow \mc{D}$ is an epimorphism of stacks, 
		\item the fibered product $\underline{X}\times_{\mc{D}}\underline{X}$ is representable by a manifold and that the projection maps
		$pr_1\colon X\times_{\mc{D}}X\rightarrow X$ and
		$pr_2\colon X\times_{\mc{D}}X\rightarrow X$ are submersions.
	\end{enumerate}
\end{proposition}
By above Proposition, to prove $p\colon \underline{X}\ra \mc{D}$ is an atlas for $\mc{D}$, it only remains to prove that $\underline{X}\times_{\mc{D}}\underline{X}$ is representable by a manifold and that the projection maps $pr_1\colon X\times_{\mc{D}}X\ra X$ and $pr_2\colon X\times_{\mc{D}}X\ra X$ are submersions. We prove them below.

For $p\colon \underline{X}\ra \mc{D}$ and for $F\circ p\colon \underline{X}\ra \mc{C}$, we have following pull
back diagrams,
\[\begin{tikzcd}
\underline{X}\times_{\mc{D}}\underline{X} \arrow[dd, "pr_1^{\mc{D}}"'] \arrow[rr, "pr_2^{\mc{D}}"] & & \underline{X} \arrow[dd, "p"] & \\
& & & \\
\underline{X} \arrow[Rightarrow, shorten >=30pt, shorten <=30pt, uurr] \arrow[rr, "p"] & & \mc{D} \arrow[rr, "F"] & & \mc{C}
\end{tikzcd} \begin{tikzcd}
\underline{X}\times_{\mc{C}}\underline{X} \arrow[dd, "pr_1^{\mc{C}}"'] \arrow[rr, "pr_2^{\mc{C}}"] & & \underline{X} \arrow[dd, "F\circ p"] \\
& & \\
\underline{X} \arrow[Rightarrow, shorten >=30pt, shorten <=30pt, uurr] \arrow[rr, "F\circ p"] & & \mc{C} 
\end{tikzcd}.\]
By uniqueness of pullback, there exists a unique morphism of stacks 
$\Psi\colon \underline{X}\times_{\mc{D}}\underline{X}\ra 
\underline{X}\times_{\mc{C}}\underline{X}$ with following $2$-commutative diagram,
\begin{equation}\label{Diagram:ExistenceofPsi}\begin{tikzcd}
\underline{X}\times_{\mc{D}}\underline{X} \arrow[rddd, "pr_1^{\mc{D}}"', bend right] \arrow[rrrd, "pr_2^{\mc{D}}", bend left] \arrow[rd, "\Psi"] & & & \\
& \underline{X}\times_{\mc{C}}\underline{X} \arrow[dd, "pr_1^{\mc{C}}"'] \arrow[rr, "pr_2^{\mc{C}}"] & & \underline{X} \arrow[dd, "F\circ p"] \\
& & & \\
& \underline{X} \arrow[rr, "F\circ p"] & & \mc{C} 
\end{tikzcd}.\end{equation}
We have assumed that the diagonal morphism $\Delta_F\colon \mc{D}\ra \mc{D}\times_{\mc{C}}\mc{D}$ is a representable surjective submersion. Consider the morphism of stacks $(p,p)\colon \underline{X}\times_{\mc{C}}\underline{X}\ra \mc{D}\times_{\mc{C}}\mc{D}$ as in equation \ref{Equation:definitionof(p,r)}. We have the following $2$-fiber product diagram,
\begin{equation}\label{Diagram:(p,p)diagonalmorphism}
\begin{tikzcd}
\mc{D}\times_{\mc{D}\times_{\mc{C}}\mc{D}}(\underline{X}\times_{\mc{C}}\underline{X}) \arrow[dd, "pr_1"'] \arrow[rr, "pr_2"] & & \underline{X}\times_{\mc{C}}\underline{X} \arrow[dd, "{(p,p)}"] \\
& & \\
\mc{D} \arrow[Rightarrow, shorten >=30pt, shorten <=30pt, uurr] \arrow[rr, "\Delta_F"] & & \mc{D}\times_{\mc{C}}\mc{D} 
\end{tikzcd}.\end{equation}
We have an isomorphism of stacks 
$\mc{D}\times_{\mc{D}\times_{\mc{C}}\mc{D}}(\underline{X}\times_{\mc{C}}\underline{X})\cong\underline{X}\times_{\mc{D}}\underline{X}$ (\cite[Corollary $69$]{Metzler}). Observe that the morphism of stacks $pr_2\colon \mc{D}\times_{\mc{D}\times_{\mc{C}}\mc{D}}
(\underline{X}\times_{\mc{C}}\underline{X})\ra \underline{X}\times_{\mc{C}}\underline{X}$ in Diagram \ref{Diagram:(p,p)diagonalmorphism} is same as the morphism of stacks $\Psi\colon \underline{X}\times_{\mc{D}}\underline{X}\ra 
\underline{X}\times_{\mc{C}}\underline{X}$ in Diagram \ref{Diagram:ExistenceofPsi}.
Thus, the above $2$-fiber product diagram can be seen as the following $2$-commutative diagram,
\begin{equation}\label{Diagram:Psi(p,p)DeltaF}
\begin{tikzcd}
\underline{X}\times_{\mc{D}}\underline{X} \arrow[dd] \arrow[rr, "\Psi"] & & \underline{X}\times_{\mc{C}}\underline{X} \arrow[dd, "{(p,p)}"] \\
& & \\
\mc{D} \arrow[Rightarrow, shorten >=30pt, shorten <=30pt, uurr] \arrow[rr, "\Delta_F"] & & \mc{D}\times_{\mc{C}}\mc{D} 
\end{tikzcd}.\end{equation} 
As the diagonal morphism $\Delta_F\colon \mc{D}\ra \mc{D}\times_{\mc{C}}\mc{D}$ is a representable surjective submersion, the $2$-fiber product $\mc{D}\times_{\mc{D}\times_{\mc{C}}\mc{D}}(\underline{X}\times_{\mc{C}}\underline{X})$ is representable by a manifold and the projection map $pr_2\colon \mc{D}\times_{\mc{D}\times_{\mc{C}}\mc{D}}(\underline{X}\times_{\mc{C}}\underline{X})\ra \underline{X}\times_{\mc{C}}\underline{X}$ is a surjective submersion at the level of manifolds. Thus, we see that $\underline{X}\times_{\mc{D}}\underline{X}$ is representable by manifold and $\Psi\colon \underline{X}\times_{\mc{D}}\underline{X}\ra 
\underline{X}\times_{\mc{C}}\underline{X}$ is a surjective submersion at the level of manifolds.

Both being compositions of surjective submersions, we see that $pr_1^{\mc{D}}=pr_1^{\mc{C}}\circ \Psi$ and 
$pr_2^{\mc{D}}=pr_2^{\mc{C}}\circ \Psi$ are surjective submersions at the level of manifolds. Thus, $p\colon \underline{X}\ra \mc{D}$ is an atlas for $\mc{D}$. So, we have shown the following:
\begin{lemma}\label{Lemma:extraassumptiononDiagonalmorphism}
	Let $F\colon \mc{D}\ra \mc{C}$ be a gerbe over a stack. Further assume that, the diagonal morphism $\Delta_F\colon \mc{D}\ra \mc{D}\times_{\mc{C}}\mc{D}$ is a representable surjective submersion. Then, the morphism of stacks $p\colon \underline{X}\ra \mc{D}$ mentioned in Diagram \ref{Diagram:FpToq} is an atlas for the stack $\pi_{\mc{D}}\colon \mc{D}\ra \text{Man}$. In particular, there exists an atlas $p\colon \underline{X}\ra \mc{D}$ for $\mc{D}$ and an atlas $q\colon \underline{X}\ra \mc{C}$ with a $2$-commutative diagram as in \ref{Diagram:FpToq}.
	  \end{lemma}

\subsection{A gerbe over a stack gives a Lie groupoid extension}
\label{Subsection:GoSgivingLiegroupoidextension}
Let $F\colon \mc{D}\ra \mc{C}$ be a gerbe over a stack. We further assume that the diagonal morphism $\Delta_F\colon \mc{D}\ra \mc{D}\times_{\mc{C}}\mc{D}$ is a representable surjective submersion. By Lemma \ref{Lemma:extraassumptiononDiagonalmorphism}, there exists an atlas $p\colon \underline{X}\ra \mc{D}$ for $\mc{D}$ and an atlas $q\colon \underline{X}\ra \mc{C}$ with a $2$-commutative diagram as in \ref{Diagram:FpToq}. For atlases $p\colon \underline{X}\ra\mc{D}$ and $q\colon \underline{X}\ra\mc{C}$, we have respectively associated the Lie groupoids $\mc{G}_p=(X\times_{\mc{D}}X\rightrightarrows X)$ and $\mc{G}_q=(X\times_{\mc{C}}X\rightrightarrows Y)$. The morphism of stacks $F\colon \mc{D}\ra \mc{C}$ induces a morphism of stacks $\Psi\colon \underline{X}\times_{\mc{D}}\underline{X}\ra \underline{X}\times_{\mc{C}}\underline{X}$ as in Diagram \ref{Diagram:ExistenceofPsi}.
Explicitly, $\Psi\colon \underline{X}\times_{\mc{D}}\underline{X}\ra \underline{X}\times_{\mc{C}}\underline{X}$ is given at the level of objects by 
\[(y,z,\alpha\colon p(y)\rightarrow p(z))\mapsto (y,z,F(\alpha)\colon F(p(y))\rightarrow F(p(z))).\]
At the level of morphisms, an arrow \[(u,v)\colon \big(y,z,\alpha\colon p(y)\rightarrow p(z))\ra \big(y',z',\alpha'\colon p(y')\rightarrow p(z')\big)\]
in $\underline{X}\times_{\mc{D}}\underline{X}$ is mapped to the arrow \[(u,v)\colon \big(y,z,F(\alpha)\colon F(p(y))\rightarrow F(p(z)))\ra \big(y',z',F(\alpha')\colon F(p(y'))\rightarrow F(p(z'))\big)\]
in $\underline{X}\times_{\mc{C}}\underline{X}$.
This morphism of stacks $\Psi\colon \underline{X}\times_{\mc{D}}\underline{X}\rightarrow
\underline{X}\times_{\mc{C}}\underline{X}$ is compatible with projection maps $pr_1,pr_2\colon \underline{X}\times_{\mc{D}}\underline{X}
\rightarrow \underline{X}$ and $pr_1,pr_2\colon \underline{X}\times_{\mc{C}}\underline{X}
\rightarrow \underline{X}$ in the sense that following diagram is a commutative diagram of morphisms of stacks, 
\begin{equation}\label{Diagram:PsiXDXtoXCX}\begin{tikzcd}
\underline{X}\times_{\mc{D}}\underline{X} \arrow[dd,xshift=0.75ex,"pr_2"]\arrow[dd,xshift=-0.75ex,"pr_1"'] \arrow[rr,"\Psi"] & & \underline{X}\times_{\mc{C}}\underline{X} \arrow[dd,xshift=0.75ex,"pr_2"]\arrow[dd,xshift=-0.75ex,"pr_1"'] \\
& & \\
\underline{X} \arrow[rr,"\text{Id}"] & & \underline{X}
\end{tikzcd}.\end{equation}
Recall from  the Section \ref{Subsection:Liegroupoidassociatedtoatlas} that, the morphisms of stacks $pr_1\colon \underline{X}\times_{\mc{D}}\underline{X}\ra \underline{X}$ and $pr_2\colon \underline{X}\times_{\mc{D}}\underline{X}\ra \underline{X}$ respectively corresponds to the source and target maps of the Lie groupoid $\mc{G}_p=(X\times_{\mc{D}}X\rightrightarrows X)$. Likewise for the Lie groupoid $\mc{G}_q=(X\times_{\mc{C}}X\rightrightarrows X)$. Let $\Theta\colon X\times_{\mc{D}}X\ra X\times_{\mc{C}}X$ be the map of manifolds associated to the morphism of stacks $\Psi\colon \underline{X}\times_{\mc{D}}\underline{X}\ra \underline{X}\times_{\mc{C}}\underline{X}$.
Then, the diagram \ref{Diagram:PsiXDXtoXCX}, gives the following diagram of morphism of Lie groupoids,
\begin{equation}\label{Diagram:ThetaXDXtoXCX}\begin{tikzcd}
X\times_{\mc{D}}X \arrow[dd,xshift=0.75ex,"t"]\arrow[dd,xshift=-0.75ex,"s"'] \arrow[rr,"\Theta"] & & X\times_{\mc{C}}X \arrow[dd,xshift=0.75ex,"t"]\arrow[dd,xshift=-0.75ex,"s"'] \\
& & \\
X \arrow[rr,"\text{Id}"] & & X
\end{tikzcd}.\end{equation}
We have the following result.
\begin{proposition}\label{Proposition:GoSgivingMorphismofLiegroupoids}
	Let $F\colon\mc{D}\ra \mc{C}$ be a gerbe over a stack. Assume further that the diagonal morphism $\Delta_F\colon \mc{D}\rightarrow \mc{D}\times_{\mc{C}}\mc{D}$ is a representable surjective submersion. Then the morphism of stacks $F\colon \mc{D}\rightarrow \mc{C}$ gives a morphism of Lie groupoids \[\mc{G}_p\ra \mc{G}_q\colon (X\times_{\mc{D}}X\rightrightarrows X)\rightarrow(X\times_{\mc{C}}X\rightrightarrows X),\]
	where $p\colon \underline{X}\ra \mc{D}$ and $q\colon \underline{X}\ra \mc{C}$ are as in Lemma \ref{Lemma:extraassumptiononDiagonalmorphism}. 
\end{proposition}
Observe that the morphism of stacks $\Psi\colon \underline{X}\times_{\mc{D}}\underline{X}
\rightarrow
\underline{X}\times_{\mc{C}}\underline{X}$ is a surjective submersion at the level of manifolds (discussion after Diagram \ref{Diagram:Psi(p,p)DeltaF}); that is, the map $\Theta\colon X\times_{\mc{D}}X\ra X\times_{\mc{C}}X$ is a surjective submersion. Thus, $(\Theta,\text{Id})\colon (X\times_{\mc{D}}X\rightrightarrows X)\rightarrow (X\times_{\mc{C}}X\rightrightarrows X)$ (Diagram \ref{Diagram:ThetaXDXtoXCX}) is a Lie groupoid extension.
\begin{lemma}\label{Lemma:MoSgivingLiegroupoidExtension}
	Let $F\colon \mc{D}\ra \mc{C}$ be a gerbe over a stack. Assume  that the diagonal morphism $\Delta_F\colon \mc{D}\rightarrow \mc{D}\times_{\mc{C}}\mc{D}$ is a representable surjective submersion. Then the morphism of stacks $F\colon \mc{D}\rightarrow \mc{C}$ gives a Lie groupoid extension \[\mc{G}_p\ra \mc{G}_q\colon (X\times_{\mc{D}}X\rightrightarrows X)\rightarrow(X\times_{\mc{C}}X\rightrightarrows X)\]
	where $p\colon \underline{X}\ra \mc{D}$ and $q\colon \underline{X}\ra \mc{C}$ are as in Lemma \ref{Lemma:extraassumptiononDiagonalmorphism}.
	  \end{lemma}
\subsection{Uniqueness of a Lie groupoid extension associated to a gerbe over a stack}\label{Subsection:uniquenessOfExtension} Let $F\colon \mc{D}\rightarrow \mc{C}$ be a gerbe over a stack. We further assume that $\Delta_F\colon \mc{D}\rightarrow \mc{D}\times_{\mc{C}}\mc{D}$ is a representable surjective submersion. For atlases $p\colon \underline{X}\ra \mc{D}$ and $q\colon \underline{X}\ra \mc{C}$ mentioned in Lemma \ref{Lemma:extraassumptiononDiagonalmorphism}, we have assigned a Lie groupoid extension 
\[ (X\times_{\mc{D}}X\rightrightarrows X)\rightarrow(X\times_{\mc{C}}X\rightrightarrows X)\] in Lemma \ref{Lemma:MoSgivingLiegroupoidExtension}. We prove in this subsection that,  up to a Morita equivalence, this Lie groupoid extension does not depend on the choice of atlases.

Let $q_Y\colon \underline{Y}\rightarrow \mc{C}$ be another atlas for $\mc{C}$ and $p_Y\colon \underline{Y}\rightarrow \mc{D}$ be the corresponding atlas for $\mc{D}$ as in Lemma \ref{Lemma:extraassumptiononDiagonalmorphism}. These atlases $p_Y\colon \underline{Y}\ra \mc{D},q_Y\colon \underline{Y}\ra\mc{C}$ gives a Lie groupoid extension \[(Y\times_{\mc{D}}Y\rightrightarrows Y)\rightarrow 
(Y\times_{\mc{C}}Y\rightrightarrows Y)\] as in Lemma \ref{Lemma:MoSgivingLiegroupoidExtension}.

We prove that 
$(X\times_{\mc{D}}X\rightrightarrows X)\rightarrow 
(X\times_{\mc{C}}X\rightrightarrows X)$ and $(Y\times_{\mc{D}}Y\rightrightarrows Y)\rightarrow 
(Y\times_{\mc{C}}Y\rightrightarrows Y)$ are Morita equivalent Lie groupoid extensions. 
Recall by Definition \ref{Definition:MoritamorphismofLiegroupoidExtension}, it means that  there exists a Lie groupoid extension $(\mc{G}_1\rightrightarrows \mc{G}_0)\rightarrow 
(\mc{H}_1\rightrightarrows \mc{H}_0)$ and a pair of Morita morphisms of Lie groupoid extensions 
\begin{align*}
\big((\mc{G}_1\rightrightarrows \mc{G}_0)\rightarrow (\mc{H}_1\rightrightarrows \mc{G}_0)\big)&\ra \big((X\times_{\mc{D}}X\rightrightarrows X)\rightarrow 
(X\times_{\mc{C}}X\rightrightarrows X)\big),\\
\big((\mc{G}_1\rightrightarrows \mc{G}_0)\rightarrow (\mc{H}_1\rightrightarrows \mc{G}_0)\big)&\ra \big((Y\times_{\mc{D}}Y\rightrightarrows Y)\rightarrow 
(Y\times_{\mc{C}}Y\rightrightarrows Y)\big).
\end{align*}

Diagrammatically that means we have to find a Lie groupoid extension and a pair of Morita morphisms of Lie groupoid extensions, 
\begin{equation}\text{ from }\begin{tikzcd}
\mc{G}_1 \arrow[dd,xshift=0.75ex,"t"]\arrow[dd,xshift=-0.75ex,"s"'] \arrow[rr] & & \mc{H}_1 \arrow[dd,xshift=0.75ex,"t"]\arrow[dd,xshift=-0.75ex,"s"'] \\
& & \\
\mc{G}_0 \arrow[rr] & & \mc{G}_0
\end{tikzcd} ~ \text{to} \begin{tikzcd}
X\times_{\mc{D}}X \arrow[dd,xshift=0.75ex,"t"]\arrow[dd,xshift=-0.75ex,"s"'] \arrow[rr] & & X\times_{\mc{C}}X \arrow[dd,xshift=0.75ex,"t"]\arrow[dd,xshift=-0.75ex,"s"'] \\
& & \\
X \arrow[rr] & & X
\end{tikzcd}\end{equation}
and 
\begin{equation} \text{ from }\begin{tikzcd}
\mc{G}_1 \arrow[dd,xshift=0.75ex,"t"]\arrow[dd,xshift=-0.75ex,"s"'] \arrow[rr] & & \mc{H}_1 \arrow[dd,xshift=0.75ex,"t"]\arrow[dd,xshift=-0.75ex,"s"'] \\
& & \\
\mc{G}_0 \arrow[rr] & & \mc{G}_0
\end{tikzcd} ~ \text{to} \begin{tikzcd}
Y\times_{\mc{D}}Y \arrow[dd,xshift=0.75ex,"t"]\arrow[dd,xshift=-0.75ex,"s"'] \arrow[rr] & & Y\times_{\mc{C}}Y \arrow[dd,xshift=0.75ex,"t"]\arrow[dd,xshift=-0.75ex,"s"'] \\
& & \\
Y \arrow[rr] & & Y
\end{tikzcd}.\end{equation}
For this, first we need (Section \ref{Subsection:MoritamorphismofGroupoidExtensions}) a smooth manifold $\mc{G}_0$ and a pair of smooth maps $\mc{G}_0\rightarrow X, \mc{G}_0\rightarrow Y$. 

For the morphisms of stacks $p_X\colon \underline{X}\rightarrow \mc{D}$ and $p_Y\colon \underline{Y}\rightarrow \mc{D}$, consider the following $2$-fiber product diagram,
\[\begin{tikzcd}
\underline{X}\times_{\mc{D}}\underline{Y} \arrow[dd, "pr_1"'] \arrow[rr, "pr_2"] & & \underline{Y} \arrow[dd, "p_Y"] \\
& & \\
\underline{X} \arrow[Rightarrow, shorten >=20pt, shorten <=20pt, uurr] \arrow[rr, "p_X"'] & & \mc{D} 
\end{tikzcd}.\]

As $p_X\colon \underline{X}\rightarrow \mc{D}$ is an atlas, the morphism of stacks $pr_2\colon \underline{X}\times_{\mc{D}}\underline{Y}\rightarrow \underline{Y}$ induces a surjective submersion, 
$g\colon X\times_{\mc{D}}Y\rightarrow Y$ at the level of manifolds.
Similarly 
$pr_1\colon \underline{X}\times_{\mc{D}}\underline{Y}\rightarrow \underline{X}$ induces a surjective submersion $f\colon X\times_{\mc{D}}Y\rightarrow X$ at the level of manifolds. 

We now construct a Lie groupoid extension 
of the form \[(**\rightrightarrows X\times_{\mc{D}}Y)\rightarrow (***\rightrightarrows X\times_{\mc{D}}Y).\]
Here for the time being we denote the respective morphism sets by $**$ and $***$.

Next, we find a Morita morphisms of Lie groupoid extensions, 
\begin{equation}\text{ from }
\begin{tikzcd}
** \arrow[dd,xshift=0.75ex,"t"]\arrow[dd,xshift=-0.75ex,"s"'] \arrow[rr] & & ***\arrow[dd,xshift=0.75ex,"t"]\arrow[dd,xshift=-0.75ex,"s"'] \\
& & \\
X\times_{\mc{D}}Y \arrow[rr] & & X\times_{\mc{D}}Y
\end{tikzcd}\text{to} \begin{tikzcd}
X\times_{\mc{D}}X \arrow[dd,xshift=0.75ex,"t"]\arrow[dd,xshift=-0.75ex,"s"'] \arrow[rr] & & X\times_{\mc{C}}X \arrow[dd,xshift=0.75ex,"t"]\arrow[dd,xshift=-0.75ex,"s"'] \\
& & \\
X \arrow[rr] & & X
\end{tikzcd} 
\end{equation} 
and 
\begin{equation} \text{ from } \begin{tikzcd}
** \arrow[dd,xshift=0.75ex,"t"]\arrow[dd,xshift=-0.75ex,"s"'] \arrow[rr] & & *** \arrow[dd,xshift=0.75ex,"t"]\arrow[dd,xshift=-0.75ex,"s"'] \\
& & \\
X\times_{\mc{D}}Y \arrow[rr] & & X\times_{\mc{D}}Y
\end{tikzcd}\text{to} \begin{tikzcd}
Y\times_{\mc{D}}Y \arrow[dd,xshift=0.75ex,"t"]\arrow[dd,xshift=-0.75ex,"s"'] \arrow[rr] & & Y\times_{\mc{C}}Y \arrow[dd,xshift=0.75ex,"t"]\arrow[dd,xshift=-0.75ex,"s"'] \\
& & \\
Y \arrow[rr] & & Y
\end{tikzcd}\end{equation} 
That means, as per Definition \ref{Definition:MoritamorphismofLiegroupoidExtension}, we need a pair of Morita morphisms of Lie groupoids 
\begin{equation} \begin{tikzcd}
** \arrow[dd,xshift=0.75ex,"t"]\arrow[dd,xshift=-0.75ex,"s"'] \arrow[rr] & & X\times_{\mc{D}}X \arrow[dd,xshift=0.75ex,"t"]\arrow[dd,xshift=-0.75ex,"s"'] \\
& & \\
X\times_{\mc{D}}Y \arrow[rr,"f"] & & X
\end{tikzcd}\end{equation} and \begin{equation}\begin{tikzcd}
*** \arrow[dd,xshift=0.75ex,"t"]\arrow[dd,xshift=-0.75ex,"s"'] \arrow[rr] & & X\times_{\mc{C}}X \arrow[dd,xshift=0.75ex,"t"]\arrow[dd,xshift=-0.75ex,"s"'] \\
& & \\
X\times_{\mc{D}}Y \arrow[rr,"f"] & & X
\end{tikzcd}\end{equation}
which are compatible with maps $**\rightarrow ***$ and $ X\times_{\mc{D}}X\rightarrow X\times_{\mc{C}}X$. 
Similarly, we need a pair of Morita morphisms of Lie groupoids 
\begin{equation}
\begin{tikzcd}
** \arrow[dd,xshift=0.75ex,"t"]\arrow[dd,xshift=-0.75ex,"s"'] \arrow[rr] & & Y\times_{\mc{D}}Y \arrow[dd,xshift=0.75ex,"t"]\arrow[dd,xshift=-0.75ex,"s"'] \\
& & \\
X\times_{\mc{D}}Y \arrow[rr,"g"] & & Y
\end{tikzcd}\end{equation}
and
\begin{equation}
\begin{tikzcd}
*** \arrow[dd,xshift=0.75ex,"t"]\arrow[dd,xshift=-0.75ex,"s"'] \arrow[rr] & & Y\times_{\mc{C}}Y \arrow[dd,xshift=0.75ex,"t"]\arrow[dd,xshift=-0.75ex,"s"'] \\
& & \\
X\times_{\mc{D}}Y \arrow[rr,"g"] & & Y
\end{tikzcd}\end{equation} 
compatible with maps $**\rightarrow ***$ and $ Y\times_{\mc{D}}Y\rightarrow Y\times_{\mc{C}}Y$.

Our task is to find $**$ and $***$. 
Recalling the set up of Morita morphisms of Lie groupoids (Definition \ref{Definition:MoritamorphismofLiegroupoids}), we see that given a surjective submersion $f\colon M\rightarrow N$ and a Lie groupoid $\mc{G}_1\rightrightarrows N$, the pullback Lie groupoid (Section \ref{Subsection:PullbackLiegroupoid}) as in below diagram,
\begin{equation} 
\begin{tikzcd}
f^*\mc{G}_1 \arrow[dd,xshift=0.75ex,"t"]\arrow[dd,xshift=-0.75ex,"s"'] \arrow[rr] & & \mc{G}_1 \arrow[dd,xshift=0.75ex,"t"]\arrow[dd,xshift=-0.75ex,"s"'] \\
& & \\
M \arrow[rr,"f"] & & N
\end{tikzcd} \end{equation}
gives a Morita morphism of Lie groupoids $(f^*\mc{G}_1\rightrightarrows M)\ra (\mc{G}_1\rightrightarrows N)$.

Let $(\mc{G}_1 \rightrightarrows X\times_{\mc{D}}Y)$ be the pullback of the Lie groupoid $(X\times_{\mc{D}}X\rightrightarrows X)$ along $f\colon X\times_{\mc{D}}Y\rightarrow X$ and $(\mc{G}_1'\rightrightarrows X\times_{\mc{D}}Y) $ be the pullback of the Lie groupoid $(Y\times_{\mc{D}}Y\rightrightarrows Y)$ along $g\colon X\times_{\mc{D}}Y\rightarrow Y$. 
However, as we will shortly see, we do not have to distinguish between these two pullbacks as they are isomorphic.
We have following diagrams representing the pullback groupoids,
\begin{equation} \begin{tikzcd}
\mc{G}_1 \arrow[dd,xshift=0.75ex,"t"]\arrow[dd,xshift=-0.75ex,"s"'] \arrow[rr] & & X\times_{\mc{D}}X \arrow[dd,xshift=0.75ex,"t"]\arrow[dd,xshift=-0.75ex,"s"'] \\
& & \\
X\times_{\mc{D}}Y \arrow[rr,"f"] & & X
\end{tikzcd}
\begin{tikzcd}
\mc{G}_1' \arrow[dd,xshift=0.75ex,"t"]\arrow[dd,xshift=-0.75ex,"s"'] \arrow[rr] & & Y\times_{\mc{D}}Y \arrow[dd,xshift=0.75ex,"t"]\arrow[dd,xshift=-0.75ex,"s"'] \\
& & \\
X\times_{\mc{D}}Y \arrow[rr,"g"] & & Y
\end{tikzcd}.
\end{equation} 
Similarly, we write $(\mc{H}_1 \rightrightarrows X\times_{\mc{D}}Y) $ for pullback of the Lie groupoid $(X\times_{\mc{C}}X\rightrightarrows X)$ along $f\colon X\times_{\mc{D}}Y\rightarrow X$ and $(\mc{H}_1'\rightrightarrows X\times_{\mc{D}}Y) $ for pullback of the Lie groupoid $(Y\times_{\mc{C}}Y\rightrightarrows Y)$ along $g\colon X\times_{\mc{D}}Y\rightarrow Y$. As before, the Lie groupoids 
$\mc{H}_1 \rightrightarrows X\times_{\mc{D}}Y $ and $\mc{H}_1'\rightrightarrows X\times_{\mc{D}}Y $ will be isomorphic. We have following diagrams representing the pullback groupoids,
\begin{equation} \begin{tikzcd}
\mc{H}_1 \arrow[dd,xshift=0.75ex,"t"]\arrow[dd,xshift=-0.75ex,"s"'] \arrow[rr] & & X\times_{\mc{C}}X \arrow[dd,xshift=0.75ex,"t"]\arrow[dd,xshift=-0.75ex,"s"'] \\
& & \\
X\times_{\mc{D}}Y \arrow[rr,"f"] & & X
\end{tikzcd}
\begin{tikzcd}
\mc{H}_1' \arrow[dd,xshift=0.75ex,"t"]\arrow[dd,xshift=-0.75ex,"s"'] \arrow[rr] & & Y\times_{\mc{C}}Y \arrow[dd,xshift=0.75ex,"t"]\arrow[dd,xshift=-0.75ex,"s"'] \\
& & \\
X\times_{\mc{D}}Y \arrow[rr,"g"] & & Y
\end{tikzcd}.
\end{equation} 
Now we give an isomorphism between $\mc{G}_1\rightrightarrows X\times_{\mc{D}}Y$ and $\mc{G}_1'\rightrightarrows X\times_{\mc{D}}Y$. The construction of isomorphism between $\mc{H}_1\rightrightarrows X\times_{\mc{D}}Y$ and $\mc{H}_1'\rightrightarrows X\times_{\mc{D}}Y$ is very much same. We define map $\mc{G}_1\rightarrow \mc{G}_1'$ by giving a morphism of stacks $\underline{\mc{G}_1}\rightarrow \underline{\mc{G}_1'}$, where
\[\underline{\mc{G}_1}=(\underline{X}\times_{\mc{D}}\underline{Y})\times_{f,X,s}(\underline{X}\times_{\mc{D}}\underline{X})\times_{f,X,t}(\underline{X}\times_{\mc{D}}\underline{Y})\]
and 
\[\underline{\mc{G}_1'}=(\underline{X}\times_{\mc{D}}\underline{Y})\times_{g,X,s}(\underline{Y}\times_{\mc{D}}\underline{Y})\times_{g,X,t}(\underline{X}\times_{\mc{D}}\underline{Y}\big).\]
A typical element in the object set of $\underline{\mc{G}_1}$ is of the form 
\[\bigg(\big(m,n,\alpha\colon p(m)\rightarrow q(n)\big),\big(a,b,p(a)\rightarrow p(b)\big),\big(m',n',\alpha'\colon p(m')\rightarrow q(n')\big)\bigg)\]
such that $m=s\big(a,b,p(a)\rightarrow p(b)\big)=a$, $n=t(a,b,p(a)\rightarrow p(b))=b$ and $p(a)\rightarrow p(b)$ is just $p(m)\rightarrow p(m')$. So, this demands a typical element to be of the form 
\[\bigg(\big(m,n,\alpha\colon p(m)\rightarrow q(n)\big),\big(m,m',p(m)\rightarrow p(m')\big),\big(m',n',\alpha'\colon p(m')\rightarrow q(n')\big)\bigg).\]
The corresponding image in $\underline{\mc{G}_1}'$ is 
\[\bigg(\big(m,n,\alpha\colon p(m)\rightarrow q(n)\big),\big(n,n',q(n)\rightarrow q(n')\big),\big(m',n',\alpha'\colon p(m')\rightarrow q(n')\big)\bigg)\]
This gives a map of stacks $\underline{\mc{G}_1}\rightarrow \underline{\mc{G}_1'}$ at the level of objects. The map at the level of morphisms can be defined similarly. This gives an isomorphism of stacks $\underline{\mc{G}_1}\rightarrow \underline{\mc{G}_1'}$, which in turn induces an isomorphism of Lie groupoids $\mc{G}_1\rightrightarrows X\times_{\mc{D}}Y$ and $\mc{G}_1'\rightrightarrows X\times_{\mc{D}}Y$. Hence, the pullbacks are isomorphic. 
\begin{lemma}\label{Lemma:LiegroupoidExtnIsIndependentofq}
	Let $F\colon \mc{D}\rightarrow \mc{C}$ be a gerbe over a stack. 
	Assume that the diagonal morphism $\Delta_F\colon \mc{D} \rightarrow \mc{D} \times_{\mc{C}}\mc{D}$ is a representable surjective submersion. 
	Then, upto a Morita equivalence, the Lie groupoid extension in Lemma \ref{Lemma:MoSgivingLiegroupoidExtension} does not depend on the choice of $q\colon \underline{X}\rightarrow \mc{C}$.
	  \end{lemma}
Thus, using Lemmas \ref{Lemma:Liegroupoidrepresentingstack}, \ref{Lemma:atlasindependent}, \ref{Lemma:existenceofatlasforC}, \ref{Lemma:pisanepimorphism}, \ref{Lemma:extraassumptiononDiagonalmorphism}, \ref{Lemma:MoSgivingLiegroupoidExtension} and \ref{Lemma:LiegroupoidExtnIsIndependentofq} we have the following result.
\begin{theorem}\label{Theorem:GoSgivesLiegroupoidextension}
	Let $F\colon \mc{D}\rightarrow \mc{C}$ be a gerbe over a stack. Assume that the diagonal morphism $\Delta_F\colon \mc{D} \rightarrow \mc{D} \times_{\mc{C}}\mc{D}$ is a representable surjective submersion. Then there exists an atlas $p\colon \underline{X}\ra \mc{D}$ for $\mc{D}$ and an atlas $q\colon \underline{X}\ra \mc{C}$, as in Lemma \ref{Lemma:extraassumptiononDiagonalmorphism}, producing a Lie groupoid extension $\phi\colon \mc{G}\rightarrow \mc{H}$, where 
	$\mc{G}=(X\times_{\mc{D}}X\rightrightarrows X)$ and 
	$\mc{H}=(X\times_{\mc{C}}X\rightrightarrows X)$. Explicitly, the morphism of stacks $\Phi\colon B\mc{G}\ra B\mc{H}$ associated to $\phi\colon \mc{G}\ra \mc{H}$ (Lemma \ref{Construction:morphismofLiegroupoidstoMorphismofstacks}) along with the morphism of stacks $F\colon \mc{D}\ra \mc{H}$ forms following $2$-commutative diagram,
	\[\begin{tikzcd}
	B\mc{G} \arrow[dd, "\cong"'] \arrow[rr,"\Phi"] & & B\mc{H} \arrow[dd, "\cong"] \\
	& & \\
	\mc{D} \arrow[Rightarrow, shorten >=20pt, shorten <=20pt, uurr] \arrow[rr, "F"] & & \mc{C} 
	\end{tikzcd}.\]
	Here the isomorphisms $\mc{D}\cong B\mc{G}$ and $\mc{C}\cong B\mc{H}$ are as mentioned in Lemma \ref{Lemma:Liegroupoidrepresentingstack}. Further, if there exists another gerbe over the stack $F'\colon \mc{D}'\ra \mc{C}'$  isomorphic to the gerbe $F\colon \mc{D}\ra \mc{C}$, then the Lie groupoid extensions associated to $F'\colon \mc{D}'\ra \mc{C}'$ and $F\colon \mc{D}\ra \mc{C}$ are Morita equivalent. 
\end{theorem} 
\begin{remark}\label{Remark:FullCapacityofDeltaF}
	Observe that we have not made full use of the condition $\Delta_F$ being a surjective submersion. We have only used the following. The morphism of stacks $p\colon \underline{X}\ra \mc{D}$ obtained in Lemma \ref{Lemma:extraassumptiononDiagonalmorphism} is such that, $\underline{X}\times_{\mc{D}}\underline{X}$ is representable by a manifold and the morphism of stacks 
	$\Psi\colon \underline{X}\times_{\mc{D}}\underline{X}
	\rightarrow
	\underline{X}\times_{\mc{C}}\underline{X}$ is a surjective submersion at the level of manifolds. 
	  \end{remark}

\section{A Gerbe over a stack associated to a Lie groupoid extension}
\label{Section:GoSassociatedtoALiegrupoidExtension}
In this section, we describe the construction of a gerbe over a stack from a given Lie groupoid extension. 

Outline of this section is as follows: 
\begin{enumerate}
	\item Given a morphism of Lie groupoids $\phi\colon \mc{G}\rightarrow \mc{H}$,
	we associate a morphism of stacks $F\colon B\mc{G}\rightarrow B\mc{H}$ (Lemma \ref{Construction:morphismofLiegroupoidstoMorphismofstacks}).
	\item If the morphism of Lie groupoids $\phi\colon \mc{G}\rightarrow \mc{H}$ in step $(1)$ is a Lie groupoid extension, 
	then we prove that the associated morphism of stacks $F\colon B\mc{G}\rightarrow B\mc{H}$ is a gerbe over a stack (Theorem \ref{Theorem:BGBHisaGoS}).
\end{enumerate}

\subsection{A Morphism of stacks associated to a Morphism of Lie groupoids}
\label{Subsection:MorphismofLGgivingMorphismofStacks} Given a morphism of Lie groupoids $\phi\colon \mc{G}\rightarrow \mc{H}$, we associate a morphism of stacks $F\colon B\mc{G}\rightarrow B\mc{H}$ in two steps: 
\begin{enumerate}
	\item Given a morphism of Lie groupoids $\phi\colon \mc{G}\rightarrow \mc{H}$, we associate a $\mc{G}-\mc{H}$ bibundle $\left<\phi\right>\colon \mc{G}\ra \mc{H}$ (Remark $3.24$ and Remark $3.27$ in \cite{MR2778793}).
	\item Given a $\mc{G}-\mc{H}$ bibundle $P\colon \mc{G}\ra \mc{H}$, we associate a morphism of stacks $BP\colon B\mc{G}\rightarrow B\mc{H}$ (Remark $3.30$ and Section $4$ in \cite{MR2778793}).
\end{enumerate}

\subsubsection{A morphism of Lie groupoids $\mc{G}\ra \mc{H}$ gives a $\mc{G}-\mc{H}$ bibundle}\label{SubSubsection:bibundleAssociatedtoMorphismofLiegroupoids} Given a morphism of Lie groupoids $\phi\colon \mc{G}\rightarrow \mc{H}$, 
we associate a $\mc{G}-\mc{H}$ bibundle $\left<\phi\right>\colon \mc{G}\ra \mc{H}$. 

Recall that, 
for a Lie groupoid $\mc{H}$, the target map $t\colon \mc{H}_1\ra \mc{H}_0$ is a principal $\mc{H}$-bundle (Example \ref{Example:tG1G0isprincipalbundle}). 
Consider the pullback of the principal $\mc{H}$-bundle $t\colon \mc{H}_1\rightarrow \mc{H}_0$ along the map $\phi_0\colon \mc{G}_0\rightarrow \mc{H}_0$ to get a principal $\mc{H}$-bundle over $\mc{G}_0$ (Section \ref{Subsection:PullbackofPrincipalLiegroupoidbundle}); as explained in the diagram below:
\begin{equation} \begin{tikzcd}
& \mc{G}_0\times_{\mc{H}_0}\mc{H}_1 \arrow[ld, "pr_1"'] \arrow[rd, "pr_2"] & & \mc{H}_1 \arrow[dd,xshift=0.75ex,"t"]
\arrow[dd,xshift=-0.75ex,"s"'] \\
\mc{G}_0 \arrow[rd, "\phi_0"] & & \mc{H}_1 \arrow[ld, "t"'] \arrow[rd, "s"] & \\
& \mc{H}_0 & & \mc{H}_0
\end{tikzcd}.\end{equation}
The maps $\mu\colon (\mc{G}_0\times_{\mc{H}_0}\mc{H}_1)
\times_{s\circ pr_2,\mc{H}_0,t}\mc{H}_1\ra \mc{G}_0\times_{\mc{H}_0}\mc{H}_1$, $((u,h),\tilde{h})\mapsto (u,h\circ \tilde{h})$ and $\tilde{\mu}\colon \mc{G}_1\times_{s,\mc{G}_0,pr_1}
(\mc{G}_0\times_{\mc{H}_0}\mc{H}_1) \ra(\mc{G}_0\times_{\mc{H}_0}\mc{H}_1)$,
$(g,(u,h))\mapsto (t(g),\phi(g)\circ h)$ respectively give \textit{a right action} of $\mc{H}$ on 
$\mc{G}_0\times_{\mc{H}_0}\mc{H}_1$ and left action of $\mc{G}$ on $\mc{G}_0\times_{\mc{H}_0}\mc{H}_1$.
Thus, the manifold $\mc{G}_0\times_{\mc{H}_0} \mc{H}_1$ along with maps $pr_1\colon \mc{G}_0\times_{\mc{H}_0} \mc{H}_1\ra \mc{G}_0, s\circ pr_2\colon \mc{G}_0\times_{\mc{H}_0} \mc{H}_1\ra \mc{H}_0$ produce a $\mc{G}-\mc{H}$ bibundle. 
This $\mc{G}-\mc{H}$ bibundle is described by the following diagram,
\begin{equation} \begin{tikzcd}
\mc{G}_1 \arrow[dd,xshift=0.75ex,"t"]\arrow[dd,xshift=-0.75ex,"s"'] & & \mc{H}_1 \arrow[dd,xshift=0.75ex,"t"]\arrow[dd,xshift=-0.75ex,"s"'] \\
& \mc{G}_0\times_{\mc{H}_0} \mc{H}_1 \arrow[rd, "s\circ pr_2"] \arrow[ld, "pr_1"'] & \\
\mc{G}_0 & & \mc{H}_0
\end{tikzcd}.\end{equation} 
\begin{construction}\label{Consturction:bibundlefromMorphismofLiegroupoids}
	Given a morphism of Lie groupoids $\phi\colon \mc{G}\rightarrow \mc{H}$, the manifold $\mc{G}_0\times_{\mc{H}_0} \mc{H}_1$ along with the maps $pr_1\colon \mc{G}_0\times_{\mc{H}_0} \mc{H}_1\ra \mc{G}_0, s\circ pr_2\colon \mc{G}_0\times_{\mc{H}_0} \mc{H}_1\ra \mc{H}_0$ is a $\mc{G}-\mc{H}$ bibundle.
	We denote the manifold $\mc{G}_0\times_{\mc{H}_0}\mc{H}_1$ by $\phi^*\mc{H}_1$ and the $\mc{G}-\mc{H}$ bibundle by $\left<\phi\right>\colon \mc{G}\ra\mc{H}$. 
	  
\end{construction}
\begin{remark}
	As a special case, when $\phi\colon \mc{G}\ra \mc{H}$ is a Lie groupoid extension, the $\mc{G}-\mc{H}$ bibundle associated to $\phi\colon \mc{G}\ra \mc{H}$ in Construction \ref{Consturction:bibundlefromMorphismofLiegroupoids} is  explained by the following diagram,
	\begin{equation}\label{Diagram:GHbibundleH1} \begin{tikzcd}
	\mc{G}_1 \arrow[dd,xshift=0.75ex,"t"]\arrow[dd,xshift=-0.75ex,"s"'] & & \mc{H}_1 \arrow[dd,xshift=0.75ex,"t"]\arrow[dd,xshift=-0.75ex,"s"'] \\
	& \mc{H}_1 \arrow[rd, "s"] \arrow[ld, "t"'] & \\
	\mc{G}_0 & & \mc{H}_0
	\end{tikzcd}.\end{equation} 
	  \end{remark}
\begin{lemma}\cite[Lemma $3.34$]{MR2778793}
	\label{Lemma:<f>isGprincipalbibundle}
	A morphism of Lie groupoids $f\colon \mc{G}\rightarrow \mc{H}$ is a Morita morphism of Lie groupoids if and only if the corresponding $\mc{G}-\mc{H}$ bibundle $\left<f\right>\colon \mc{G}\ra \mc{H}$ (mentioned in Construction \ref{Consturction:bibundlefromMorphismofLiegroupoids}) is a $\mc{G}$-principal bibundle (Remark \ref{Remark:GPrincipalbibunde}).
	  \end{lemma} 
\subsubsection{A bibundle gives a morphism of stacks}\label{SubSubsection:MorphismofStacksassociatedtobibundle} In this subsection, given a $\mc{G}-\mc{H}$ bibundle $P\colon \mc{G}\rightarrow \mc{H}$, we associate a morphism of stacks $BP\colon B\mc{G}\rightarrow B\mc{H}$. 

Before we describe the general construction, let us consider a special situation. Suppose that the Lie groupoids $\mc{G}$ and $\mc{H}$ are of the form $\mc{G}=(G\rightrightarrows *)$ and $\mc{H}=(H\rightrightarrows *)$ for Lie groups $G$ and $H$.
In this case, $B\mc{G}$ is the collection of principal $G$-bundles, and $B\mc{H}$ is likewise.

In this set up, a $\mc{G}-\mc{H}$ bibundle is given by a smooth manifold $P$ with an action of $G$ from left side and an action of $H$ from right side as in the following diagram,
\begin{equation} \begin{tikzcd}
G \arrow[dd,xshift=0.75ex,"t"]\arrow[dd,xshift=-0.75ex,"s"'] & & H \arrow[dd,xshift=0.75ex,"t"]\arrow[dd,xshift=-0.75ex,"s"'] \\
& P \arrow[ld] \arrow[rd] & \\
* & & *
\end{tikzcd}.\end{equation} 
The condition that $P\rightarrow *$ is a principal $H$-bundle implies that $P=H$ (up to an isomorphism). So, in this setup, a $\mc{G}-\mc{H}$ bibundle is nothing but an action of $G$ on $H$ from left. 
Given a left action of $G$ on $H$, our task is to associate a morphism of stacks $B\mc{G}\rightarrow B\mc{H}$; that is, a morphism of stacks $BG\ra BH$.

There is a classical construction of a principal $H$-bundle for a given principal $G$-bundle and an action of $G$ on $H$ (Chapter $1$ in \cite{KobNom1}). Here, we briefly recall the construction given in \cite{KobNom1}. 
Given a principal $G$-bundle $\pi\colon Q\rightarrow M$ and a left action of $G$ on $H$, we have an action of $G$ on $Q\times H$, given by $ g\cdot (q,h)=(qg,g^{-1}h)$.
The projection map $pr_1\colon Q\times H\ra Q$ induces the map $\widetilde{pr_1}\colon (Q\times H)/G\ra Q/G\cong M$. This produces a principal $H$-bundle $(Q\times H)/G\ra M$.
The following diagram illustrates the construction,
\begin{equation}\label{Diagram:BGBHLiegroupsmap}
\begin{tikzcd}
(Q\times H)/G \arrow[dd, "\widetilde{pr_1}"'] & & Q\times H \arrow[ldd, "pr_1"'] \arrow[rdd, "pr_2"] \arrow[ll, "\text{quotient}"'] & & \\
& & G \arrow[dd,xshift=0.75ex,"t"]\arrow[dd,xshift=-0.75ex,"s"'] & & H \arrow[dd,xshift=0.75ex,"t"]\arrow[dd,xshift=-0.75ex,"s"'] \\
Q/G \arrow[d, "\cong","\widetilde{\pi}"'] & Q \arrow[ld] \arrow[rd] \arrow[ld, "\pi"'] \arrow[l, "\text{quotient}"'] & & H \arrow[ld] \arrow[rd] & \\
M & & * & & * 
\end{tikzcd}
\end{equation}
See the above principal $H$-bundle as,
\begin{equation} \begin{tikzcd}
& & H \arrow[dd,xshift=0.75ex,"t"]
\arrow[dd,xshift=-0.75ex,"s"'] \\
& (Q\times H)/G \arrow[rd] \arrow[ld,"\widetilde{\pi}\circ \widetilde{pr_1}"'] & \\
M & & *
\end{tikzcd}.\end{equation}
The functor $BG\ra BH$ at the level of morphisms is obvious.
This defines a morphism of stacks $BP\colon BG\rightarrow BH$. 

Now, we consider the general construction of a morphism of stacks $BP\colon B\mc{G}\ra B\mc{H}$ from a $\mc{G}-\mc{H}$ bibundle $P\colon \mc{G}\ra \mc{H}$. Let $\pi\colon Q\rightarrow M$ be a principal $\mc{G}$-bundle.
The following diagram gives a principal $\mc{H}$-bundle,
\[\begin{tikzcd}
(Q\times_{\mc{G}_0}P)/\mc{G}_1 \arrow[dd, "\widetilde{pr_1}"'] & & Q\times_{\mc{G}_0}P \arrow[ldd, "pr_1"'] \arrow[rdd, "pr_2"] \arrow[ll, "\text{quotient}"'] & & \\
& & \mc{G}_1 \arrow[dd,xshift=0.75ex,"t"]\arrow[dd,xshift=-0.75ex,"s"'] & & \mc{H}_1 \arrow[dd,xshift=0.75ex,"t"]\arrow[dd,xshift=-0.75ex,"s"'] \\
Q/\mc{G}_1 \arrow[d, "\cong","\widetilde{\pi}"'] & Q \arrow[ld] \arrow[rd] \arrow[ld, "\pi"'] \arrow[l, "\text{quotient}"'] & & P \arrow[ld,"a_{\mc{G}}"'] \arrow[rd] & \\
M & & \mc{G}_0 & & \mc{H}_0 
\end{tikzcd}.\]
For our convenience, we interpret the above diagram as
\begin{equation}\label{Diagram:BGBHLiegroupoidsmap}
\begin{tikzcd}
& & (Q\times_{\mc{G}_0} P)/\mc{G}_1 \arrow[ldd] \arrow[rdd] & & \\
& & \mc{G}_1 \arrow[dd,xshift=0.75ex,"t"]\arrow[dd,xshift=-0.75ex,"s"'] & & \mc{H}_1 \arrow[dd,xshift=0.75ex,"t"]\arrow[dd,xshift=-0.75ex,"s"'] \\
& Q \arrow[ld,"\pi"'] \arrow[rd] & & P \arrow[ld,"a_{\mc{G}}"'] \arrow[rd,"a_{\mc{H}}"] & \\
M & & \mc{G}_0 & & \mc{H}_0
\end{tikzcd}.\end{equation} 
At the level of objects, the morphism of stacks $BP\colon B\mc{G}\rightarrow B\mc{H}$ defined as 
\begin{equation}\label{Equation:DefinitionofBP}
BP(\pi\colon Q\rightarrow M)=(\widetilde{\pi}\circ \widetilde{pr_1}\colon (Q\times_{\mc{G}_0}P)/\mc{G}_1\rightarrow M).
\end{equation}
At the level of morphisms, it is defined similarly as in the case of $\mc{G}=(G\rightrightarrows *)$ and $\mc{H}=(H\rightrightarrows *)$.
Thus, given a $\mc{G}-\mc{H}$ bibundle $P\colon \mc{G}\rightarrow \mc{H}$ we have associated a morphism of stacks $BP\colon B\mc{G}\rightarrow B\mc{H}$. 
\begin{construction}\label{Construction:BPassociatedtobibundle}
	Given a $\mc{G}-\mc{H}$ bibundle $P\colon \mc{G}\rightarrow \mc{H}$, we have a morphism of stacks $BP\colon B\mc{G}\rightarrow B\mc{H}$ defined as in Equation \ref{Equation:DefinitionofBP}.
\end{construction}
Combining Constructions \ref{Consturction:bibundlefromMorphismofLiegroupoids} and \ref{Construction:BPassociatedtobibundle}, we have the following result.
\begin{construction}\label{Construction:morphismofLiegroupoidstoMorphismofstacks}
	Given a morphism of Lie groupoids $f\colon \mc{G}\rightarrow \mc{H}$, we have a morphism of stacks $BP\colon B\mc{G}\rightarrow B\mc{H}$.  
\end{construction}
Next, we want to construct a weak $2$-category whose objects are Lie groupoids, and morphisms are bibundles. We need the notion of composition of bibundles.
The idea of composition of bibundles is same as that of constructing $BP\colon B\mc{G}\ra B\mc{H}$ from a  given a $\mc{G}-\mc{H}$ bibundle $P\colon \mc{G}\ra \mc{H}$. 

Let $P\colon \mc{G}\ra \mc{H}$ be a $\mc{G}-\mc{H}$ bibundle and 
$Q\colon \mc{H}\ra\mc{H}'$ be a $\mc{H}-\mc{H}'$ bibundle. We have the
following diagrams for bibundles,
\begin{equation} 
\begin{tikzcd}
\mc{G}_1 \arrow[dd,xshift=0.75ex,"t"]\arrow[dd,xshift=-0.75ex,"s"'] & & \mc{H}_1 \arrow[dd,xshift=0.75ex,"t"]\arrow[dd,xshift=-0.75ex,"s"'] & & \mc{H}_1'\arrow[dd,xshift=0.75ex,"t"]\arrow[dd,xshift=-0.75ex,"s"'] \\
& P \arrow[ld] \arrow[rd] & & Q \arrow[ld] \arrow[rd ] & \\
\mc{G}_0 & & \mc{H}_0 & & \mc{H}'_0 
\end{tikzcd}.
\end{equation}
Ignoring the  action of $\mc{G}$ on $P$, we can consider $a_{\mc{G}}\colon P\ra \mc{G}_0$ as a principal $\mc{H}$-bundle,
\begin{equation} 
\begin{tikzcd}
& & \mc{H}_1 \arrow[dd,xshift=0.75ex,"t"]\arrow[dd,xshift=-0.75ex,"s"'] & & \mc{H}_1'\arrow[dd,xshift=0.75ex,"t"]\arrow[dd,xshift=-0.75ex,"s"'] \\
& P \arrow[ld,"a_{\mc{G}}"'] \arrow[rd] & & Q \arrow[ld] \arrow[rd] & \\
\mc{G}_0 & & \mc{H}_0 & & \mc{H}'_0 
\end{tikzcd}.
\end{equation}
Given 
a principal $\mc{H}$-bundle and a $\mc{H}-\mc{H}'$ bibundle, we know (equation \ref{Equation:DefinitionofBP}) how to associate a principal $\mc{H}'$-bundle. For the principal $\mc{H}$-bundle $a_{\mc{G}}\colon P\ra \mc{G}_0$, we associate the principal $\mc{H}'$-bundle $BQ(a_{\mc{G}})\colon (P\times_{\mc{H}_0}Q)/\mc{H}_1\ra \mc{G}_0$. The following diagram illustrates the construction,
\begin{equation}
\begin{tikzcd}
& & (P\times_{\mc{H}_0} Q)/\mc{H}_1 \arrow[ldd] \arrow[rdd] & & \\
& & \mc{H}_1 \arrow[dd,xshift=0.75ex,"t"]\arrow[dd,xshift=-0.75ex,"s"'] & & \mc{H}'_1 \arrow[dd,xshift=0.75ex,"t"]\arrow[dd,xshift=-0.75ex,"s"'] \\
& P \arrow[ld] \arrow[rd] & & Q \arrow[ld] \arrow[rd] & \\
\mc{G}_0 & & \mc{H}_0 & & \mc{H}'_0
\end{tikzcd}.\end{equation} 

Action of $\mc{G}$ on $P$ induces an action of $\mc{G}$ on $(P\times_{\mc{H}_0}Q)/\mc{H}_1$, producing the following $\mc{G}-\mc{H}'$ bibundle,
\begin{equation}
\label{Diagram:compositionOfbibundles}\begin{tikzcd}
\mc{G}_1 \arrow[dd,xshift=0.75ex,"t"]\arrow[dd,xshift=-0.75ex,"s"'] & & \mc{H}'_1 \arrow[dd,xshift=0.75ex,"t"]\arrow[dd,xshift=-0.75ex,"s"'] \\
& (P\times_{\mc{H}_0} Q)/\mc{H}_1 \arrow[ld] \arrow[rd] & \\
\mc{G}_0 & & \mc{H}'_0 
\end{tikzcd}.\end{equation}
\begin{definition}\label{Definition:compositionOfbibundles}
	Let $P\colon \mc{G}\ra \mc{H}$ be a $\mc{G}-\mc{H}$ bibundle and $Q\colon \mc{H}\ra \mc{H}'$ be a 
	$\mc{H}-\mc{H}'$ bibundle. We define \textit{the composition of $Q$ with $P$} to be the $\mc{G}-\mc{H}'$ bibundle
	\begin{equation}\label{Equation:compositionOfbibundles}
	Q\circ P=(P\times_{\mc{H}_0} Q)/\mc{H}_1
	\end{equation}
	as in the Diagram \ref{Diagram:compositionOfbibundles}.
	  \end{definition}
Recall (Theorem \ref{Theorem:BGBHisomorphicimpliesGHareME}) that, for Lie groupoids $\mc{G}$ and $\mc{H}$, if the stacks $B\mc{G}$ and $B\mc{H}$ are isomorphic, then $\mc{G}$ and $\mc{H}$ are Morita equivalent Lie groupoids. Now we prove that, if $\mc{G}$ and $\mc{H}$ are Morita equivalent Lie groupoids, then the stacks $B\mc{G}$ and $B\mc{H}$ are isomorphic. 

\begin{proposition}
	\label{Proposition:Moritaequivalentimpliesisomorphicstacks}
	Let $\mc{H},\mc{H}'$ be Morita equivalent Lie groupoids (Definition \ref{Definition:MoritaEquivalentLiegroupoids}), then, the stacks $B\mc{H}$ and $B\mc{H}'$ are isomorphic.
	\begin{proof} Let $\mc{H}$ and $\mc{H}'$ be Morita equivalent Lie groupoids; that is, there exists a Lie groupoid $\mc{G}$ and a pair of Morita morphisms of Lie groupoids $f\colon \mc{G}\rightarrow \mc{H}$ and $ g\colon \mc{G}\rightarrow \mc{H}'$. With this data, we produce an isomorphism of stacks $B\mc{H}\rightarrow B\mc{H}'$.
		
		Recall that (Lemma \ref{Lemma:mapBGtoBHdeterminesGHbibundle}), giving a morphism of stacks $B\mc{H}\rightarrow B\mc{H}'$ is same as giving a $\mc{H}-\mc{H}'$ bibundle. Here, we take a $\mc{H}-\mc{H}'$ bibundle to represent a morphism of stacks $B\mc{H}\rightarrow B\mc{H}'$.
		The morphism of Lie groupoids $f\colon \mc{G}\rightarrow \mc{H}$ gives the following $\mc{G}-\mc{H}$ bibundle (Construction \ref{Consturction:bibundlefromMorphismofLiegroupoids}),
		\begin{equation} \begin{tikzcd}
		\mc{G}_1 \arrow[dd,xshift=0.75ex,"t"]\arrow[dd,xshift=-0.75ex,"s"'] & & \mc{H}_1 \arrow[dd,xshift=0.75ex,"t"]\arrow[dd,xshift=-0.75ex,"s"'] \\
		& \left<f\right> \arrow[rd, "s\circ pr_2"] \arrow[ld, "pr_1"'] & \\
		\mc{G}_0 & & \mc{H}_0
		\end{tikzcd}.\end{equation} 
		As $f\colon \mc{G}\rightarrow \mc{H}$ is a Morita morphism of Lie groupoids, Lemma \ref{Lemma:<f>isGprincipalbibundle} says that $\left<f\right>\colon \mc{G}\ra \mc{H}$ is a $\mc{G}$-principal bibundle. Thus, $\left<f\right>\colon \mc{G}\ra \mc{H}$ can be considered as a $\mc{H}-\mc{G}$ bibundle,
		\begin{equation}\label{Diagram:<f>asHGbibundle}\begin{tikzcd}
		\mc{H}_1 \arrow[dd,xshift=0.75ex,"t"]\arrow[dd,xshift=-0.75ex,"s"'] & & \mc{G}_1 \arrow[dd,xshift=0.75ex,"t"]\arrow[dd,xshift=-0.75ex,"s"'] \\
		& \left<f\right> \arrow[rd, "pr_1"] \arrow[ld, "s\circ pr_2"'] & \\
		\mc{H}_0 & & \mc{G}_0
		\end{tikzcd}.\end{equation} 
		The morphism of Lie groupoids $g\colon \mc{G}\rightarrow \mc{H}'$ gives the following $\mc{G}-\mc{H}'$ bibundle,
		\begin{equation}\label{Diagram:<g>asGH'bibundle} \begin{tikzcd}
		\mc{G}_1 \arrow[dd,xshift=0.75ex,"t"]\arrow[dd,xshift=-0.75ex,"s"'] & & \mc{H}'_1 \arrow[dd,xshift=0.75ex,"t"]\arrow[dd,xshift=-0.75ex,"s"'] \\
		& \left<g\right> \arrow[ld, "pr_1"] \arrow[rd, "s\circ pr_2"'] & \\
		\mc{G}_0 & & \mc{H}'_0
		\end{tikzcd}.\end{equation} 
		Composing the $\mc{H}-\mc{G}$ bibundle (Diagram \ref{Diagram:<f>asHGbibundle}) $\left<f\right>\colon \mc{H}\ra \mc{G}$ with the $\mc{G}-\mc{H}'$ bibundle (Diagram \ref{Diagram:<g>asGH'bibundle}) $\left<g\right>\colon \mc{G}\ra \mc{H}'$, 
		we get the $\mc{H}-\mc{H}'$ bibundle $\left<g\right>\circ \left<f\right>\colon \mc{H}\rightarrow \mc{H}'$ (equation \ref{Equation:compositionOfbibundles}), as explained  in the following diagram,
		\begin{equation} \begin{tikzcd}
		& & \left<g\right>\circ \left<f\right> \arrow[ldd] \arrow[rdd] & & \\
		\mc{H}_1 \arrow[dd,xshift=0.75ex,"t"]\arrow[dd,xshift=-0.75ex,"s"'] & & \mc{G}_1 \arrow[dd,xshift=0.75ex,"t"]\arrow[dd,xshift=-0.75ex,"s"'] & & \mc{H}'_1 \arrow[dd,xshift=0.75ex,"t"]\arrow[dd,xshift=-0.75ex,"s"'] \\
		& \left<f\right> \arrow[ld] \arrow[rd] & & \left<g\right> \arrow[ld] \arrow[rd] & \\
		\mc{H}_0 & & \mc{G}_0 & & \mc{H}'_0
		\end{tikzcd}.\end{equation} 
		As $g\colon \mc{G}\rightarrow \mc{H}'$ is also a Morita morphism of Lie groupoids, interchanging $f$ and $g$ we obtain a $\mc{H}'-\mc{H}$ bibundle as follows,
		\begin{equation} \begin{tikzcd}
		& & \left<f\right>\circ \left<g\right> \arrow[ldd] \arrow[rdd] & & \\
		\mc{H}'_1 \arrow[dd,xshift=0.75ex,"t"]\arrow[dd,xshift=-0.75ex,"s"'] & & \mc{G}_1 \arrow[dd,xshift=0.75ex,"t"]\arrow[dd,xshift=-0.75ex,"s"'] & & \mc{H}_1 \arrow[dd,xshift=0.75ex,"t"]\arrow[dd,xshift=-0.75ex,"s"'] \\
		& \left<g\right> \arrow[ld] \arrow[rd] & & \left<f\right> \arrow[ld] \arrow[rd] & \\
		\mc{H}'_0 & & \mc{G}_0 & & \mc{H}_0
		\end{tikzcd}.\end{equation} 
		The $\mc{H}-\mc{H}'$ bibundle $\left<g\right>\circ \left<f\right>\colon \mc{H}\rightarrow \mc{H}'$ gives a morphism of stacks $B\mc{H}\rightarrow B\mc{H}'$ and the $\mc{H}'-\mc{H}$ bibundle $\left<f\right>\circ \left<g\right>\colon \mc{H}'\rightarrow \mc{H}$ gives a morphism of stacks $B\mc{H}'\rightarrow B\mc{H}$. It is easy to see that the maps $B\mc{H}'\rightarrow B\mc{H}$ and $B\mc{H}\rightarrow B\mc{H}'$ are inverses to each other, giving an isomorphism of stacks $B\mc{H}\rightarrow B\mc{H}'$. Thus, the stacks $B\mc{H}$ and $B\mc{H}'$ are isomorphic.
	\end{proof}
\end{proposition} 

\subsection{A Lie groupoid extension gives a gerbe over a stack}
\label{Subsection:LieGpdExtngivesGoS}
Let $\phi\colon (\mc{G}_1\rightrightarrows \mc{G}_0)\rightarrow (\mc{H}_1\rightrightarrows \mc{H}_0)$ be a Lie groupoid extension. 
We have described a construction of a morphism of stacks $F\colon B\mc{G}\ra B\mc{H}$ from a morphism of Lie groupoids $f\colon \mc{G}\ra \mc{H}$ (Construction \ref{Construction:morphismofLiegroupoidstoMorphismofstacks}). In particular, given a Lie groupoid extension $\phi\colon (\mc{G}_1\rightrightarrows \mc{G}_0)\ra (\mc{H}_1\rightrightarrows \mc{H}_0)$, we have a morphism of stacks $F\colon B\mc{G}\ra B\mc{H}$, which at the level of objects have the following description,
\begin{equation}\label{Equation:definitionofFBGBH}
F(\pi\colon Q\ra M)=(\widetilde{\pi}\circ \widetilde{pr_1}\colon (Q\times_{\mc{G}_0}\mc{H}_1)/\mc{G}_1\rightarrow M).
\end{equation} 
The following diagram (using Diagram \ref{Diagram:GHbibundleH1}) explains the same,
\begin{equation} \begin{tikzcd}
& & (Q\times_{\mc{G}_0} \mc{H}_1)/\mc{G}_1 \arrow[ldd] \arrow[rdd] & & \\
& & \mc{G}_1 \arrow[dd,xshift=0.75ex,"t"]\arrow[dd,xshift=-0.75ex,"s"'] & & \mc{H}_1 \arrow[dd,xshift=0.75ex,"t"]\arrow[dd,xshift=-0.75ex,"s"'] \\
& Q \arrow[ld] \arrow[rd] & & \mc{H}_1 \arrow[ld,"t"'] \arrow[rd,"s"] & \\
M & & \mc{G}_0 & & \mc{G}_0
\end{tikzcd}.\end{equation} 
Here we prove that the morphism of stacks $F\colon B\mc{G}\ra B\mc{H}$ is a gerbe over a stack (Definition \ref{Definition:GerbeoverstackasinMR2817778}). That is, the morphism of stacks $F\colon B\mc{G}\ra B\mc{H}$ and the diagonal morphism $\Delta_F\colon B\mc{G}\ra B\mc{G}\times_{B\mc{H}}B\mc{G}$ associated to $F$ are epimorphisms of stacks.

\subsubsection{Proof that $F\colon B\mc{G}\ra B\mc{H}$ is an epimorphism}
\label{SubSubsection:FisanEpimorphism} Given a manifold $U$ and a morphism of stacks $q\colon \underline{U}\ra B\mc{H}$, we prove that, there exists an open cover $\{U_i\}$ of $U$ and a morphism of stacks $l_i\colon \underline{U_i}\ra B\mc{G}$, for each $i$, giving the following $2$-commutative diagram,
\begin{equation}\label{Diagram:UiUBGBH}
\begin{tikzcd}
\underline{U_i} \arrow[dd, "l_i"'] \arrow[rr,"\Phi=\text{inclusion}"] & & \underline{U} \arrow[dd, "q"] \\
& & \\
B\mc{G} \arrow[Rightarrow, shorten >=20pt, shorten <=20pt, uurr] \arrow[rr, "F"] & & B\mc{H}
\end{tikzcd}.
\end{equation}
This will prove that $F\colon B\mc{G}\ra B\mc{H}$ is an epimorphism of stacks (Definition \ref{Definition:epimorphismofStacks}).

Let $\pi\colon P\ra U$ be the principal $\mc{H}$-bundle associated to the morphism of stacks $q\colon \underline{U}\ra B\mc{H}$ (Lemma \ref{Lemma:MorphismfromMtoBG}). 
For  $\pi\colon P\ra U$, there exists  an open cover $\{U_i\}$ of $U$ and a map $r_i\colon U_i\ra \mc{H}_0=M$, such that $\pi|_{\pi^{-1}(U_i)}\colon \pi^{-1}(U_i)\ra U_i$ is the pullback of $t\colon \mc{H}_1\ra \mc{H}_0=M$ along $r_i\colon U_i\ra M$ as explained by the following diagram (see Corollary \ref{Corollary:LocalPropertyofprincipalLiegroupoidbundle}),
\begin{equation}\label{Diagram:pullbackoftalongri}
\begin{tikzcd}
\pi^{-1}(U_i) \arrow[dd,"\pi|_{\pi^{-1}(U_i)}"'] \arrow[rr] & & \mc{H}_1 \arrow[dd, "t"] \\
& & \\
U_i \arrow[rr, "r_i"] & & \mc{H}_0=M 
\end{tikzcd}. 
\end{equation}
Now, pullback the principal $\mc{G}$-bundle $t\colon \mc{G}_1\ra \mc{G}_0=M$ along $r_i\colon U_i\ra \mc{G}_0=M$ to get the principal $\mc{G}$-bundle $l_i\colon W_i\ra U_i$,
\begin{equation}\label{Diagram:WiG1UiM}
~~~~~~~\begin{tikzcd}
W_i \arrow[dd,"l_i"'] \arrow[rr] & & \mc{G}_1 \arrow[dd, "t"] \\
& & \\
U_i \arrow[rr, "r_i"] & & \mc{G}_0=M 
\end{tikzcd}. 
\end{equation}
This principal $\mc{G}$-bundle $l_i\colon W_i\ra U_i$ gives a morphism of stacks $\underline{U_i} \ra B\mc{G}$, which we denote by $l_i$ (by abuse of notation). So, we have the morphism of stacks $l_i\colon \underline{U_i}\ra B\mc{G}$ for each $i$. This gives a pair of  compositions of morphisms of stacks $q\circ \Phi\colon \underline{U}_i\ra\underline{U}\ra B\mc{H}$ and $F\circ l_i\colon \underline{U}_i\ra B\mc{G}\ra B\mc{H}$.
We prove that these two compositions give the $2$-commutative diagram \ref{Diagram:UiUBGBH}, which would then imply that $F\colon B\mc{G}\ra B\mc{H}$ is an epimorphism of stacks. For that, we need the following lemma.
\begin{lemma}\label{Lemma:compatabilityofFwithpullback}
	Let $\pi\colon Q\ra U$ be the pullback of the principal $\mc{G}$-bundle $t\colon \mc{G}_1\ra M$ along a smooth map $r\colon U\ra M$. Then $F(\pi\colon Q\ra U)$ is the pullback of the principal $\mc{H}$-bundle $F(t\colon \mc{G}_1\ra M)$ along the smooth map $r\colon U\ra M$.
	\begin{proof}
		Consider the following pullback diagram,
		\begin{equation}\label{Diagram:pullbackoftalongr}\begin{tikzcd}
		Q \arrow[dd, "\pi"'] \arrow[rr] & & \mc{G}_1 \arrow[dd, "t"] \\
		& & \\
		U \arrow[rr, "r"] & & M
		\end{tikzcd}.\end{equation}
		We have $F(t\colon \mc{G}_1\ra M)=(\widetilde{t_{\mc{G}}}\circ \widetilde{pr_1}\colon (\mc{G}_1\times_M \mc{H}_1)/\mc{G}_1\ra M)$ (Equation \ref{Equation:definitionofFBGBH}) which can be expressed by the following diagram,
		\begin{equation}\begin{tikzcd}
		& & (\mc{G}_1\times_M \mc{H}_1)/\mc{G}_1 \arrow[ldd] \arrow[rdd] & & \\
		& & \mc{G}_1 \arrow[dd,xshift=0.75ex,"t"]\arrow[dd,xshift=-0.75ex,"s"'] & & \mc{H}_1 \arrow[dd,xshift=0.75ex,"t"]\arrow[dd,xshift=-0.75ex,"s"'] \\
		& \mc{G}_1 \arrow[ld, "t_{\mc{G}}"'] \arrow[rd, "s_{\mc{G}}"] & & \mc{H}_1 \arrow[ld, "t_{\mc{H}}"'] \arrow[rd] & \\
		M & & M & & M 
		\end{tikzcd}.\end{equation}
		Adjoining the pullback diagram (Diagram \ref{Diagram:pullbackoftalongr}) with the above diagram, we have the following diagram,
		\begin{equation}\begin{tikzcd}
		& & & (\mc{G}_1\times_M \mc{H}_1)/\mc{G}_1 \arrow[ldd] \arrow[rdd] & & \\
		& Q \arrow[ld, "\pi"'] \arrow[rd] & & \mc{G}_1 \arrow[dd,xshift=0.75ex,"t"]\arrow[dd,xshift=-0.75ex,"s"'] & & \mc{H}_1 \arrow[dd,xshift=0.75ex,"t"]\arrow[dd,xshift=-0.75ex,"s"'] \\
		U \arrow[rd, "r"] & & \mc{G}_1 \arrow[ld, "t_{\mc{G}}"'] \arrow[rd, "s_{\mc{G}}"] & & \mc{H}_1 \arrow[ld, "t_{\mc{H}}"'] \arrow[rd] & \\
		& M & & M & & M 
		\end{tikzcd}.\end{equation}
		Thus, we have $F(\pi\colon Q\ra U)=((Q\times_M\mc{H}_1)/\mc{G}_1\ra U)$. 
		As $Q=U\times_M\mc{G}_1$, we have
		\begin{align*}
		(Q\times_M\mc{H}_1)/\mc{G}_1=
		(U\times_M\mc{G}_1\times_M\mc{H}_1)/\mc{G}_1=U\times_M (\mc{G}_1\times_M\mc{H}_1)/\mc{G}_1.
		\end{align*}
		Note that $U\times_M (\mc{G}_1\times_M\mc{H}_1)/\mc{G}_1$ is precisely the pullback of $(\mc{G}_1\times_M \mc{H}_1)/\mc{G}_1$ along $r\colon U\ra M$. Thus, $F(\pi\colon Q\ra U)$ is the pullback of $F(t\colon \mc{G}_1\ra M)$ along $r\colon U\ra M$.
	\end{proof}
\end{lemma}
As the Diagrams \ref{Diagram:pullbackoftalongri} and \ref{Diagram:WiG1UiM} are pullback diagrams, observe that 
\[F(l_i:W_i\ra U_i)=(\pi|_{\pi^{-1}(U_i)}\colon \pi^{-1}(U_i)\ra U_i),\] and 
\[(\pi|_{\pi^{-1}(U_i)}\colon \pi^{-1}(U_i)\ra U_i)=
q(U_i\ra U)=(q\circ \Phi)(\text{Id}\colon U_i\ra U_i).\] Here, $(W_i\ra U_i)=l_i(\text{Id}\colon U_i\ra U_i)$.
So, $(F\circ l_i) (\text{Id}\colon U_i\ra U_i)$ is equal to $(q\circ \Phi)(\text{Id}\colon U_i\ra U_i)$. So, there is an isomorphism $(F\circ l_i)(\text{Id}\colon U_i\ra U_i)\rightarrow (q\circ \Phi)(\text{Id}\colon U_i\ra U_i)$. For the same reason, it turns out that there is an isomorphism 
$(F\circ l_i)(f\colon N\ra U_i)\rightarrow (q\circ \Phi)(f\colon N\ra U_i)$ for each $f\colon N\ra U$ in $\underline{U_i}$. Thus, we have the following $2$-commutative diagram, 
\begin{equation}
\begin{tikzcd}
\underline{U_i} \arrow[dd, "l_i"'] \arrow[rr,"\Phi=\text{inclusion}"] & & \underline{U} \arrow[dd, "q"] \\
& & \\
B\mc{G} \arrow[Rightarrow, shorten >=20pt, shorten <=20pt, uurr] \arrow[rr, "F"] & & B\mc{H}
\end{tikzcd}.
\end{equation}
Thus, $F\colon B\mc{G}\ra B\mc{H}$ is an epimorphism of stacks. We summarize the discussion as follows.
\begin{proposition}\label{Proposition:BGBHisanepimorphism}
	Given a Lie groupoid extension $f\colon (\mc{G}_1\rightrightarrows M)\ra (\mc{H}_1\rightrightarrows M)$, the corresponding morphism of stacks $F\colon B\mc{G}\ra B\mc{H}$ is an epimorphism of stacks.
\end{proposition}
\subsubsection{Proof that the diagonal morphism $\Delta_F\colon B\mc{G}\ra B\mc{G}\times_{B\mc{H}}B\mc{G}$ is an epimorphism.}
\label{SubSubsection:diagonalisEpimorphism} 

As $\phi\colon \mc{G}\ra \mc{H}$ is a Lie groupoid extension (in particular, $\phi_1\colon \mc{G}_1\ra \mc{H}_1$ is a submersion), the $2$-fiber product $\mc{G}\times_{\mc{H}}\mc{G}$ is a Lie groupoid (Definition \ref{Definition:2-fiberproductinLiegroupoids}). As stackification and Yoneda embedding preserves the $2$-fiber product (\cite[$\text{I}.2.4$]{Carchedi}, \cite[\href{https://stacks.math.columbia.edu/tag/04Y1}{Tag 04Y1}]{Johan}), we see that \[B\mc{G}\times_{B\mc{H}}B\mc{G}\cong B(\mc{G}\times_{\mc{H}}\mc{G}).\]
Further, the diagonal morphism of stacks $\Delta_F\colon B\mc{G}\ra B\mc{G}\times_{B\mc{H}}B\mc{G}$ is the morphism of stacks associated to the diagonal morphism of Lie groupoids $\Delta_{\phi}\colon \mc{G}\ra \mc{G}\times_{\mc{H}}\mc{G}$, given by $\Delta_{\phi}(a)=(a,\text{Id}\colon a\ra a,a)$ and $\Delta_{\phi}(g)=(g,g)$ for $a\in \mc{G}_0$ and $g\in \mc{G}_1$.
We have the following morphism of Lie groupoids, 
\[\begin{tikzcd}
\mc{G}_1 \arrow[dd,xshift=0.75ex,"t"]\arrow[dd,xshift=-0.75ex,"s"'] \arrow[rr] & & (\mc{G}\times_{\mc{H}}\mc{G})_1 \arrow[dd,xshift=0.75ex,"t"]\arrow[dd,xshift=-0.75ex,"s"'] \\
& & \\
M \arrow[rr] & & (\mc{G}\times_{\mc{H}}\mc{G})_0 
\end{tikzcd}.\]
As the above morphism of Lie groupoids is not identity on base space, one can not immediately use Proposition \ref{Proposition:BGBHisanepimorphism} to conclude that $\Delta_F\colon B\mc{G}\ra B\mc{G}\times_{B\mc{H}}B\mc{G}$ is an epimorphism of stacks. We need a few more results to conclude the same.
\begin{remark}\label{Remark:descriptionof(GHG)0}
	As $\phi\colon (\mc{G}_1\rightrightarrows \mc{G}_0)\ra (\mc{H}_1\rightrightarrows \mc{H}_0)$ is a Lie groupoid extension, the map $\mc{G}_0\ra \mc{H}_0$ is an identity map. Thus, from equation \ref{Equation:GKH0}, we have \[(\mc{G}\times_{\mc{H}}\mc{G})_0= \mc{G}_0\times_{Id,\mc{H}_0,s}\mc{H}_1\times_{t\circ pr_2,\mc{K}_0,Id}\mc{H}_0=\mc{H}_1.\]
	  \end{remark}
\begin{lemma}\label{Lemma:GHGisTransitiveLieGroupoid}
	The Lie groupoid $\mc{G}\times_{\mc{H}}\mc{G}$ is a transitive Lie groupoid (Definition \ref{Definition:TransitiveLiegroupoid}). 
	\begin{proof}
		To prove $\mc{G}\times_{\mc{H}}\mc{G}$ is a transitive Lie groupoid, we prove that, for $h_1,h_2\in \mc{H}_1=(\mc{G}\times_{\mc{H}}\mc{G})_0$ (Remark \ref{Remark:descriptionof(GHG)0}) there exists 
		$(g_1,h,g_2)\in (\mc{G}\times_{\mc{H}}\mc{G})_1$ such that 
		$s(g_1,h,g_2)=h\circ \phi(g_1)=h_1$ and $t(g_1,h,g_2)=\phi(g_2)\circ h=h_2$ (equation \ref{Equation:DefinitionOfsourcemapforGHK} and equation \ref{Equation:DefinitionOftargetmapforGHK}).
		
		As $\phi\colon \mc{G}_1\ra \mc{H}_1$ is surjective, we can choose $g_1\in \mc{G}_1$ to be such that $\phi(g_1)=1_{s(h_1)}$. Choose such a $g_1\in \mc{G}_1$.
		Choose $h=h_1$ and $g_2\in \mc{G}_2$ such that 
		$\phi(g_2)=h_2\circ h_1^{-1}$. 
		So, given $h_1,h_2\in \mc{H}_1=(\mc{G}\times_{\mc{H}}\mc{G})_0$, there exists $(g_1,h,g_2)\in (\mc{G}\times_{\mc{H}}\mc{G})_1$ such that $s(g_1,h,g_2)=h\circ \phi(g_1)=h_1$ and $t(g_1,h,g_2)=\phi(g_2)\circ h=h_2$. Thus, $\mc{G}\times_{\mc{H}}\mc{G}$ is a transitive Lie groupoid.
	\end{proof}
\end{lemma}
\begin{lemma}\label{Lemma:TransitiveLiegroupoidisMEtoLiegroup}
	Any transitive Lie groupoid $\mc{G}$ (Definition \ref{Definition:isotropygroup}) is Morita equivalent to the Lie group $\mc{G}_x$ for any $x\in \mc{G}_0$, that is, 
	the Lie groupoid $(\mc{G}_1\rightrightarrows \mc{G}_0)$ is Morita equivalent to the Lie groupoid $(\mc{G}_x\rightrightarrows *)$.
	\begin{proof}
		Given a Lie groupoid $\mc{G}$ and an object $x$ in $\mc{G}_0$, we have a morphism of Lie groupoids $\psi\colon (\mc{G}_x\rightrightarrows *)\rightarrow (\mc{G}_1\rightrightarrows \mc{G}_0)$ given by $*\mapsto x$ at the level of objects and $g\mapsto g$ at the level of morphisms. The following diagram expresses this morphism,
		\[\begin{tikzcd}
		\mc{G}_x \arrow[dd,xshift=0.75ex,"t"]\arrow[dd,xshift=-0.75ex,"s"'] \arrow[rr,"\psi_1"] & & \mc{G}_1 \arrow[dd,xshift=0.75ex,"t"]\arrow[dd,xshift=-0.75ex,"s"'] \\
		& & \\
		* \arrow[rr,"\psi_0"] & & \mc{G}_0 
		\end{tikzcd}.\]
		Now we prove that $\psi\colon (\mc{G}_x\rightrightarrows *)\rightarrow (\mc{G}_1\rightrightarrows \mc{G}_0)$ is a Morita morphism of Lie groupoids. This implies that $(\mc{G}_1\rightrightarrows \mc{G}_0)$ is Morita equivalent to the Lie groupoid $(\mc{G}_x\rightrightarrows *)$.
		
		Observe that the morphism set of the pullback groupoid (Section \ref{Subsection:PullbackLiegroupoid}) is
		\begin{align*}
		*\times_{\psi_0,\mc{G}_0,s}\mc{G}_1\times_{t\circ pr_2,\mc{G}_0,\psi_0}*&=\{(*,g,*)|\psi_0(*)=s(g), (t\circ pr_2)(*,g)=\psi_0(*)\}\\
		&=\{ (*,g,*)|s(g)=x, t(g)=x \}\\
		&=\{(*,g,x)|g\in s^{-1}(x)\bigcap t^{-1}(x)\}\\
		&=\{(*,g,*)|g\in \mc{G}_x\}=\mc{G}_x.
		\end{align*}
		So, the pullback groupoid of the Lie groupoid $(\mc{G}_1\rightrightarrows \mc{G}_0)$ along the map $*\rightarrow \mc{G}_0$ is the Lie groupoid $ (\mc{G}_x\rightrightarrows *)$. Thus, we have a Morita equivalence of Lie groupoids $\psi\colon (\mc{G}_x\rightrightarrows *)\rightarrow (\mc{G}_1\rightrightarrows \mc{G}_0)$. 
	\end{proof}
\end{lemma}
Combining Lemma \ref{Lemma:GHGisTransitiveLieGroupoid} and Lemma \ref{Lemma:TransitiveLiegroupoidisMEtoLiegroup} we see that $\mc{G}\times_{\mc{H}}\mc{G}$ is Morita equivalent to a Lie groupoid of the form $(K\rightrightarrows *)$. Thus, by Proposition \ref{Proposition:Moritaequivalentimpliesisomorphicstacks}, 
the stacks $B(\mc{G}\times_{\mc{H}}\mc{G})$ and $B(K\rightrightarrows *)$ are isomorphic. So, the morphism of stacks $\Delta_F\colon B\mc{G}\ra B\mc{G}\times_{B\mc{H}}B\mc{G}$ is isomorphic to the map $B\mc{G}\ra B(K\rightrightarrows *)$.

Using an argument similar to the proof of  Proposition \ref{Proposition:BGBHisanepimorphism}, we conclude that for any morphism of Lie groupoids $(\mc{G}_1\rightrightarrows M)\ra (K\rightrightarrows *)$, the corresponding morphism of stacks $B\mc{G}\ra B\mc{K}$ is an epimorphism of stacks. Thus, $\Delta_F\colon B\mc{G}\ra B\mc{G}\times_{B\mc{H}}B\mc{G}$ is an epimorphism of stacks. Therefore we obtain the following:
\begin{proposition}\label{Proposition:DiagonalmapisanEpimorphism}
	Given a Lie groupoid extension $\phi\colon (\mc{G}_1\rightrightarrows M)\ra (\mc{H}_1\rightrightarrows M)$ the diagonal morphism of stacks $\Delta_F\colon B\mc{G}\ra B\mc{G}\times_{B\mc{H}}B\mc{G}$ is an epimorphism of stacks.
\end{proposition}
Finally we conclude the following.
\begin{theorem}\phantomsection
	\label{Theorem:BGBHisaGoS} 
	\begin{enumerate}
		\item{Given a Lie groupoid extension $\phi\colon (\mc{G}_1\rightrightarrows M)\ra (\mc{H}_1\rightrightarrows M)$, the corresponding morphism of stacks $F\colon B\mc{G}\ra B\mc{H}$ is a gerbe over the stack $B\mc{H}$.}
		\item{ Let $\phi\colon \mc{G}_1\ra \mc{H}_1\rightrightarrows M$ and $\phi''\colon \mc{G}''_1\ra \mc{H}''_1\rightrightarrows M''$ be Morita equivalent Lie groupoid extensions. Let $\Phi\colon B\mc{G}\ra B\mc{H}$ and $\Phi''\colon B\mc{G}''\ra B\mc{H}''$ be the respective morphism of stacks corresponding to $\phi$ and $\phi''$. Then, $\Phi$ and $\Phi''$ are isomorphic in the sense that, the following diagram is $2$-commutative,
			\[\begin{tikzcd}
			B\mc{G}'' \arrow[rr,"\Phi''"] \arrow[dd,"\cong"'] & & B\mc{H}'' \arrow[dd,"\cong"] \\
			& & \\
			B\mc{G} \arrow[rr,"\Phi"] \arrow[Rightarrow, shorten >=20pt, shorten <=20pt, uurr]& & B\mc{H} . 
			\end{tikzcd}\]}
	\end{enumerate}
	\begin{proof}
		\begin{enumerate}
			\item Immediate from Propositions \ref{Proposition:BGBHisanepimorphism} and \ref{Proposition:DiagonalmapisanEpimorphism}.
			\item Let the Morita equivalence be given by the Lie groupoid extension $\phi' \colon\mc{G}'_1\ra \mc{H}'_1\rightrightarrows M'$. In particular, that means we have a Morita morphism (Definition \ref{Definition:MoritamorphismofLiegroupoidExtension}) from the Lie groupoid extension $\phi' \colon\mc{G}'_1\ra \mc{H}'_1\rightrightarrows M'$ to $\phi \colon\mc{G}_1\ra \mc{H}_1\rightrightarrows M$, expressed by the following diagram, 
			\begin{equation}
			\label{Diagram:MELiegroupoidextensions}
			\begin{tikzcd}
			\mc{G}_1 \arrow[dd,xshift=0.75ex,"t"]\arrow[dd,xshift=-0.75ex,"s"'] \arrow[rrrrrrr, "\phi_1", bend left] & & \mc{G}_1' \arrow[ll, "\psi_{\mc{G}}"'] \arrow[rrr, "\phi_1'"] \arrow[dd,xshift=0.75ex,"t"]\arrow[dd,xshift=-0.75ex,"s"'] & & & \mc{H}_1' \arrow[rr, "\psi_{\mc{H}}"] \arrow[dd,xshift=0.75ex,"t"]\arrow[dd,xshift=-0.75ex,"s"'] & & \mc{H}_1 \arrow[dd,xshift=0.75ex,"t"]\arrow[dd,xshift=-0.75ex,"s"'] \\
			& & & & & & & \\
			M \arrow[rrrrrrr, "\text{Id}", bend right] & & M' \arrow[rrr, "\text{Id}"] \arrow[ll, "f"'] & & & M' \arrow[rr, "f"] & & M
			\end{tikzcd}.
			\end{equation}
			Here, $(\psi_\mc{G},f)\colon (\mc{G}'_1\rightrightarrows M')\rightarrow (\mc{G}_1\rightrightarrows M)$ and 
			$(\psi_\mc{H},f)\colon (\mc{H}'_1\rightrightarrows M')\rightarrow (\mc{H}_1\rightrightarrows M)$ are Morita morphisms of Lie groupoids. Then, by Proposition \ref{Proposition:Moritaequivalentimpliesisomorphicstacks} $B\mc{G}'\cong B\mc{G}$ and $B\mc{H}'\cong B\mc{H}$, and commutativity of \ref{Diagram:MELiegroupoidextensions} gives the following commutative diagram
			\[\begin{tikzcd}
			B\mc{G} \arrow[rrrrrr, bend left,"\Phi"] & & B\mc{G}' \arrow[ll,"\cong"'] \arrow[rr,"\Phi'"] & & B\mc{H}' \arrow[rr,"\cong"] & & B\mc{H}.
			\end{tikzcd}\]
			Reorganizing the above diagram we obtain, 
			\[\begin{tikzcd}
			B\mc{G}' \arrow[rr,"\Phi'"] \arrow[dd,"\cong"'] & & B\mc{H}' \arrow[dd,"\cong"] \\
			& & \\
			B\mc{G} \arrow[rr,"\Phi"]& & B\mc{H}. 
			\end{tikzcd}\] 
			Thus, the gerbe $\Phi\colon B\mc{G}\ra B\mc{H}$ is isomorphic to the gerbe $\Phi'\colon B\mc{G}'\ra B\mc{H}'$. Repeating the same argument for a Morita morphism from the Lie groupoid extension $\phi' \colon\mc{G}'_1\ra \mc{H}'_1\rightrightarrows M'$ to $\phi''\colon \mc{G}''_1\ra \mc{H}''_1\rightrightarrows M''$, we complete the proof.
		\end{enumerate}
	\end{proof}
\end{theorem}
\begin{remark}\label{Remark:PropertyofBGtoBH}
	Let $\mc{D}\ra \mc{C}$ be a gerbe over a stack. Assume further that $\mc{D}\ra \mc{D}\times_{\mc{C}}\mc{D}$ is a representable surjective submersion. In particular, this means there exists  atlases $\underline{X}\ra \mc{C}$ and $\underline{X}\ra \mc{D}$ respectively for the stacks $\mc{C}$ and $\mc{D}$ such that the smooth map $X\times_{\mc{D}}X\ra X\times_\mc{C}X$ is a surjective submersion (Remark \ref{Remark:FullCapacityofDeltaF}). 
	Now we make the following observation:
	
	Let $(\Phi,1_M):(\mc{G}_1\rightrightarrows M)\ra (\mc{H}_1\rightrightarrows M)$ be a Lie groupoid extension and 
	$B\mc{G}\ra B\mc{H}$ be the associated morphism of stacks. Then
	\begin{enumerate}
		\item the morphism of stacks $B\mc{G}\ra B\mc{H}$ is a gerbe over the stack $B\mc{H}$ (Theorem \ref{Theorem:BGBHisaGoS}).
		\item there exists atlases $\underline{M}\ra B\mc{H}$ and $\underline{M}\ra B\mc{G}$ satisfying $\underline{M}\times_{B\mc{G}}\underline{M}=\mc{G}_1$ and $\underline{M}\times_{B\mc{H}}\underline{M}=\mc{H}_1$ (see Example \ref{Example:AtlasforBG}). Moreover, the smooth map $\Phi:\mc{G}_1\ra \mc{H}_1$ associated to the morphism of stacks $\underline{M}\times_{B\mc{G}}\underline{M}\ra \underline{M}\times_{B\mc{H}}\underline{M}$ is a surjective submersion.
	\end{enumerate}
	  \end{remark}
\subsection*{Acknowledgement} S Chatterjee acknowledges research support from SERB, DST, Govt of India grant MTR/2018/000528. P Koushik would like to thank Jochen Heinloth and David Michael Roberts for useful discussions in Mathoverflow, and for responding to queries by e-mails. Authors gratefully acknowledge e-mail communications received from Camille Laurent-Gengoux for the authors' queries. Authors sincerely thank the anonymous referee for making several useful and important suggestions, which helped a lot in improving the presentation and clarity of the article. 

\bibliography{ontwonotionsArxiv}
\bibliographystyle{plain}
\end{document}